\documentclass[11pt]{article}

\usepackage{amsmath,amsthm}

\usepackage{amssymb,latexsym}

\usepackage{enumerate}

\newtheorem{tw}{Theorem}[subsection]
\newtheorem{lm}[tw]{Lemma}
\newtheorem{wn}[tw]{Corollary}
\newtheorem{stw}[tw]{Proposition}
\newenvironment{dow}{\it Proof.\rm}{\hfill $\Box$}

\theoremstyle{definition}
\newtheorem*{df}{Definition}
\newtheorem{uw}[tw]{Remark}

\newcommand{\BL}{{\mathbb L}}
\newcommand{\BN}{{\mathbb N}}
\newcommand{\BR}{{\mathbb R}}

\newcommand{\FF}{{\mathcal{F}}}
\newcommand{\GG}{{\mathcal{G}}}
\newcommand{\HH}{{\mathcal{H}}}

\newcommand{\DM}{{\mathcal{D}}}
\newcommand{\SM}{{\mathcal{S}}}

\newcommand{\MM}{{\mathcal{M}}}

\newcommand{\intt}{{\int_{t}^{T}}}

\newcommand{\into}{{\int_{0}^{T}}}
\newcommand{\intot}{{\int_{0}^{t}}}

\newcommand{\BRD}{{\mathbb{R}^{d}}}

\newcommand{\upe}{{u^{p}_{\varepsilon}}}

\newcommand{\upde}{{u^{p-2}_{\varepsilon}}}
\newcommand{\upce}{{u^{p-4}_{\varepsilon}}}

\newcommand{\esssup}{\mathop{\mathrm{ess\,sup}}}

\newcommand{\nsubsection}{\setcounter{equation}{0}\subsection}

\begin{document}

\title {BSDEs with monotone generator and two irregular
reflecting barriers}
\author {Tomasz Klimsiak \smallskip\\
{\small Faculty of Mathematics and Computer Science,
Nicolaus Copernicus University} \\
{\small  Chopina 12/18, 87--100 Toru\'n, Poland}\\
{\small e-mail: tomas@mat.uni.torun.pl}}
\date{}
 \maketitle
\begin{abstract}
We consider BSDEs with two reflecting  irregular barriers. We give
necessary and sufficient conditions for existence and uniqueness
of $\mathbb{L}^{p}$ solutions for equations with generators
monotone with respect to $y$ and Lipschitz continuous with respect
to $z$, and with data in $\mathbb{L}^{p}$ spaces for $p\ge 1$. We
also prove that the solutions can be approximated via penalization
method.
\end{abstract}


\footnotetext{{\em Key words or phrases:} Reflected backward
stochastic differential equation, monotone generator,
$\mathbb{L}^{p}$-solutions.}

\footnotetext{Research supported by the Polish Minister of Science
and Higher Education under Grant N N201 372 436.}

\nsubsection{Introduction}

Nonlinear reflected BSDEs with one continuous barrier were
introduced in \cite{EKPPQ} as a generalization of the notion of
nonlinear BSDEs considered for the first time in \cite{PP}. At the
same time in \cite{KS} a nonlinear reflected BSDE with two
continuous barriers were introduced. Already in these initial
papers it was shown that reflected BSDEs have natural connections
with the optimal stopping problem, mixed control problem,
viscosity solutions of the obstacle problem for PDEs and Dynkin
games. In subsequent works on reflected BSDEs these connections
were used successfully to investigation of the problems mentioned
above and to the mixed game problem, risk-sensitive optimal
problem, switching problem and other optimality problems (see,
e.g., \cite{H,HL,HZ,J,PengXu1}). The connections with the obstacle
problem for PDEs allowed to give new existence results and
numerical schemes for solutions of PDEs and appeared powerful in
investigation of homogenization problems and regularity properties
of solutions of PDEs (see, e.g., \cite{AA,BHM,HH1,J,Kl2,Kl3,MZ,MX}).

Let $B$ be a standard $d$-dimensional Brownian motion defined on
some complete probability space $(\Omega,\FF,P)$ and let
$\{\FF_t\}$ be the standard augmentation of the filtration
generated by $B$. Suppose we are given two real progressively
measurable processes $U,L$ on $[0,T]$ such that $L\le U$, an
$\FF_T$\,-measurable random variable $\xi$ such that $L_T\le\xi\le
U_T$ and $f:\Omega\times
[0,T]\times\mathbb{R}\times\mathbb{R}^{d}\rightarrow\mathbb{R}$
such that  $f(\cdot,y,z)$ is progressively measurable. Let us
recall that if $U,L$ are continuous then by a solution of the
reflected BSDE with generator $f$, terminal condition $\xi$ and
barriers $L,U$ we mean a triple $(Y,Z,R)$ of progressively
measurable processes such that $t\mapsto
f(\cdot,Y_t,Z_t)\in\mathbb{L}^1(0,T)$,
$P(\int^T_0|Z_t|^2\,dt<\infty)=1$ and $R$ is a continuous  process
finite variation such that $R_{0}=0$ and
\begin{equation}
\label{eq1.1} \left\{
\begin{array}{l} Y_{t}=\xi+\intt
f(s,Y_{s},Z_{s})\,ds
-\intt dR_{s} -\intt Z_{s}\,dB_{s},\quad t\in [0,T],\medskip\\
L_{t}\le Y_{t}\le U_{t},\quad t\in [0,T],\medskip\\ \int_{0}^{T}
(Y_{t}-L_{t})\,dR^{+}_{t}=\int_{0}^{T}(U_{t}-Y_{t})\,dR^{-}_{t}=0,
\end{array}
\right.
\end{equation}
where $R^-,R^+$ are increasing processes such that $R=R^+-R^-$.

Because of many applications of reflected BSDEs many attempts have
been made to relax the assumptions on the data considered in the
pioneering papers \cite{KS,EKPPQ}, i.e. continuity of $U,L$,
linear growth of $f$ and  Lipschitz-continuity of $f$ with respect
to $y,z$, square-integrability of the data. Among the attempts one
can distinguish three main directions. First of all, many papers
are devoted to proving existence of solutions under weaker
assumptions on regularity of the generator. For instance, in
\cite{BHM,EH,HH11,K} generators having quadratic growth with respect to
$z$ are considered. In \cite{LMX,X} it is assumed that $f$ is
continuous and monotone with respect to $y$ and satisfies the
general growth condition, in \cite{LS,C}  it is only required that
$f$ is continuous with respect to $y,z$ and satisfies the linear
growth condition, and in \cite{ZZ} the generator is
left-continuous and monotone with respect to $y$. In the second
group of papers less regular barriers are considered. In
\cite{H,HH} the barriers are c\`adl\`ag whereas in
\cite{PengXu,PengXu1} they are merely mesaurable,
square-integrable and satisfy the so-called Mokobodzki condition
which roughly speaking says that between barriers there exists a
semimartingale having some integrability properties. It is worth
mentioning that in \cite{HH} (see also \cite{HH1}) the Mokobodzki
condition is replaced by the following one: $L_{t}<U_{t},\,
L_{t-}<U_{t-}\,$, $t\in[0,T]$. In the the third group existence
and uniqueness results for reflecting BSDEs with data that are not
square-integrable are proved (see \cite{Kl4,S} for results for
equations with data in $\mathbb{L}^{p}$ with $p\in [1,2)$ and
\cite{A,HP} for the case where $p\in (1,2)$). Finally, let us note
that to our knowledge at present there are only few papers, all on
equations with one reflecting barrier, that cover two of the three
cases described above (see \cite{A,Kl4,S}) and there is no paper
which covers all the three.

In the present paper we consider reflected BSDEs with data in
$\mathbb{L}^{p}$ spaces, $p\in [1,2)$, and with two merely
measurable barriers satisfying the generalized Mokobodzki
condition. Regarding the generator, we focus attention on its
dependence on the variable $y$. We assume that $f$ is monotone
with respect to $y$,  Lipschitz continuous with respect to $z$ and
satisfies a very general growth condition with respect to $y$
which is weaker then the so-called general growth condition
considered in \cite{BDHPS} in case of nonreflected BSDEs and in
\cite{Kl4,S} in case of BSDEs with one continuous reflecting
barrier. This growth condition has the form
\[
\forall_{r>0}\quad  \sup_{|y|\le r}|f(\cdot,y,0)-f(\cdot,0,0)| \in
\mathbb{L}^{1}(0,T).
\]
It seems to be the best possible growth condition on $f$ with
respect to $y$; it is widely used in the theory of  partial
differential equations (see \cite{Betal.} and the references given
there). Motivated by future applications to PDEs, we consider
reflected BSDEs more general than (\ref{eq1.1}). Suppose that in
addition to $\xi,f, U,L$ we are given a progressively measurable
c\`adl\`ag process $V$ such that $V_0=0$. The problem is to find a
triple $(Y,Z,R)$ of progressively measurable processes such that
$Z,f(\cdot,Y,Z)$ have the same integrability properties as in case
of equation (\ref{eq1.1}) and $R$ is a c\`adl\`ag process of
finite variation such that $R_{0}=0$ and
\begin{equation}
\label{eq1.2} \left\{
\begin{array}{l} Y_{t}=\xi+\intt
f(s,Y_{s},Z_{s})\,ds+\intt dV_{s}+\intt dR_{s}
-\intt Z_{s}\,dB_{s},\quad t\in [0,T],\medskip\\
L_{t}\le Y_{t}\le U_{t}\quad \mbox{for a.e. }t\in [0,T],\medskip\\
\int_{0}^{T} (Y_{t-}-\hat{L}_{t-})\,dR^{+}_{t}
=\int_{0}^{T}(\check{U}_{t-}-Y_{t-})\,dR^{-}_{t}=0
\end{array}
\right.
\end{equation}
for any progressively measurable c\`adl\`ag processes
$\hat{L},\check{U}$ such that $L_{t}\le \hat{L}_{t}\le Y_{t}\le
\check{U}_{t}\le U_{t}$ for a.e. $t\in [0,T]$.
%The name ``RBSDE
%with measure energy" follows from the fact that we would like to
%apply this kind of stochastic equations to obstacle problem for
%PDEs with measure data. In this framework a finite variation
%process $V$ appearing in RBSDEs will be in Revuz duality with some
%soft (smooth) measure on the right hand side of partial
%differential equation.

In the paper we  give  existence, uniqueness and comparison
results for equations of the form (\ref{eq1.2}). We also provide
Snell's envelope representation and prove that for every $p\in
[1,2)$ solutions of (\ref{eq1.2}) can be approximated by the
penalization method.

In \cite{PengXu} it is shown that in case $p=2$  there exists a
solution of BSDE with two reflecting barriers for Lipschitz
continuous  generators satisfying the linear growth
condition. The remarkable feature of the theory of $\mathbb{L}^p$
solutions of equations with monotone generators is the fact that
conditions ensuring existence of solutions of nonreflecting BSDEs
together with the Mokobodzki condition are insufficient for
existence of solutions of reflected BSDEs in the given class of
integrability. To get existence we introduce a generalized
Mokobodzki condition which contrary to the standard one also
relates the growth of the generator with that of the barriers. In
fact, we formulate two versions of the generalized condition: one
for $p>1$ and another one for $p=1$. One of our main results says
that under some minimal assumptions ensuring existence of
solutions of nonreflected BSDE (see \cite[Theorem 4.2, Theorem
6.3]{BDHPS}) the generalized Mokobodzki condition is necessary and
sufficient for existence of a solution of reflected BSDE in a
given class of integrability. In the proof of that result we use
among others things the method of supersolutions applied earlier
in \cite{PengXu} and the technique of stopping times used in
\cite{HH1,Kl4,S}.

In the last part of the paper we answer the question what happens
if despite the fact that we consider reflected BSDEs with monotone
generator we only assume the standard Mokobodzki condition. From
the comments given in the preceding paragraph it follows that in
that case in general we can not expect existence of
$\mathbb{L}^{p}$ solutions. Nevertheless, we show that there
always exists a unique solution of (\ref{eq1.2}). It may happen,
however, that some of its components are  nonintegrable  for every
$p>0$.

The paper is organized as follows. In Section \ref{sec2} we
provide basic notation used in the paper and we prove It\^o's
formula for c\`adl\`ag processes and the function $x\rightarrow
|x|^{p},\, p\in [1,2)$. In Section \ref{sec3} we prove existence,
comparison results a priori estimates  for solutions of
(\ref{eq1.2}). In Section \ref{sec4} we prove some properties of
supersolutions of (\ref{eq1.2}) and some useful lemmas required in
proofs of existence of solutions of reflected BSDEs. In
particular, we prove a generalization of the monotone convergence
theorem for BSDEs,  Snell's envelope representation of
supersolutions and a lemma on behavior of their jumps. In Section
\ref{sec5} we prove existence and uniqueness results for BSDEs
with one reflecting barrier whereas in Section \ref{sec6} for equations
with two barriers. In Section \ref{sec7} we consider the case of
nonintegrable solutions.

\nsubsection{Preliminary results}
\label{sec2}

Let $B=\{B_{t}, t\ge 0\}$ be a standard $d$-dimensional Brownian
motion defined on some complete filtered probability space
$(\Omega,\FF,P)$, where $\FF=\{\FF_{t}, t\ge 0\}$ is the augmented
filtration generated by $B$. In the whole paper all notions whose
definitions are related to some filtration are understood with
respect to the filtration $\FF$.

Given a stochastic process $X$ on $[0,T]$ with values in
$\mathbb{R}^{n}$ we set $X^{*}_{t}=\sup_{0\le s\le t}|X_{s}|$,
$t\in[0,T]$, where $|\cdot|$ denotes the Euclidean norm on
$\mathbb{R}^{n}$.  By $\mathcal{S}$ (resp. $\DM$) we denote the
set of all progressively measurable continuous (resp. c\'adl\`ag)
processes. For $p>0$ we denote by $\mathcal{S}^{p}$ (resp.
$\DM^{p}$) the set of all processes $X\in\mathcal{S}$ (resp.
$\DM$) such that
\[
E\sup_{t\in[0,T]}|X_{t}|^{p}<\infty.
\]
$M$ is the set of all progressively measurable
processes $X$ such that
\[
P(\int_{0}^{T}|X_{t}|^{2}\,dt<\infty)=1
\]
and for $p>0$, $M^{p}$ is the set of all processes $X\in M$ such
that
\[
E(\int_{0}^{T}|X_{t}|^{2}\,dt)^{p/2}<\infty.
\]
For $p,q>0$, $\mathbb{L}^{p,q}(\FF)$ denotes the set of all
progressively measurable processes $X$ such that
\[
E(\int_{0}^{T}|X_{t}|^{p}\,dt)^{q/(1\wedge 1/p)}<\infty.
\]
For brevity we denote $\mathbb{L}^{p,p}(\FF)$ by
$\mathbb{L}^{p}(\FF)$.

For a given measurable space $(X,\mu,\GG)$ we denote  by
$\mathbb{L}^{p}(X,\mu,\GG)$ the set of all $\GG$-measurable real
valued functions $f$ such that
$\int_{X}|f|^{p}(x)\,d\mu(x)<\infty.$ If it is clear from the
context which measure and $\sigma$-field is considered we omit
them in the notation.

$\MM_{c}$ (resp. $\MM^{loc}_{c}$) is the set of all continuous martingales (resp. local martingales) and
$\MM^{p}_{c}$, $p\ge1$, is the set of all martingales
$M\in\MM_{c}$ such that $E(\langle M \rangle_{T})^{p/2}<\infty$.
$\mathcal{V}_{c}$ (resp. $\mathcal{V}^{+}_{c}$) is the set of all
continuous progressively measurable processes of finite variation
(resp. increasing processes) such that $V_{0}=0$ and
$\mathcal{V}^{p}_{c}$ (resp. $\mathcal{V}^{+,p}_{c}$) is the set
of all processes $V\in\mathcal{V}_{c}$ (resp.
$V\in\mathcal{V}^{+}_{c}$) such that $E|V|^{p}_{T}<\infty$.
$\mathcal{V}$ (resp. $\mathcal{V}^{+}$) is the set of all
c\`adl\`ag progressively measurable processes of finite variation
(resp. increasing processes) such that $V_{0}=0$ and
$\mathcal{V}^{p}$ (resp. $\mathcal{V}^{+,p}$) is the set of all
processes $V\in\mathcal{V}$ (resp. $V\in\mathcal{V}^{+}$) such
that $E|V|^{p}_{T}<\infty$. $\HH^{p}=\MM_{c}^{p}+\mathcal{V}^{p}$,
$\HH^{p}_{c}=\MM_{c}^{p}+\mathcal{V}_{c}^{p}$. For a given process
$V\in\mathcal{V}$ by $dV$ we denote the random measure generated
by its trajectories.

By $\mathcal{T}$ we denote the set of all stopping times with
values in  $[0,T]$ and by $\mathcal{T}_t$ the set of all stopping
times with values in $[t,T]$. We say that a sequence
$\{\tau_{k}\}\subset \mathcal{T}$ is stationary if
\[
P(\liminf_{k\rightarrow+\infty} \{\tau_{k}=T\})=1.
\]
For a given measurable process $Y$ of class (D) we denote
\[
\|Y\|_{1}=\sup\{E|Y_{\tau}|,\tau\in\mathcal{T}\}.
\]
%We also put $f_{t}=|f(t,0,0)|,\, t\in [0,T]$.

For $X\in\DM$ we set $X_{t-}=\lim_{s\nearrow t} X_{s}$ and $\Delta
X_{t}= X_{t}-X_{t-}$ with the convention that $X_{0-}=0$. Let
$\{X^{n}\}\subset \DM$, $X\in\DM$. We say that $X^{n}\rightarrow
X$ in ucp if $\sup_{t\in [0,T]}|X^{n}_{t}-X_{t}|\rightarrow 0$ in
probability $P$.

In the whole paper all equalities and inequalities and other
relations between random elements are understood to hold $P$-a.s..
To avoid ambiguity we stress that writing $X_{t}=Y_{t}$, $t\in
[0,T]$ we mean that $X_{t}=Y_{t}$, $t\in [0,T]$, $P$-a.s., while
writing $X_{t}=Y_{t}$ for a.e. (resp. for every) $t\in [0,T]$ we
mean that $X_{t}=Y_{t}, P$-a.s. for a.e. (resp. for every) $t\in
[0,T]$. We also adopt the convention that
$\int_{a}^{b}=\int_{(a,b]}$.

%$Prog$ denotes the $\sigma$-field of progressive subsets of
%$[0,T]\times\Omega$.

$T_{k}(x)=\min\{k, \max\{-k,x\}\}$, $x\in\BR$,
$x^{+}=\max\{x,0\}$, $x^{-}=\max\{-x,0\}$ and
\[
\hat x=\hat{\mbox{sgn}}(x),\quad
\hat{\mbox{sgn}}(x)=\mathbf{1}_{x\neq 0}\frac{x}{|x|},\quad
x\in\mathbb{R}^d.
\]

One of our basic tools will be the following It\^o-Tanaka formula.
The formula is probably known, but we do not have a reference.

\begin{stw}
\label{prop2.1} Let $p\ge 1$ and let $X$ be a progressively
measurable process of the form
\begin{align}\label{eq2.1}
X_{t}=X_{0}+\intot dK_{s}+\intot Z_{s}\,dB_{s},\quad t\in [0,T],
\end{align}
where $K\in \mathcal{V}$ and $Z\in M$. Then
\begin{align}
\label{eq2.02} |X_{t}|^{p}-|X_{0}|^{p}&=p\intot
|X_{s-}|^{p-1}\hat{X}_{s-}\,dK_{s}
+p\intot|X_{s}|^{p-1}\hat{X}_{s}\,dB_{s}\nonumber\\
&\quad+\frac12 p(p-1)\mathbf{1}_{X_{s}\neq 0}|X_{s}|^{p-2}|Z_{s}|^{2}\,ds
+L_{t}\mathbf{1}_{\{p=1\}}+I_{t}(p),
\end{align}
where $L\in \mathcal{V}^{+}_{c}$ and
\[
I_{t}(p)=\sum_{0<s\le t}\Delta|X_{s}|^{p}
-\sum_{0<s\le t}p|X_{s-}|^{p-1}\hat{X}_{s-}\Delta X_{s}
\]
is a nonnegative increasing pure jump process.
\end{stw}
\begin{dow}
Write $u_{\varepsilon}^{p}(x)=(|x|^{2}+\varepsilon^{2})^{p/2}$,
$x\in\mathbb{R}$. Then
\[
(\frac{d}{dx}u^{p}_{\varepsilon})(x)=pu^{p-2}(x)x, \quad
(\frac{d^{2}}{dx^{2}}u^{p}_{\varepsilon})(x)=pu^{p-2}(x)
+p(p-2)u^{p-4}_{\varepsilon}(x)x^{2}
\]
for $x\in\mathbb{R}$. By It\^o's formula,
\begin{align}\label{eq2.2}
\nonumber\upe(X_{t})-\upe(X_{0})&=\intot\frac{d}{dx}\upe(X_{s-})\,dX_{s}
+\frac12\intot\frac{d^{2}}{dx^{2}} \upe(X_{s})d[X,X]^{c}_{s}\\&
\quad\nonumber+\sum_{0<s\le t} \{\Delta \upe(X_{s})
-\frac{d}{dx}\upe(X_{s-})\Delta X_{s}\}\\&\nonumber
= \intot p\upde(X_{s-})X_{s-}\,dK_{s}
+\intot p\upde(X_{s-})X_{s-}Z_{s}\,dB_{s}\\
&\quad +\frac12 \intot p\upde(X_{s})+p(p-2)
\upce(X_{s})X^{2}_{s}Z^{2}_{s}\,ds\nonumber\\
&\quad+ \sum_{0<s\le t} \{\Delta \upe(X_{s})-p\upde(X_{s-})
X_{s-}\Delta X_{s}\}.
\end{align}
Observe that $\upe\rightarrow |\cdot |^{p}$ uniformly on compact
subsets of $\mathbb{R}$. Hence
\begin{equation}
\label{eq2.3} \upe(X)-\upe(X_{0})\rightarrow |X|^{p}-|X_{0}|^{p}
\quad\mbox{in ucp}.
\end{equation}
By elementary computation, for $\varepsilon\le 1$ we have
\begin{equation}
\label{eq2.03} |\Delta\upe(X_{t})|\le
|\frac{d\upe}{dx}(X^{*}_{t})||\Delta X_{t}| \le
(|X^{*}_{t}|^{2}+1)^{p/2}|\Delta X_{t}|.
\end{equation}
Since $p\upde(x)x\rightarrow p|x|^{p-1}\hat{x}$ for  $x\in\BR$
and, by (\ref{eq2.1}), $\sum_{0< s\le t}|\Delta X_{s}|\le
|K|_{t}$, applying the Lebesgue dominated convergence theorem we
see that as $\varepsilon\rightarrow 0^{+}$ then  $P$-a.s.,
\begin{equation}
\label{eq2.4}
\sum_{0<s\le t} \{\Delta \upe(X_{s})-p\upde(X_{s-}) X_{s-}\Delta
X_{s}\}\rightarrow \sum_{0<s\le t}
\{\Delta|X_{s}|^{p}-p|X_{s-}|^{p-1}\hat{X}_{s-}\Delta X_{s}\}
\end{equation}
for $t\in [0,T]$. Using  once again the second inequality in
(\ref{eq2.03}) and the Lebesgue dominated convergence theorem we
conclude that
\begin{equation}
\label{eq2.5} \intot p\upde(X_{s-})X_{s-}\,dK_{s}
\rightarrow\intot p|X_{s-}|^{p-1}\hat{X}_{s-}\,dK_{s},\quad t\in
[0,T]
\end{equation}
$P$-a.s. and
\begin{equation}
\label{eq2.6} \int_{0}^{\cdot} p\upde(X_{s})X_{s} Z_{s}\,dB_{s}
\rightarrow \int_{0}^{\cdot} p|X_{s}|\hat{X}_{s}
Z_{s}\,dB_{s}\quad\mbox{in ucp}
\end{equation}
as $\varepsilon\rightarrow 0^{+}$. For every $q\in\mathbb{R}$ and
$x\in\BR$,
\begin{align*}
u^{q}_{\varepsilon}(x)|x|^{2}&=(|x|^{2}+\varepsilon^{2})^{q/2}|x|^{2}
=(|x|^{2}+\varepsilon^{2})^{q/2}(|x|^{2}
+\varepsilon^{2})-\varepsilon^{2}u_{\varepsilon}^{q}(x)\\
&=u_{\varepsilon}^{q+2}(x)-\varepsilon^{2}u^{q}_{\varepsilon}(x).
\end{align*}
Hence
\begin{align*}
S_{t}(\varepsilon)&\equiv \frac12\intot
(p\upde(X_{s})+p(p-2)\upce(X_{s})X^{2}_{s}Z^{2}_{s})\,ds\\
&= \frac12\intot p(\upce(X_{s})|X_{s}|^{2}
+\varepsilon^{2}\upce(X_{s}))|Z_{s}|^{2}\,ds\\
&\quad+\frac12\intot
p(p-2)\upce(X_{s})|X_{s}|^{2}|Z_{s}|^{2}\,ds\\
&=\frac12 \intot p(p-1)\upce(X_{s})|X_{s}|^{2}|Z_{s}|^{2}\,ds
+\frac12\intot p\varepsilon^{2}\upce(X_{s})|Z_{s}|^{2}\,ds,
\end{align*}
that is
\begin{equation}
\label{eq2.8} S_{t}(\varepsilon)=\frac12\intot
p(p-1)\upce(X_{s})|X_{s}|^{2}|Z_{s}|^{2}\,ds
+L_{t}^{\varepsilon}(p),\quad t\in [0,T],
\end{equation}
where $L_{t}^{\varepsilon}(p)\equiv\frac12\intot
p\varepsilon^{2}\upce(X_{s})|Z_{s}|^{2}\,ds$. Since
$\frac{|X_{s}|}{u_{\varepsilon}(X_{s})}\nearrow
\mathbf{1}_{\{X_{s}\neq 0\}}$, applying the  monotone convergence
theorem gives
\begin{align}
\label{eq2.9} \int_{0}^{\cdot}
\upce(X_{s})|X_{s}|^{2}|Z_{s}|^{2}\,ds&=\int_{0}^{\cdot} \left(
\frac{|X_{s}|}{u_{\varepsilon}(X_{s})}\right)^{4-p}
|X_{s}|^{p-2}|Z_{s}|^{2}\mathbf{1}_{\{X_{s}\neq 0\}}\,ds\nonumber\\
&\quad\rightarrow \int_{0}^{\cdot} \mathbf{1}_{\{X_{s}\neq 0\}}
|X_{s}|^{p-2}|Z_{s}|^{2}\,ds\quad\mbox{in ucp}.
\end{align}
Putting (\ref{eq2.2})--(\ref{eq2.9}) together we get
\begin{align*}
|X_{t}|^{p}-|X_{0}|^{p}&=p\intot |X_{s-}|^{p-1}\hat{X}_{s-}\,dK_{s}
+p\intot|X_{s}|^{p-1}\hat{X}_{s}\,dB_{s}\\
&\quad+\frac12 p(p-1)\mathbf{1}_{X_{s}\neq
0}|X_{s}|^{p-2}|Z_{s}|^{2}\,ds +L_{t}(p)+I_{t}(p)
\end{align*}
with $L(p)\in\mathcal{V}^{+}$ defined by $L_{t}(p)=
\lim_{\varepsilon\rightarrow 0^{+}}L_{t}^{\varepsilon}(p)$, $t\in
[0,T]$. An elementary computation analogous to that in the proof
of \cite[Lemma 2.2]{BDHPS} shows that in fact
$L_{t}(p)=L_{t}(p)\mathbf{1}_{\{p=1\}}$. Therefore putting
$L_{t}=L_{t}(1)$, $t\in[0,T]$ we get (\ref{eq2.02}). Finally,
comparing It\^o's formula proved in \cite[Section 2]{Protter} with
formula (\ref{eq2.02}) with $p=1$ shows that $L$ is a symmetric
local time at zero of the process $X$. In particular, $L$ is
continuous, and the proof is complete.
\end{dow}

\nsubsection{Existence and uniqueness of solutions of BSDEs}
\label{sec3}

In this section we study existence and uniqueness of solutions of
equations of the form
\begin{equation}
\label{eq3.1}
Y_{t}=\xi+\int_{t}^{T}f(s,Y_{s},Z_{s})\,ds
+\intt dV_{s}-\int_{t}^{T} Z_{s}\,dB_{s},\quad t\in[0,T],
\end{equation}
where $V\in\mathcal{V}$. In what
follows it will be convenient to denote equations of the form
(\ref{eq3.1}) by BSDE$(\xi,f+dV)$.

With formula (\ref{eq2.02}) at hand, to prove existence and
uniqueness of solutions of (\ref{eq3.1}) it suffices to repeat
step by step, with some obvious changes, the proofs of
corresponding results from \cite{BDHPS} obtained for
multidimensional equations of the form (\ref{eq3.1}) with $V=0$.
In \cite{BDHPS}, however, to prove existence of solutions of
(\ref{eq3.1}) the authors use some important results from other
papers proved in case $V=0$. Instead of repeating arguments from
all these papers and from \cite{BDHPS} we decided to take
advantage of the fact that we are concerned with one dimensional
equations and provide new proof which is based only on results (or
slightly modified results) obtained in \cite{BDHPS}.

Let $p\ge 1$. We will need the following hypotheses.

\begin{enumerate}
\item[(H1)] $E|\xi|^{p}+E(\int_{0}^{T}|f(t,0,0)|\,dt)^{p}
+E(\into d|V|_{s})^{p}<\infty$.
\item[(H2)]There exists $\lambda>0$ such that $|f(t,y,z)-f(t,y,z')|
\le \lambda |z-z'|$
for every $t\in [0,T], y\in\mathbb{R}, z,z'\in\BRD$.
\item[(H3)]There exists $\mu\in\mathbb{R}$ such that
$(f(t,y,z)-f(t,y',z))(y-y')\le \mu(y-y')^{2}$
for every $t\in [0,T], y\in\mathbb{R}, z,z'\in\BRD$.
\item [(H4)]For every $(t,z)\in[0,T]\times\BRD$ the mapping
$\mathbb{R}\ni y\rightarrow f(t,y,z)$
is continuous.
\item[(H5)]For every $r>0$ the mapping
$[0,T]\ni t\rightarrow\sup_{|y|\le r}|f(t,y,0)-f(t,0,0)|$ belongs
to  $\mathbb{L}^{1}(0,T)$.
\item[(A)]There exist $\mu\in\mathbb{R}$ and $\lambda \ge 0$ such
that
\[
\hat{y}f(t,y,z)\le f_{t}+\mu|y|+\lambda |z| ,
\]
for every $(t,y,z)\in [0,T]\times\mathbb{R}\times\BRD$, where
$\hat{y}=\mathbf{1}_{\{y\neq 0\}}\frac{y}{|y|}$ and $f_{t}$ is a
nonnegative progressively measurable process.
\item[(Z)]There exist $\alpha\in (0,1)$, $\gamma\ge 0$ and
nonnegative process $g\in\mathbb{L}^{1}(\FF)$ such that
\[
|f(t,y,z)-f(t,y,0)|\le \gamma (g_{t}+|y|+|z|)^{\alpha}
\]
for every $(t,y,z)\in [0,T]\times\mathbb{R}\times\BRD$.
\end{enumerate}

We begin with proving comparison result for BSDEs of the form (\ref{eq3.1}).

\begin{stw}\label{prop3.4}
Let $(Y^{i},Z^{i})$, $i=1,2$, be a solution of
BSDE$(\xi^{i},f^{i}+dV^{i})$. Assume that
$(Y^{1}-Y^{2})^{+}\in\mathcal{D}^{q}$ for some $q>1$. If
$\xi^{1}\le\xi^{2}$, $dV^{1}\le dV^{2}$, $f^1,f^2$ satisfy
\mbox{\rm(H3)} and either
\begin{equation}
\label{eq3.11} f^{2}\mbox{ satisfies {\rm(H2)}},\quad
\mathbf{1}_{\{Y^{1}_{t}>Y^{2}_{t}\}}
(f^{1}(t,Y^{1}_{t},Z^{1}_{t})-f^{2}(t,Y^{1}_{t},Z^{1}_{t}))\le 0
\end{equation}
for a.e. $t\in [0,T]$ or
\begin{equation}
\label{eq3.12} f^{1}\mbox{satisfies {\rm(H2)}},\quad
\mathbf{1}_{\{Y^{1}_{t}>Y^{2}_{t}\}}
(f^{1}(t,Y^{2}_{t},Z^{2}_{t})-f^{2}(t,Y^{2}_{t},Z^{2}_{t}))\le 0
\end{equation}
for a.e. $t\in [0,T]$ then $Y^{1}_{t}\le Y^{2}_{t}$, $t\in [0,T]$.
\end{stw}
\begin{dow}
Without loss of generality we may assume that $\mu\le 0$. Assume
that (\ref{eq3.11}) is satisfied and let $p\in(1,q)$. Then by the
It\^o-Tanaka formula (see \cite[Chapter IV, Section 7]{Protter})
and Proposition \ref{prop2.1}, for every stopping time
$\tau\in\mathcal{T}$ we have
\begin{align}\label{eq3.13}
&|(Y^{1}_{t\wedge \tau}-Y^{2}_{t\wedge \tau})^{+}|^{p}
+\frac{p(p-1)}{2}\int_{t\wedge\tau}^{\tau}\mathbf{1}_{\{Y^{1}_{s}
\neq Y^{2}_{s}\}}|(Y^{1}_{s}-Y^{2}_{s})^{+}|^{p-2}
|Z^{1}_{s}-Z^{2}_{s}|^{2}\,ds \nonumber\\
&\quad\le|(Y^{1}_{\tau}-Y^{2}_{\tau})^{+}|^{p} +p\int_{t\wedge
\tau}^{\tau}|(Y^{1}_{s}-Y^{2}_{s})^{+}|^{p-1}
(f^{1}(s,Y^{1}_{s},Z^{1}_{s})-f^{2}(s,Y^{2}_{s},Z^{2}_{s}))\,ds
\nonumber\\
&\qquad+p\int_{t\wedge\tau}^{\tau}
|(Y^{1}_{s-}-Y^{2}_{s-})^{+}|^{p-1 }\,(dV^{1}_{s}-dV^{2}_{s})\nonumber\\
&\qquad -p\int_{t\wedge\tau}^{\tau}
|(Y^{1}_{s}-Y^{2}_{s})^{+}|^{p-1} (Z^{1}_{s}-Z^{2}_{s})\,dB_{s}.
\end{align}
By the assumptions,
\begin{align*}
&\int_{t\wedge\tau}^{\tau}|(Y^{1}_{s-}-Y^{2}_{s-})^{+}|^{p-1}\,
(dV^{1}_{s}-dV^{2}_{s})\le 0.
\end{align*}
Moreover, using (\ref{eq3.11}), (H2), (H3) one can check that for
a.e. $t\in [0,T]$,
\begin{align*}
&\mathbf{1}_{\{Y^{1}_{t}>Y^{2}_{t}\}}
(f^{1}(t,Y^{1}_{t},Z^{1}_{t})-f^{2}(t,Y^{2}_{t},Z^{2}_{t}))
\le\lambda\mathbf{1}_{\{Y^{1}_{t}>Y^{2}_{t}\}}
|Z^{1}_{t}-Z^{2}_{t}| .
\end{align*}
Therefore from (\ref{eq3.13}) it follows that for a.e. $t\in
[0,T]$,
\begin{align}
\label{eq3.07}
&|(Y^{1}_{t\wedge\tau}-Y^{2}_{t\wedge\tau})^{+}|^{p}
+\frac{p(p-1)}{2} \int_{t\wedge\tau}^{\tau}\mathbf{1}_{\{Y^{1}_{s}
\neq Y^{2}_{s}\}}|(Y^{1}_{s}-Y^{2}_{s})^{+}|^{p-2}
|Z^{1}_{s}-Z^{2}_{s}|^{2}\,ds\nonumber\\
&\qquad\le |(Y^{1}_{\tau}-Y^{2}_{\tau})^{+}|^{p}
+p\lambda\int_{t\wedge\tau}^{\tau}|(Y^{1}_{s}
-Y^{2}_{s})^{+}|^{p-1}|Z^{1}_{s}-Z^{2}_{s}|\,ds\nonumber\\
&\qquad\quad-p\int_{t\wedge\tau}^{\tau}
|(Y^{1}_{s}-Y^{2}_{s})^{+}|^{p-1}
(Z^{1}_{s}-Z^{2}_{s})\,dB_{s}\nonumber\\
&\qquad\le|(Y^{1}_{\tau}-Y^{2}_{\tau})^{+}|^{p}
+\frac{p\lambda^{2}}{p-1}\int_{t\wedge\tau}^{\tau}
|(Y^{1}_{s}-Y^{2}_{s})^{+}|^{p}\,ds\nonumber\\
&\qquad\quad+\frac{p(p-1)}{4}\int_{t\wedge\tau}^{\tau}
|(Y^{1}_{s}-Y^{2}_{s})^{+}|^{p-2}
|Z^{1}_{s}-Z^{2}_{s}|^{2}\,ds\nonumber\\
&\qquad\quad-p\int_{t\wedge\tau}^{\tau}
|(Y^{1}_{s}-Y^{2}_{s})^{+}|^{p-1}(Z^{1}_{s}-Z^{2}_{s})\,dB_{s}.
\end{align}
Set $\tau_{k}=\inf\{t\in [0,T];\int_{0}^{t}
|(Y^{1}_{s}-Y^{2}_{s})^{+}|^{2(p-1)}
|Z^{1}_{s}-Z^{2}_{s}|^{2}\,ds\ge k\}\wedge T$. Then from
(\ref{eq3.07}) with $\tau=\tau_{k}$ we obtain
\[
E|(Y^{1}_{t}-Y^{2}_{t})^{+}|^{p} \le
E|(Y^{1}_{\tau_{k}}-Y^{2}_{\tau_{k}})^{+}|^{p}
+\frac{p\lambda^{2}}{p-1}E\int_{t\wedge\tau_{k}}^{\tau_{k}}
|(Y^{1}_{s}-Y^{2}_{s})^{+}|^{p}\,ds.
\]
Since $(Y^{1}-Y^{2})^{+}\in \mathcal{D}^{q}$ for some $q>1$, the
process $(Y^{1}-Y^{2})^{+}$ is of class (D). Therefore letting
$k\rightarrow+\infty$ in the above inequality and using the fact
that $\xi^{1}\le\xi^{2}$ we get
\[
E|(Y^{1}_{t}-Y^{2}_{t})^{+}|^{p} \le
\frac{p\lambda^{2}}{p-1}E\int_{t}^{T}
|(Y^{1}_{s}-Y^{2}_{s})^{+}|^{p}\,ds,\quad t\in [0,T],
\]
so the desired result follows by Gronwall's lemma. Since the proof
in case (\ref{eq3.12}) is satisfied is analogous, we omit it.
\end{dow}

\begin{wn}
\label{cor3.10} Assume \mbox{\rm(Z)}. Let $(Y^{i},Z^{i})$,
$i=1,2$, be a solution of BSDE$(\xi^{i},f^{i}+dV^{i})$ such that
$(Y^{i},Z^{i})\in \mathbb{L}^{q}(\FF)\otimes\mathbb{L}^{q}(\FF)$
for some $q>\alpha$. If $(Y^{1}-Y^{2})^{+}$ is of class
\mbox{\rm(D)}, $\xi^{1}\le\xi^{2}$, $dV^{1}\le dV^{2}$, $f^1,f^2$
satisfy \mbox{\rm(H2)} and \mbox{\rm(\ref{eq3.11})} or
\mbox{\rm(\ref{eq3.12})} is satisfied then $Y^{1}_{t}\le
Y^{2}_{t}$, $t\in [0,T]$.
\end{wn}
\begin{dow}
We only consider the case where (\ref{eq3.11}) is satisfied. As
usual, without loss of generality we may assume that $\mu\le 0$.
Due to Proposition \ref{prop3.4} it suffices to show that
$(Y^{1}-Y^{2})^{+}\in\mathcal{D}^{p}$ for some $p>1$. By the
It\^o-Tanaka  formula and the inequality $dV^{1}\le dV^{2}$, for
every stopping time $\tau\in\mathcal{T}$,
\begin{align}
\label{eq3.14}
\nonumber(Y^{1}_{t\wedge\tau}-Y^{2}_{t\wedge\tau})^{+}
&\le (Y^{1}_{\tau}-Y^{2}_{\tau})^{+}
+\int_{t\wedge\tau}^{\tau} \mathbf{1}_{\{Y^{1}_{s} >
Y^{2}_{s}\}}(f^{1}(s,Y^{1}_{s},Z^{1}_{s})
-f^{2}(s,Y^{2}_{s},Z^{2}_{s}))\,ds\\
&\quad-\int_{t\wedge\tau}^{\tau} \mathbf{1}_{\{Y^{1}_{s} >
Y^{2}_{s}\}}(Z^{1}_{s}-Z^{2}_{s})\,dB_{s}.
\end{align}
Write
\begin{equation}
\label{eq3.09} I_{t}\equiv \mathbf{1}_{\{Y^{1}_{t}
> Y^{2}_{t}\}}(f^{1}(t,Y^{1}_{t},Z^{1}_{t})
-f^{2}(t,Y^{2}_{t},Z^{2}_{t})).
\end{equation}
By (\ref{eq3.11}),
\begin{align*}
I_{t}&= \mathbf{1}_{\{Y^{1}_{t}> Y^{2}_{t}\}}
\{f^{1}(t,Y^{1}_{t},Z^{1}_{t})
-f^{2}(t,Y^{1}_{t},Z^{1}_{t})+f^{2}(t,Y^{1}_{t},Z^{1}_{t})
-f^{2}(t,Y^{2}_{t},Z^{2}_{t})\}\\
&\le \mathbf{1}_{\{Y^{1}_{t}> Y^{2}_{t}\}}
(f^{2}(t,Y^{1}_{t},Z^{1}_{t}) -f^{2}(t,Y^{2}_{t},Z^{2}_{t})),
\end{align*}
and by monotonicity of $f^{2}$ with respect to $y$,
\begin{align*}
& \mathbf{1}_{\{Y^{1}_{t}> Y^{2}_{t}\}}
(f^{2}(t,Y^{1}_{t},Z^{1}_{t}) -f^{2}(s,Y^{2}_{t},Z^{2}_{t}))
=\mathbf{1}_{\{Y^{1}_{t}> Y^{2}_{t}\}}
(f^{2}(t,Y^{1}_{t},Z^{1}_{t})
-f^{2}(t,Y^{1}_{t},0))\\
&\quad+\mathbf{1}_{\{Y^{1}_{t}>Y^{2}_{t}\}}(f^{2}(t,Y^{1}_{t},0)
-f^{2}(t,Y^{2}_{t},0))+\mathbf{1}_{\{Y^{1}_{t}> Y^{2}_{t}\}}
(f^{2}(t,Y^{2}_{t},0)-f^{2}(s,Y^{2}_{t},Z^{2}_{t}))\\
&\le \mathbf{1}_{\{Y^{1}_{t}> Y^{2}_{t}\}}
(f^{2}(t,Y^{1}_{t},Z^{1}_{t})
-f^{2}(t,Y^{1}_{t},0))+\mathbf{1}_{\{Y^{1}_{t}> Y^{2}_{t}\}}
(f^{2}(t,Y^{2}_{t},0)-f^{2}(t,Y^{2}_{t},Z^{2}_{t})).
\end{align*}
Using assumption (Z) we conclude from the above inequality that
\[
I_{t}\le2\gamma(g_{t}+|Y^{1}_{t}|+|Z_{t}^{1}|
+|Y_{t}^{2}|+|Z_{t}^{2}|)^{\alpha}
\]
for a.e. $t\in[0,T]$. Therefore taking the conditional expectation with
respect to $\FF_t$ of both sides of (\ref{eq3.14}) with $\tau$
replaced by $\tau_k=\inf\{t\in [0,T];\int_{0}^{t}
|Z^{1}_{s}-Z^{2}_{s}|^{2}\,ds\ge k\}\wedge T$,  letting
$k\rightarrow +\infty$ and using the fact that $(Y^{1}-Y^{2})^{+}$
is of class (D) we obtain
\[
(Y^{1}_{t}-Y^{2}_{t})^{+} \le 2\gamma
E^{\mathcal{F}_{t}}\{\int_{0}^{T}
(g_{t}+|Y^{1}_{t}|+|Z_{s}^{1}|+|Y^{2}_t|+|Z^{2}_t|)^{\alpha}\,dt\}.
\]
Using now the assumptions of the corollary and applying Doob's
inequality gives the desired result.
\end{dow}

\begin{uw}
\label{uwni} Observe that if $f$ does not depend on $z$ then in
Corollary \ref{cor3.10} assumption (Z) and the assumptions that
$(Y^{i},Z^{i})\in \mathbb{L}^{q}(\FF)\otimes\mathbb{L}^{q}(\FF)$
for some $q>\alpha$ are superfluous. This follows from the fact
that $I_{t}$ defined by (\ref{eq3.09}) is less or equal to zero if $f$
does not depend on $z$.
\end{uw}

The proofs of the following lemma and proposition are analogous to
those of Lemma 3.1 and Proposition 3.2 in \cite{BDHPS}, the only
difference being in the fact that we use the It\^o's formula
proved in Proposition \ref{prop2.1} instead of It\^o's formula
proved in \cite{BDHPS}.

\begin{lm}\label{lm3.1}
Let assumption \mbox{\rm(A)}  hold and let $(Y,Z)$ be a solution
of BSDE$(\xi,f+dV)$ with $f,V$ such that
\begin{equation}
\label{eq3.8} E(\into f_{s}\,ds)^{p}+E|V|^{p}_{T}<\infty
\end{equation}
for some $p>0$. If $Y\in \mathcal{D}^{p}$ then $Z\in M^{p}$ and
there exists $C$ depending only on $p$ such that for every $a\ge
\mu+\lambda^{2}$,
\begin{align*}
E(\int_{0}^{T} e^{2as}|Z_{s}|^{2}\,ds)^{p/2} &\le
CE\bigg(\sup_{t\le T} e^{apt}|Y_{t}|^{p}
+(\int_{0}^{T}e^{as}f_{s}\,ds)^{p}+(\into e^{as}
d|V|_{s})^{p}\bigg).
\end{align*}
\end{lm}

\begin{stw}\label{prop3.1}
Let assumption \mbox{\rm(A)} hold and let $(Y,Z)$ be a solution of
BSDE$(\xi,f+dV)$ with $f,V$ satisfying \mbox{\rm(\ref{eq3.8})} for
some $p>1$. If $Y\in \mathcal{D}^{p}$ then there exists $C$
depending only on $p$ such that for every $a\ge
\mu+\lambda^{2}/[1\wedge (p-1)]$,
\begin{align*}
&E\sup_{t\le T}e^{apt}|Y_{t}|^{p}+E(\int_{0}^{T}
e^{2as}|Z_{s}|^{2}\,ds)^{p/2}\\
&\qquad\le CE\bigg(e^{apT}|\xi|^{p}
+(\int_{0}^{T}e^{as}f_{s}\,ds)^{p} +(\into
e^{as}\,d|V|_{s})^{p}\bigg).
\end{align*}
\end{stw}

\begin{stw}
\label{prop.gen} \mbox{\rm(i)} Let $p>1$ and let $(Y,Z)\in
\DM^{p}\otimes M^{p}$ be a solution of BSDE$(\xi,f+dV)$ with
$\xi,f,V$ satisfying \mbox{\rm(H1)--(H3)}. Then there exists $C$
depending only on $\mu^{+},\lambda,T, p$ such that
\begin{align*}
&E(\into|f(s,Y_{s},Z_{s})|\,ds)^{p}\\
&\qquad\le CE\left( |Y^{*}_{T}|^{p} +(\into|Z_{s}|^{2}\,ds)^{p/2}
+(\into|f(s,0,0)|\,ds)^{p}+|V|^{p}_{T}\right).
\end{align*}
\mbox{\rm(ii)} Let $p=1$ and let $(Y,Z)$ such that
$(Y,Z)\in\DM^{q}\otimes M^{q}$ for  $q\in (0,1)$ and  $Y$ is of
class \mbox{\rm(D)} be a solution of BSDE$(\xi,f+dV)$ with
$\xi,f,V$ satisfying \mbox{\rm(H1)--(H3), (Z)}. Then there exists
$C$ depending only on $\mu^{+},\lambda,T$ such that
\begin{align*}
E\into|f(s,Y_{s},Z_{s})|\,ds&\le CE\left( \|Y\|_{1}
+\gamma\into (g_{t}+|Y_{s}|+|Z_{s}|)^{\alpha}\,ds\right.\\
&\quad+\left.\into|f(s,0,0)|\,ds+|V|_{T}\right).
\end{align*}
\end{stw}
\begin{dow}
By It\^o's formula,
\begin{equation}
\label{eq4.4}
-\intot \hat{Y}_{s}f(s,Y_{s},Z_{s})\,ds\le |Y_{t}|-|Y_{0}|
+\intot \hat{Y}_{s-}\,dV_{s}-\intot \hat{Y}_{s}\,dB_{s}.
\end{equation}
Since by (H3), $-\hat{Y}_{s}(f(s,Y_{s},0)-\mu Y_{s})\ge
\hat{Y}_{s}f(s,0,0)$, we have
\[
-\hat{Y}_{s}f(s,Y_{s},Z_{s})+\hat{Y}_{s}
(f(s,Y_{s},Z_{s})-f(s,Y_{s},0))-\hat{Y}_{s}f(s,0,0)+\mu
\hat{Y}_{s} Y_{s}\ge0.
\]
Hence, by (\ref{eq4.4}),
\begin{align*}
\intot|f(s,Y_{s},Z_{s})|\,ds&\le |Y_{t}|+|V|_{t} -\intot
\hat{Y}_{s}\,dB_{s}
+2\intot|f(s,Y_{s},Z_{s})-f(s,Y_{s},0)|\,ds\\
&\quad+ 2\intot |f(s,0,0)|\,ds +2\mu \intot |Y_{s}|\,ds\\
&\le |Y_{t}|+|V|_{t}-\intot\hat{Y}_{s}\,dB_{s} +2\intot|Z_{s}|\,ds
+ 2\intot |f(s,0,0)|\,ds\\
&\quad +2\mu \intot |Y_{s}|\,ds,
\end{align*}
from which one can easily get (i). To prove (ii) we use assumption
(Z) to estimate the integral involving
$|f(s,Y_{s},Z_{s})-f(s,Y_{s},0)|$. We then get
\begin{align*}
\nonumber\intot|f(s,Y_{s},Z_{s})|\,ds&\le |Y_{t}|+|V|_{T}
-\intot \hat{Y}_{s}\,dB_{s}
+2\gamma\intot(g_{t}+|Y_{s}|+|Z_{s}|)^{\alpha}\,ds\\
&\quad+ 2\intot|f(s,0,0)|\,ds +2\mu \intot |Y_{s}|\,ds,
\end{align*}
from which (ii) immediately follows.
\end{dow}

\begin{tw}
\label{prop3.2} Let $p>1$. Under assumptions \mbox{\rm(H2), (H3)}
there exists at most one solution $(Y,Z)\in\mathcal{D}^{p}\otimes
M^{p}$ of BSDE$(\xi,f+dV)$.
\end{tw}
\begin{dow}
Follows from Proposition \ref{prop3.4}.
\end{dow}

\begin{tw}
\label{th3.1} Let $p>1$ and \mbox{\rm(H1)--(H5) hold}. Then there
exists a unique solution  $(Y,Z)\in \mathcal{D}^{p}\otimes M^{p}$
of BSDE$(\xi,f+dV)$.
\end{tw}
\begin{dow}
Without loss of generality we may assume that $\mu\le 0$. Let us
assume that $f$ is bounded. By the representation property of the
Brownian filtration there exists a unique process
$(\overline{Y},\overline{Z})\in \mathcal{D}^{p}\otimes M^{p}$ such
that
\[
\overline{Y}_{t}=\intt dV_{s}-\intt \overline{Z}_{s}\,dB_{s},
\quad t\in [0,T].
\]
Put $\tilde{f}(t,y,z)=f(t,y+\overline{Y}_{t},z+\overline{Z}_{t})$
and observe that the data $(\xi,\tilde{f})$ satisfy assumptions
(H1)--(H5). Therefore by \cite[Theorem 4.2]{BDHPS} there exists a
unique solution $(\tilde{Y},\tilde{Z})\in\SM^{p}\otimes M^{p}$ of
the BSDE
\[
\tilde{Y}_{t}=\xi+\intt\tilde{f}(s,\tilde{Y}_{s},\tilde{Z}_{s})\,ds
-\intt \tilde{Z}_{s}\,dB_{s},\quad t\in [0,T].
\]
Clearly the pair $(\tilde{Y}+\overline{Y},\tilde{Z}+\overline{Z})$
is a unique solution of BSDE$(\xi,f+dV)$.

Now suppose that $f$ is bounded from below. Write $f_{n}=f\wedge
n$. Then by the first step of the proof there exists a unique
solution $(Y^{n},Z^{n})$ of BSDE$(\xi,f_{n}+ dV)$. By Proposition
\ref{prop3.4}, $Y^{n}_{t}\le Y^{n+1}_{t}$, $t\in [0,T]$ for
$n\in\mathbb{N}$. Therefore defining $Y_{t}=\sup_{n\ge 0} Y^{n}_{t},\, t\in [0,T]$
we have that
\begin{equation}
\label{limbo1}
Y^{n}_{t}\nearrow Y_{t},\quad t\in [0,T].
\end{equation}
Moreover, by Proposition \ref{prop3.1}, there
exists $C>0$ not depending on $n$ such that
\begin{equation}
\label{eq3.2}
E\sup_{0\le t\le T}|Y^{n}_{t}|^{p}
+E(\into|Z^{n}_{s}|^{2}\,ds)^{p/2}\le C.
\end{equation}
Hence, by Proposition \ref{prop.gen},
\begin{equation}
\label{eq3.3} \sup_{n\in\mathbb{N}}
E(\into|f_{n}(s,Y^{n}_{s},Z^{n}_{s})|\,ds)^{p}<\infty.
\end{equation}
By It\^o's formula,
\begin{align}
\label{eq3.4} &|Y^{n}_{t}-Y^{m}_{t}|^{p}+\frac12 p(p-1) \intt
|Y^{n}_{s}-Y^{m}_{s}|^{p-2} \mathbf{1}_{\{Y^{n}_{s}\neq
Y^{m}_{s}\}} |Z^{n}_{s}-Z^{m}_{s}|^{2}\,ds\nonumber \\
&\qquad =p\intt(f_{n}(s,Y^{n}_{s},Z^{n}_{s})
-f_{m}(s,Y^{m}_{s},Z^{m}_{s}))|Y^{n}_{s}-Y^{m}_{s}|^{p-1}
\,\hat{\mbox{sgn}}({Y^{n}_{s}-Y^{m}_{s}})\,ds\nonumber \\
&\quad\qquad+p\intt
(Z^{n}_{s}-Z^{m}_{s})|Y^{n}_{s}-Y^{m}_{s}|^{p-1}\,
\hat{\mbox{sgn}}({Y^{n}_{s}-Y^{m}_{s}})\,dB_{s},\quad t\in [0,T].
\end{align}
By the Burkholder-Davis-Gundy inequality, (H2) and H\"older's
inequality,
\begin{align}
\label{limbo2} \nonumber E\sup_{0\le t\le T}
|Y^{n}_{t}-Y^{m}_{t}|^{p} &\le \left(E(\into
|Z^{n}_{s}|^{2}\,ds)^{p/2} +E(\into
|Z^{m}_{s}|^{2}\,ds)^{p/2}\right)^{1/p}\\
&\qquad\quad\times\left(E(\into|Y^{n}_{s}-Y^{m}_{s}|^{2(p-1)}\,ds)
^{p/2(p-1)}\right)^{(p-1)/p}\nonumber\\
&\quad+\left(E(\into|f_{n}(s,Y^{n}_{s},0)-f_{m}
(s,Y^{m}_{s},0)|\,ds)^{p}\right)^{1/p}\nonumber \\
&\qquad\quad\times \left(E\sup_{0\le t\le T}
|Y^{n}_{t}-Y^{m}_{t}|^{p}\right)^{(p-1)/p}.
\end{align}
By (\ref{limbo1}) and (\ref{eq3.2}),
\begin{equation}
\label{eq3.5} \lim_{n,m\rightarrow+\infty}
E(\into|Y^{n}_{s}-Y^{m}_{s}|^{2(p-1)}\,ds)^{p/2(p-1)}=0.
\end{equation}
Therefore by (\ref{eq3.2}) the first term on the right-hand side
of (\ref{limbo2}) converges to zero. By monotonicity of $f_{n},f$
with respect to $y$ and monotonicity of the sequence $\{Y^{n}\}$,
\[
f_{1}(t,Y_{t},0)\le f_{n}(t,Y^{n}_{t},0)\le f(t,Y^{1}_{t},0).
\]
Therefore from (H4), (H5), (\ref{limbo1}), (\ref{eq3.2}), (\ref{eq3.3}) it
follows that
\begin{equation}
\label{limbo4}
E(\into|f_{n}(s,Y^{n}_{s},0)-f_{m}(s,Y^{m}_{s},0)|\,ds)^{p'}
\rightarrow0,
\end{equation}
for every $p'<p$. Without loss of generality we may assume that
(\ref{limbo4}) holds true for $p$ in place of $p'$,
%since otherwise we would do the whole above reasoning with $p'$ in place of $p$.
which when combined with (\ref{eq3.2}) implies convergence to zero
of the second term on the right-hand side of inequality
(\ref{limbo2}). Consequently, $Y\in \DM^{p}$ and
\begin{equation}
\label{limbo3}
E\sup_{0\le t\le T}|Y^{n}_{t}-Y_{t}|^{p}\rightarrow 0.
\end{equation}
Since
\begin{align*}
&E(\into (f_{n}(s,Y^{n}_{s},Z^{n}_{s})
-f_{m}(s,Y^{m}_{s},Z^{m}_{s}))
(Y^{n}_{s}-Y^{m}_{s})\,ds)^{p/2}\\&\quad \le
\left(E(\into|f_{n}(s,Y^{n}_{s},Z^{n}_{s})
-f_{m}(s,Y^{m}_{s},Z^{m}_{s})|\,ds)^{p}\right)^{1/2}
\left(E\sup_{0\le t\le T}|Y^{n}_{t}-Y^{m}_{t}|^{p}\right)^{1/2},
\end{align*}
we conclude from (\ref{eq3.3}), (\ref{eq3.4}) with $p=2$ and
(\ref{limbo3}) that
\begin{equation}
\label{eq3.6} \lim_{n,m\rightarrow+\infty}
E(\into|Z^{n}_{s}-Z^{m}_{s}|^{2}\,ds)^{p/2}=0.
\end{equation}
Therefore there exists a process $Z\in M^{p}$ such that
(\ref{eq3.6}) holds with $Z$ in place of $Z^{m}$. From this and
(\ref{limbo4}) we conclude that
\[
\lim_{n\rightarrow+\infty}
E\left(\into|f_{n}(s,Y^{n}_{s},Z^{n}_{s})-f(s,Y_{s},Z_{s})|
\,ds\right)^{p}=0,
\]
which together  with (\ref{limbo3}), (\ref{eq3.6}) shows that
$(Y,Z)$ is a solution of BSDE$(\xi,f+dV)$.

Finally, in the general case, we approximate $f$ by the sequence
$\{f_n\}$, where $f_{n}=f\vee(-n)$, $n\in\mathbb{N}$. By what has
already been proved for each $n$ there exists a unique solution
$(Y^{n},Z^{n})\in\DM^{p}\otimes M^{p}$ of BSDE$(\xi,f_{n}+dV)$.
Repeating arguments from the proof of the previous step shows that
$(Y^{n},Z^{n})$ converges in $\SM^{p}\otimes M^{p}$ to the unique
solution of BSDE$(\xi,f+dV)$.
\end{dow}

\begin{tw}
If $p=1$ and \mbox{\rm(H2), (H3), (Z)} are satisfied then there
exists at most one solution  of BSDE$(\xi,f+dV)$ such that $Y$ is
of class \mbox{\rm(D)} and $Z\in \bigcup_{\beta>\alpha}
M^{\beta}$.
\end{tw}
\begin{dow}
Follows from Corollary \ref{cor3.10}.
\end{dow}

\begin{stw}
\label{prop3.3} Assume that \mbox{\rm(H1)--(H5)} hold with $p=1$
and $f$ does not depend on $z$. Then there exists a solution
$(Y,Z)$ of BSDE$(\xi,f+dV)$ such that $Y$ is of class
\mbox{\rm(D)} and $(Y,Z)\in \bigcap _{\beta<1}
\mathcal{D}^{\beta}\otimes M^{\beta}$.
\end{stw}
\begin{dow}
Standard arguments show that without loss of generality we may
assume that $\mu\le 0$. Set
\[
\xi^{n}=T_{n}(\xi)\quad f_{n}(t,y)=f(t,y)-f(t,0)+T_{n}(f(t,0)),
\quad V^{n}_{t}
=\intot\mathbf{1}_{\{|V|_{s}\le n\}}\,dV_{s}.
\]
By Theorem \ref{th3.1},  for every $n\in\mathbb{N}$ there exists a
solution $(Y^{n},Z^{n})\in \DM^{2}\otimes M^{2}$ of
BSDE$(\xi^{n},f_{n}+dV^{n})$. Let $m\ge n$. Write
$\delta Y=Y^{m}-Y^{n}$,
$\delta Z=Z^{m}-Z^{n}$, $\delta \xi= \xi^{m}-\xi^{n}$ and
\[
\tau_{k}=\inf\{t\in [0,T]; \intot |\delta Z_{s}|^{2}\,ds >k\}\wedge T.
\]
By the It\^o-Tanaka formula,
\begin{align*}
|\delta Y_{t\wedge \tau_{k}} |&\le |\delta Y_{\tau_{k}}|
+\int_{\tau_{k}\wedge t}^{\tau_{k}}\hat{\mbox{sgn}}({\delta
Y_{s}})
(f_{m}(s,Y^{m}_{s})-f_{n}(s,Y^{n}_{s}))\,ds\\
&\quad+\int_{\tau_{k}\wedge t}^{\tau_{k}}
\hat{\mbox{sgn}}(\delta Y_{s})\,d(V^{m}_{s}-V^{n}_{s})
+\int_{\tau_{k}\wedge t}^{\tau_{k}}\hat{\mbox{sgn}}({\delta
Y_{s}})\delta Z_{s}\,dB_{s}\\
&\le|\delta Y_{\tau_{k}}| +\int_{\tau_{k}\wedge t}^{\tau_{k}}
|f_{m}(s,Y^{n}_{s})-f_{n}(s,Y^{n}_{s})|\,ds\\
&\quad+\int_{\tau_{k}\wedge t}^{\tau_{k}}d|V^{m}-V^{n}|_{s}
+\int_{\tau_{k}\wedge t}^{\tau_{k}}\hat{\mbox{sgn}}({\delta
Y_{s}})\delta Z_{s}\,dB_{s} \quad t\in [0,T],
\end{align*}
the last inequality being a consequence of monotonicity of $f_{n}$
with respect to $y$. Conditioning with respect to $\FF_{t}$, using
the fact that $\delta Y$ is  of class (D) and the definitions of
$f_{n}$, $\xi^{n}$, $V^{n}$ we conclude from the above inequality
that
\[
|\delta Y_{t}|\le E^{\FF_{t}}(|\xi|\mathbf{1}_{\{|\xi|>n\}}
+\into|f(s,0)|\mathbf{1}_{\{|f(s,0)|>n\}}\,ds
+\into\mathbf{1}_{\{|V|_{s}>n\}}\,d|V|_{s}).
\]
Now repeating step by step the arguments following Eq. (12) in the
proof of \cite[Proposition 6.4]{BDHPS} we get the existence
result.
\end{dow}

\begin{tw}
\label{th3.2} Assume that $p=1$ and \mbox{\rm(H1)--(H5), (Z)} are
satisfied. Then there exists a solution $(Y,Z)$ of
BSDE$(\xi,f+dV)$ such that $Y$ is of class \mbox{\rm(D)} and
$(Y,Z)\in\bigcap _{\beta<1} \mathcal{D}^{\beta}\otimes M^{\beta}$.
\end{tw}
\begin{dow}
Without loss of generality we may assume that $\mu\le 0$. Let
$(Y^{0},Z^{0})=(0,0)$. By Proposition \ref{prop3.3} we can define
recursively the sequence $\{(Y^n,Z^n)\}$ by putting
\begin{equation}
\label{eq3.10}
Y^{n+1}_{t}=\xi+\intt f(s,Y^{n+1}_{s},Z^{n}_{s})\,ds
-\intt dV_{s}-\intt Z^{n+1}_{s}\,dB_{s},\quad t\in [0,T].
\end{equation}
Since
\[
Y^{n+1}_{t}-Y^{n}_{t}=\intt (f(s,Y^{n+1}_{s},Z^{n}_{s})
-f(s,Y^{n}_{s},Z^{n-1}_{s}))\,ds
-\intt(Z^{n+1}_{s}-Z^{n}_{s})\,dB_{s}
\]
for $t\in[0,T]$, repeating step by step the proof of Theorem 6.3
in \cite{BDHPS} shows that $(Y^{n},Z^{n})$ converges to some
process $(Y,Z)$ belonging to $\DM^{q}\otimes M^{q}$ for $q\in
(0,1)$ and that $Y^{n}\rightarrow Y$ in the norm $\|\cdot\|_{1}$.
Therefore passing to the limit in (\ref{eq3.10}) in ucp topology
we see that $(Y,Z)$ is a solution of BSDE$(\xi,f+dV)$.
\end{dow}

\nsubsection{Supersolutions of BSDEs} \label{sec4}

In this section we investigate supersolutions of BSDEs. In particular we provide a priori estimates for
supersolutions, Snell envelope representation result for minimal
supersolutions and explicit formula for its jumps. Moreover, we
prove some useful technical lemmas which  generalize known results
on monotone convergence of solutions of BSDEs and regularity
properties of monotone limits of supersolutions. The results on
supersolutions play a pivotal role in the study of reflected
RBSDEs because one can regard solution of reflected BSDE with one
barrier as a minimal supersolution of some BSDE and view solution
of reflected BSDE with two barriers as a minimal supersolution of
some BSDE of the form (\ref{eq3.1}).

Let us fix a process $V\in\mathcal{V}$.
\begin{df}
We say that a pair of processes $(Y,Z)$ is a supersolution  (resp.
subsolution) of BSDE$(\xi,f+dV)$ if
\begin{enumerate}
\item [a)]$Z\in M$, $t\mapsto f(t,Y_{t},Z_{t})\in\mathbb{L}^{1}(0,T)$,
\item [b)]There exists a process $K\in \mathcal{V}^{+}$
(resp. $K\in\mathcal{V}^{-}$) such that
\[
Y_{t}=\xi+\intt f(s,Y_{s},Z_{s})\,ds +\intt dV_{s}+\intt dK_{s}
-\intt Z_{s}\,dB_{s},\quad t\in [0,T].
\]
\end{enumerate}
\end{df}

Suppose that $(Y,Z)$ is a supersolution of some BSDE with data
$(\xi,f,V)$. In the rest of this section $K$ stands for the
increasing c\`adl\`ag process such that $K_{0}=0$ and the above
equation is satisfied.

The following Lemma \ref{lm4.1}, Lemma \ref{lm4.2} and Proposition
\ref{prop4.3} were proved in \cite{Kl4} (Lemma 3.1, Lemma 3.2 and
Proposition 3.4, respectively) in the case where the measure $dV$
is absolutely continuous with respect to the Lebesgue measure.
Using the  It\^o-Tanaka formula proved in Proposition
\ref{prop2.1} one can prove these results for general $V$ by the
same method as in \cite{Kl4}.

\begin{lm}\label{lm4.1}
Let $(Y,Z)$ be a supersolution of BSDE$(\xi,f+dV)$. Assume that
\mbox{\rm(H3)} is satisfied, there exists a progressively
measurable process $X$ such that $X_{t}\ge Y_{t}$, $t\in [0,T]$
and the mappings $[0,T]\ni t\rightarrow X^{+}_t$, $[0,T]\ni
t\rightarrow f^{-}(t,X_{t},0)$ belong to $\mathbb{L}^{1}(0,T)$,
$P$-a.s..
\begin{enumerate}
\item[\rm(i)]If \mbox{\rm(H2)} is satisfied then for every
$\tau\in\mathcal{T}$  and $a\ge \mu$,
\begin{align*}
\int_{0}^{\tau} e^{at}dK_{t}&\le |e^{a\tau}Y_{\tau}|
+|Y_{0}|+\int_{0}^{\tau}e^{as}Z_{s}dB_{s}
+\lambda\int_{0}^{\tau}e^{as}|Z_{s}|\,ds\\&\quad
+\int_{0}^{\tau}e^{as}(f^{-}(s,X_{s},0)\,ds+dV^{-}_{s})
+\int_{0}^{\tau}a^{+}e^{as}X_{s}^{+}\,ds.
\end{align*}
\item[\rm(ii)]If \mbox{\rm(Z)} is satisfied then for every $\tau\in\mathcal{T}$
and $a\ge\mu$,
\begin{align*}
\int_{0}^{\tau} e^{at}dK_{t}&\le |e^{a\tau}Y_{\tau}|+|Y_{0}|
+\int_{0}^{\tau}e^{as}Z_{s}dB_{s}+\gamma\int_{0}^{\tau}e^{as}(g_{s}
+|Y_{s}|+|Z_{s}|)^{\alpha}\,ds\\
&\quad+\int_{0}^{\tau}e^{as}(f^{-}(s,X_{s},0)\,ds+dV^{-}_{s})
+\int_{0}^{\tau}a^{+}e^{as}X_{s}^{+}\,ds.
\end{align*}
\end{enumerate}
\end{lm}

\begin{lm}\label{lm4.2}
Let $(Y,Z)$ be a supersolution of BSDE$(\xi,f+dV)$. If
\mbox{\rm(A)} is satisfied and for some $p>0$,
$Y\in\mathcal{D}^{p}$, \mbox{\rm(H1)} is satisfied and
\[
E(\int_{0}^{T}X^{+}_{s}\,ds)^{p}
+E(\int_{0}^{T}f^{-}(s,X_{s},0)\,ds)^{p}<\infty
\]
for some progressively measurable process $X$ such that $X_{t}\ge
Y_{t}$, $t\in [0,T]$, then $Z\in M^{p}$ and there exists $C$
depending only on $\lambda,p,T$ such that for every $a\ge
\mu+\lambda^{2}$,
\begin{align*}
&E((\int_{0}^{T} e^{2as}|Z_{s}|^{2}\,ds)^{p/2}+(\int_{0}^{T}
e^{as}\,dK_{s})^{p})\le CE\bigg(\sup_{t\le T} e^{apt}|Y_{t}|^{p}
+(\int_{0}^{T}e^{as}|f(s,0,0)|\,ds)^{p}\\
&\qquad\qquad\qquad\quad+(\into e^{as}
d|V|_{s})^{p}+(\int_{0}^{T}e^{as}f^{-}(s,X_{s},0)\,ds)^{p}
+(\int_{0}^{T}a^{+}e^{as}X_{s}^{+}\,ds)^{p}\bigg).
\end{align*}
\end{lm}

\begin{stw}\label{prop4.3}
Assume \mbox{\rm(A)} and let $(Y,Z)$ be a supersolution of
BSDE$(\xi,f+dV)$. If for some $p>1$, $Y\in\mathcal{D}^{p}$,
\mbox{\rm(H1)} is satisfied and
\begin{equation}
\label{eq4.10} E(\int_{0}^{T}f^{-}(s,X_{s},0)\,ds)^{p}<\infty
\end{equation}
for some $X\in \mathcal{D}^{p}$ such that $X_{t}\ge Y_{t}$, $t\in
[0,T]$, then there exists $C$ depending only on $\lambda,p,T$ such
that for every $a\ge\mu+\lambda^{2}/[1\wedge (p-1)]$ and every
$\tau\in\mathcal{T}$,
\begin{align*}
&E\sup_{t\le \tau}e^{apt}|Y_{t}|^{p}+E(\int_{0}^{\tau}
e^{2as}|Z_{s}|^{2}\,ds)^{p/2}
+E(\int_{0}^{\tau}e^{as}\,dK_{s})^{p}\\
&\qquad \le CE\bigg(e^{ap\tau}|Y_{\tau}|^{p}
+(\int_{0}^{\tau}e^{as}|f(s,0,0)|\,ds)^{p}
+(\int_{0}^{\tau}e^{as}\,d|V|_{s})^{p} +\sup_{t\le \tau}
|e^{at}X^{+}_{t}|^{p}\\
&\qquad\quad+(\int_{0}^{\tau}e^{as}f^{-}(s,X_{s},0)\,ds)^{p}
+(\int_{0}^{\tau}a^{+}e^{as}X^{+}_{s}\,ds)^{p}\bigg).
\end{align*}
Assume additionally that $f$ does not depend on $z$. Then if
\mbox{\rm(H1)} and \mbox{\rm(\ref{eq4.10})} are satisfied with
$p=1$ and $X,Y$ are of class \mbox{\rm(D)} then for every
$a\ge\mu$,
\begin{align*}
\|e^{a\cdot} Y\|_{1}& +E\int_{0}^{T}e^{as}\,dK_{s}\le
E\left(e^{aT}|\xi|+\int_{0}^{T}e^{as}|f(s,0)|\,ds\right.\\
&\quad+ \into e^{as} d|V|_{s}+\int_{0}^{T}e^{as}f^{-}(s,X_{s})\,ds
+\left.\int_{0}^{T}a^{+}e^{as}X^{+}_{s}\,ds\right)+\|e^{a\cdot}X^{+}\|_{1}.
\end{align*}
\end{stw}

\begin{lm}
\label{lm4.4}
Let $(Y,Z)$ be a supersolution of BSDE$(\xi,f+dV)$.
Assume that
\begin{enumerate}
\item [\rm(a)]$f$ does not depend on $y,z$,  \mbox{\rm(H1)}
with $p=1$ is satisfied.
\item [\rm(b)]$Y_{t}^{n}\nearrow Y_{t}$, $t\in [0,T]$, $Y^{n},Y$
are of class \mbox{\rm(D)}, $Y_{t}\ge L_{t}$ for a.e. $t\in
[0,T]$, $dV^{n}\le dV$, where
\[
Y_{t}^{n}=\xi+\intt f(s)\,ds
+\intt dV^{n}_{s}+\intt n(Y_{s}^{n}-L_{s})^{-}\,ds
-\intt Z^{n}_{s}\,dB_{s},\quad t\in [0,T].
\]
\end{enumerate}
Then $(Y,Z)$ is the smallest supersolution of BSDE$(\xi,f+dV)$
such that $L_{t}\le Y_{t}$ for a.e. $t\in [0,T]$ and $Y$ is of
class \mbox{\rm(D)}.
\end{lm}
\begin{dow}
Let $(\bar{Y},\bar{Z})$ be a supersolution of BSDE$(\xi,f+dV)$
such that $\bar{Y}$ is of class (D) and $L_{t}\le \bar{Y}_{t}$ for
a.e. $t\in [0,T]$. Then there exists
$\bar{K}\in\mathcal{V}^{+,1}$ such that
\[
\bar{Y}_{t}=\xi+\intt f(s)\,ds+\intt dV_{s}
+\intt d\bar{K}_{s}-\intt\bar{Z}_{s}\,dB_{s},\quad t\in [0,T].
\]
Since $\bar{Y}_{t}\ge L_{t}$ for a.e. $t\in [0,T]$, we have
\[
\bar{Y}_{t}=\xi+\intt f(s)\,ds+\intt dV_{s}
+\intt d\bar{K}_{s}+\intt n(\bar{Y}_{s}-L_{s})^{-}\,ds
-\intt\bar{Z}_{s}\,dB_{s},\quad t\in [0,T].
\]
By Corollary \ref{cor3.10},  $Y^{n}_{t}\le \bar{Y}_{t}$, $t\in
[0,T]$, and consequently $Y_{t}\le \bar{Y}_{t}$, $t\in [0,T]$.
\end{dow}

\begin{wn}
\label{col4.5} Let $Y$ be the process of Lemma \ref{lm4.4} and let
assumptions of Lemma \ref{lm4.4} hold. Then for every
$\hat{L}\in\mathcal{D}$ such that $L_{t}\le\hat{L}_{t}\le Y_{t}$
for a.e. $t\in[0,T]$,
\[
Y_{t}=\esssup_{\tau\in \mathcal{T}_{t}}E(\int_{t}^{\tau}f(s)\,ds
+\int_{t}^{\tau}dV_{s}+\hat{L}_{\tau}\mathbf{1}_{\{\tau<T\}}
+\xi\mathbf{1}_{\{\tau=T\}}|\FF_{t}),\quad t\in[0,T].
\]
\end{wn}
\begin{dow}
It suffices to observe that from Lemma \ref{lm4.4} it follows that
for every $\hat{L}\in\mathcal{D}$ such that $L_{t}\le
\hat{L}_{t}\le Y_{t}$ for a.e. $t\in [0,T]$ the process
$\bar{Y}_{t}\equiv Y_{t}+\intot f(s)\,ds +\intot dV_{s}$ is the
smallest supermartingale majorazing
$\bar{L}_{t}\equiv\hat{L}_{t}+\intot f(s)\,ds +\intot dV_{s}$ such
that $\bar{Y}_{T}=\xi+\into f(s)\,ds +\into dV_{s}$.
\end{dow}

\begin{lm}
\label{lm.min} Assume that \mbox{\rm(H1)} holds with $p=1$, $L\in
\DM $ is of class \mbox{\rm(D)} and $Y\in\DM$ is of the form
\[
Y_{t}=\esssup_{\tau\in \mathcal{T}_{t}}E(\int_{t}^{\tau}f(s)\,ds
+\int_{t}^{\tau}dV_{s}+L_{\tau}\mathbf{1}_{\{\tau<T\}}
+\xi\mathbf{1}_{\{\tau=T\}}|\FF_{t}),\quad t\in[0,T].
\]
Then
\[
Y_{t-}=L_{t-}\vee (Y_{t}+\Delta V_{t}), \quad t\in(0,T].
\]
\end{lm}
\begin{dow}
Let us fix $t\in (0,T]$. By the properties of the Snell envelope,
for every $s\in [0,t)$,
\[
Y_{s}=\esssup_{s\le \tau\le t}E(\int_{s}^{\tau}f(s)\,ds
+\int_{s}^{\tau} dV_{r}+L_{\tau}\mathbf{1}_{\{\tau<t\}}
+Y_{t}\mathbf{1}_{\{\tau=t\}}|\FF_{s}).
\]
Letting $s\rightarrow t^{-}$ we get the desired result.
\end{dow}

\begin{lm}
\label{lm4.5544} Let $Y$ be a nonnegative supermartingale of class
\mbox{\rm(D)}. Then there exists a stationary sequence
$\{\tau_{k}\}\subset\mathcal{T}$ of stopping times such that
$Y_{\tau_{k}}\le Y_{0}\vee k$, $k\in\mathbb{N}$.
\end{lm}
\begin{dow}
Since $Y$ is of class (D), there exists $K\in\mathcal{V}^{+,1}$
such that
\[
Y_{t}=Y_{T}+\intt dK_{s}-\intt Z_{s}\,dB_{s},\quad t\in [0,T].
\]
Let $k,l>0$ and let $Y^{l}_{T}=Y_{T}\wedge l$,
$K^{l}_{t}=K_{t}\wedge l$, $t\in[0,T]$. By Theorem \ref{th3.1},
there exists a solution
$(\bar{Y}^{l},\bar{Z}^{l})\in\DM^{2}\otimes M^{2}$ of
BSDE$(Y^{l}_{T},dK^{l})$. One can check that
$\bar{Y}^{l}_{t}\nearrow Y_{t},\, t\in [0,T]$. Since we can regard
$(\bar{Y}^{l},\bar{Z}^{l})$ as a solution of
RBSDE$(Y^{l}_{T},0,\bar{Y}^{l})$ (for the definition of the last
equation see Section \ref{sec5}), it follows from \cite[Theorem
4.2]{Peng} that $\bar{Y}^{n,l}_{t}\nearrow \bar{Y}^{l}_{t}$, $t\in
[0,T]$, where $(\bar{Y}^{n,l},\bar{Z}^{n,l})\in\SM^{2}\otimes
M^{2}$ is a solution of the BSDE
\[
\bar{Y}^{n,l}_t=Y^{l}_{T}+\intt n(\bar{Y}^{n,l}-\bar{Y}^{l})\,ds
-\intt \bar{Z}^{n,l}_{s}\,dB_{s},\quad t\in [0,T].
\]
Put $\tau_{k}=\inf\{t\in [0,T], Y_{t}>k\}\wedge T$. The sequence
$\{\tau_{k}\}$ is stationary. Moreover,
$\bar{Y}^{n,l}_{\tau_{k}}=\bar{Y}^{n,l}_{\tau_{k}-}\le
Y_{\tau_{k}-}\le k$ on $\{\tau_{k}>0\}$. Hence
$\bar{Y}_{\tau_{k}}^{n,l}\le k\vee Y_{0}$, so letting
$n\rightarrow+\infty$ and then $l\rightarrow+\infty$ we get the
desired result.
\end{dow}

\begin{lm}
\label{lm4.55} If $Y\in\mathcal{V}^{1}+\mathcal{M}_{c}^{loc}$ is of class
\mbox{\rm(D)} then there exist a stationary sequence
$\{\tau_{k}\}\subset\mathcal{T}$ and a sequence of constants
$\{c_{k}\}\subset\BR$, $c_{k}=c(k,Y_{0})$, such that
$Y^{*}_{\tau_{k}}\le c_{k}$ for $k\in\BN$.
\end{lm}
\begin{dow}
By the representation property of Brownian filtration and the
assumptions of the lemma there exists $Z\in M,V\in\mathcal{V}^{1}$
such that
\[
Y_{t}=Y_{T}+\intt dV_{s}-\intt Z_{s}\,dB_{s},\quad t\in [0,T]
\]
and $Y_{T}\in\mathbb{L}^{1}(\FF_{T})$. By Proposition
\ref{prop3.3} there exist unique solutions $(Y^{1},Z^{1})$,
$(Y^{2},Z^{2})\in \bigcap_{q<1} \DM^{q}\otimes M^{q}$ of
BSDE$(Y^{+}_{T},dV^{+})$ and BSDE$(Y^{-}_{T},dV^{-})$,
respectively, such that $Y^{1},Y^{2}$ are of class (D). Since
$Y^{1},Y^{2}$ are nonnegative supermartingales, it follows from
Lemma \ref{lm4.5544} that there exist a stationary sequence
$\{\tau_{k}\}\subset\mathcal{T}$ and a sequence
$\{b_{k}\}\subset\BR$ such that
$Y^{-,*}_{\tau_{k}}+Y^{+,*}_{\tau_{k}}\le b_{k}$. Therefore the
result follows from the fact that $Y=Y^{1}-Y^{2}$.
\end{dow}

\begin{lm}
\label{lm4.6} Let $p\ge 1$ and let $(Y,Z)$ be a supersolution of
BSDE$(\xi,f+dV)$. Assume that
\begin{enumerate}
\item [\rm(a)]$\xi,f,L,V$ satisfy assumptions \mbox{\rm(H1)--(H5)},
$L^+\in\mathbb{L}^{\infty,p}(\FF)$,
$(Y,Z)\in\mathcal{D}^{p}\otimes M^{p}$,
$E(\into|f(s,Y_{s},Z_{s})|\,ds)^{p}<\infty$ in case $p>1$ or, in
case $p=1$, $\xi,f,L,V$ satisfy \mbox{\rm(H1)--(H5)} and
\mbox{\rm(Z)}, $L^{+},Y$ are of class \mbox{\rm(D)},
$(Y,Z)\in\mathcal{D}^{q}\otimes M^{q}$ for every $q\in(0,1)$ and
$E\into|f(s,Y_{s},Z_{s})|\,ds<\infty$.
\item [\rm(b)]$V^{n}\in\mathcal{V},\, dV^{n}\le dV$,
$V^{n}_{t}\nearrow V_{t}$, $t\in [0,T]$, $|V^{n}|_{T}\le|V|_{T}$,
$Y^{n}\in\DM^{p}$ (in case $p=1$, $Y^{n}\in\DM^{q}$, $q\in (0,1)$,
$Y^{n}$ is of class \mbox{\rm(D)}) and $Y^{n}_{t}\nearrow Y_{t}$,
$t\in [0,T]$, $Z^{n}\rightarrow Z$ in $\mathbb{L}^{1}(0,T)$ in
probability $P$, where
\[
Y^{n}_{t}=\xi+\intt f(s,Y^{n}_{s},Z^{n}_{s})\,ds
+\intt dV^{n}_{s}+\intt n(Y^{n}_{s}-L_{s})^{-}\,ds
-\intt Z^{n}_{s}\,dB_{s}, t\in [0,T].
\]
\end{enumerate}
Then with the notation  $h(t)=f(t,Y_{t},Z_{t})$, we have
\begin{enumerate}
\item[\rm(i)]$(Y,Z)$ is the smallest supersolution of BSDE$(\xi,f+dV)$
such that $L_{t}\le Y_{t}$ for a.e. $t\in [0,T]$ in the class of
all processes $(Y,Z)$ from $\mathcal{D}^{p}\otimes M^{p}$ such
that $E(\into |f(s,Y_{s},Z_{s})|\,ds)^{p}<\infty$ in case $p>1$
(in the class of all processes $(Y,Z)\in\mathcal{D}^{q}\otimes
M^{q}$ for $q\in (0,1)$ such that $Y$ is of class \mbox{\rm(D)}
and $E\into |f(s,Y_{s},Z_{s})|\,ds<\infty$ in case $p=1$ ).
\item [\rm(ii)]$(Y,Z)$ is the smallest supersolution of the linear
BSDE$(\xi,h+dV)$ such that $L_{t}\le Y_{t}$ for a.e. $t\in[0,T]$
in the class of all processes $(Y,Z)$ such that $Y\in\DM^{p}$ in
case $p>1$ (in the class of processes $(Y,Z)$ such that $Y$ is of
class \mbox{\rm(D)} in case $p=1$).
\item [\rm(iii)]For every $\hat{L}\in\mathcal{D}$ such that
$L_{t}\le \hat{L}_{t}\le Y_{t}$ for a.e. $t\in [0,T]$,
\[
Y_{t}=\esssup_{\tau\in \mathcal{T}_{t}}
E(\int_{t}^{\tau}f(s,Y_{s},Z_{s})\,ds +\int_{t}^{\tau}
dV_{s}+\hat{L}_{\tau}\mathbf{1}_{\{\tau<T\}}
+\xi\mathbf{1}_{\{\tau=T\}}|\FF_{t}),\quad t\in[0,T].
\]
\item [\rm(iv)]For every $\hat{L}\in\mathcal{D}$ such
that $L_{t}\le \hat{L}_{t}\le Y_{t}$ for a.e. $t\in [0,T]$,
\begin{equation}
\label{eq4.01} Y_{t-}=\hat{L}_{t-}\vee (Y_{t}+\Delta V_{t}),\quad
t\in (0,T].
\end{equation}
\end{enumerate}
\end{lm}
\begin{dow}
For fixed $p\ge 0$ by Lemma \ref{lm3.1} and Proposition \ref{prop.gen}, $Z^{n}\in
M^{p}$ if $p>1$, and if $p=1$ then $Z^{n}\in M^{q}$, $q\in (0,1)$
and $E(\into|f(s,Y^{n}_{s},Z^{n}_{s})|\,ds)^{p}<\infty$.

(i) Suppose that $p=1$. Let $(\bar{Y},\bar{Z})$ be a supersolution
of BSDE$(\xi,f+dV)$ such that $\bar{Y}_{t}\ge L_{t}$ for a.e.
$t\in [0,T]$, $\bar{Y}$ is of class (D), $(\bar{Y},\bar{Z})\in
\DM^{q}\otimes M^{q}$, $q\in (0,1)$ and $E\into
|f(s,\bar{Y}_{s},\bar{Z}_{s})|\,ds<\infty.$ Then there exists
$\bar{K}\in\mathcal{V}^{+}$ such that
\[
\bar{Y}_{t}=\xi+\intt f(s,\bar{Y}_{s},\bar{Z}_{s})\,ds +\intt
dV_{s}+\intt n(\bar{Y}_{s}-L_{s})^{-}\,ds +\intt
d\bar{K}_{s}-\intt \bar{Z}_{s}\,dB_{s}
\]
for $t\in[0,T]$ since $\int^T_t(\bar{Y}_{s}-L_{s})^{-}\,ds=0$.
Hence, by Corollary \ref{cor3.10}, $Y^{n}_{t}\le \bar{Y}_{t}$,
$t\in [0,T]$, which implies that $Y_{t}\le\bar{Y}_{t}$,
$t\in[0,T]$. The proof of (i) in case of $p>1$ is similar, so we omit
it.

(ii) Let $p=1$. From Theorem \ref{th3.2} it follows that there
exists a unique solution $(\bar{Y}^{n},\bar{Z}^{n})$ of the BSDE
\[
\bar{Y}^{n}_{t}=\xi+\intt h(s)\,ds+\intt dV^{n}_{s} +\intt
n(\bar{Y}^{n}-L_{s})^{-}\,ds -\intt \bar{Z}^{n}_{s}\,dB_{s},\quad
t\in[0,T]
\]
such that $\bar{Y}^{n}$ is of class (D) and
$(\bar{Y}^{n},\bar{Z}^{n})\in\DM^{q}\otimes M^{q}$ for all $q\in
(0,1)$. Observe that by Corollary \ref{cor3.10},
\begin{equation}
\label{eq4.1}
\bar{Y}^{n}_{t}\le Y_{t},\quad t\in [0,T].
\end{equation}
By Lemma \ref{lm4.55} there exist a stationary sequence
$\{\tau_{k}\}\subset\mathcal{T}$ and constants $c_{k}$ such that
\[
\int_{0}^{\tau_{k}}|h(s)|\,ds +Y^{+,*}_{\tau_{k}}
+\bar{Y}^{1,*}_{\tau_{k}}+\int_{0}^{\tau_{k}} d|V|_{s}\le c_{k}
\]
for $k\in\mathbb{N}$. Therefore from (\ref{eq4.1}) and Proposition
\ref{prop4.3} it follows that for every $k\in\mathbb{N}$,
\begin{align*}
&E\sup_{t\le \tau_{k}}|\bar{Y}^{n}_{t}|^{2}
+E\int_{0}^{\tau_{k}}|\bar{Z}^{n}_{s}|^{2}\,ds+E(\int_{0}^{\tau_{k}}
d\bar{K}^{n}_{s})^{2}\\
&\qquad\le CE \left((\int_{0}^{\tau_{k}} |h(s)|\,ds)^{2}
+(\int_{0}^{\tau_{k}}d |V|_{s})^{2}
+|Y^{+,*}_{\tau_{k}}|^{2}+|\bar{Y}^{1,*}_{\tau_{k}}|^{2}\right)
\le c'_{k}.
\end{align*}
Applying \cite[Theorem 3.1]{PengXu} in each interval
$[0,\tau_{k}]$ and using stationarity of the sequence
$\{\tau_{k}\}$ we conclude that the process $\bar{Y}_{t}\equiv
\sup_{n\ge1}\bar{Y}^{n}_{t}$, $t\in [0,T]$, is c\`adl\`ag and
there exist processes $\bar{K}\in\mathcal{V}^{+}$ and $\bar{Z}\in
M$ such that
\[
\bar{Y}_{t}=\xi+ \intt h(s)\,ds
+\intt dV_{s}+\intt d\bar{K}_{s}-\intt \bar{Z}_{s}\,dB_{s}.
\]
From the above formula and integrability of $\bar{Y},\xi,h,V$ it
follows immediately that $\bar{K}\in \mathcal{V}^{+,1}$, whereas
from Lemma \ref{lm4.2} it follows that $\bar{Z}\in M^{q}$, $q\in
(0,1)$. In view of Lemma \ref{lm4.4} to complete the proof it suffices to show that
$\bar{Y}_{t}=Y_{t}$, $t\in [0,T]$. To this end, let us observe
that by It\^o's formula and monotonicity of the mapping
$x\mapsto(x-L_{t})^{-}$, for every $\tau\in\mathcal{T}$ we have
\begin{align*}
|Y^{n}_{t\wedge\tau}-\bar{Y}^{n}_{t\wedge\tau}| &\le
|Y^{n}_{\tau}-\bar{Y}^{n}_{\tau}|+\int_{\tau\wedge t}^{\tau}
(f(s,Y^{n}_{s},Z_{s}^{n})-f(s,Y_{s},Z_{s}))
\hat{\mbox{sgn}}({Y^{n}_{s}-\bar{Y}^{n}_{s}})\,ds\\
&\quad+n\int_{t\wedge\tau}^{\tau}\hat{\mbox{sgn}}({Y^{n}_{s}
-\bar{Y}^{n}_{s}})((Y^{n}_{s}-L_{s})^{-}
-(\bar{Y}^{n}_{s}-L_{s})^{-})\,ds\\
&\quad +\int_{t\wedge\tau}^{\tau}(Z^{n}_{s}-\bar{Z}^{n}_{s})
\hat{\mbox{sgn}}({Y^{n}_{s}-\bar{Y}^{n}_{s}})\,dB_{s}\\
&\le |Y^{n}_{\tau}-\bar{Y}^{n}_{\tau}|+\int_{t\wedge\tau}^{\tau}
|f(s,Y^{n}_{s},Z_{s}^{n})-f(s,Y_{s},Z_{s})|\,ds\\
&\quad +\int_{t\wedge\tau}^{\tau}(Z^{n}_{s}-\bar{Z}^{n}_{s})
\hat{\mbox{sgn}}({Y^{n}_{s}-\bar{Y}^{n}_{s}})\,dB_{s}.
\end{align*}
Therefore taking expectation of both sides of the above inequality
with $\tau$ replaced by  $\tau_{k}\equiv \inf\{t\in [0,T]$,
$\intot|Z^{n}_{s}-\bar{Z}^{n}_{s}|^{2}\,ds>k\}\wedge T$, letting
$k\rightarrow +\infty$ and using the fact that $Y^{n}-\bar{Y}^{n}$
is of class (D) we obtain
\begin{equation}
\label{eq4.2}
E|Y_{t}^{n}-\bar{Y}^{n}_{t}|
\le E\into |f(s,Y^{n}_{s},Z_{s}^{n})-f(s,Y_{s},Z_{s})|\,ds,
\quad t\in [0,T].
\end{equation}
We now show that the right-hand side of (\ref{eq4.2}) converges to
zero. To this end, let us first observe that by the assumptions on
the convergence of the sequence $\{(Y^{n},Z^{n})\}$, (H2), (H4),
(H5) and the Lebesgue dominated convergence theorem,
\begin{equation}
\label{eq4.17}
\into |f(s,Y_{s}^{n},Z^{n}_{s})-f(s,Y_{s},Z_{s})|\,ds\rightarrow 0
\end{equation}
and the same is true  with $Z$ in place of $Z^{n}$ and $Y^{n}$ in
place of $Y$. Let us note that (\ref{eq4.2}), (\ref{eq4.17}) hold
true for $p>1$ as well. The proof of (\ref{eq4.2}), (\ref{eq4.17})
for $p>1$ is analogous to the above proof for $p=1$, the only
difference being in the fact that in case $p>1$ assumption (Z) is
not needed and the processes $(\bar{Y}^{n},\bar{Z}^{n}),
(\bar{Y},\bar{Z})$ considered above belong to $\DM^{p}\otimes
\mathcal{M}^{p}$. If $p>1$ then (\ref{eq4.17}) implies convergence
of the right-hand side of (\ref{eq4.2}) to zero because by
Proposition \ref{prop4.3} with $X=Y$ and Proposition
\ref{prop.gen},
\[
\sup_{n\ge1}E(\into|f(s,Y^{n}_{s},Z^{n}_{s})|\,ds)^{p}<\infty.
\]
If $p=1$ then by (H3), monotonicity of the sequence $\{Y^{n}\}$
and (Z),
\begin{align*}
f(s,Y_{s},Z_{s})\le f(s,Y^{n}_{s},Z_{s})&\le
f(s,Y^{1}_{s},Z_{s})\\
&\le |f(s,Y^{1}_{s},Z^{1}_{s})| +2\gamma
(g_{s}+|Y^{1}_{s}|+|Z_{s}|+|Z^{1}_{s}|)^{\alpha}
\end{align*}
Hence, by (a) and the remark at the beginning of the proof,
\[
E\into |f(s,Y^{1}_{s},Z_{s})|\,ds <\infty,
\]
which when combined with the fact that $h\in\mathbb{L}^{1}(\FF)$
allows us to apply the Lebesgue dominated convergence theorem to
get
\begin{equation}
\label{eq4.12}
E\into|f(s,Y^{n}_{s},Z_{s})-f(s,Y_{s},Z_{s})|\,ds\rightarrow 0.
\end{equation}
By (Z) we also have
\[
|f(s,Y^{n}_{s},Z^{n}_{s})-f(s,Y^{n}_{s},Z_{s})|
\le 2\gamma (g_{s}+|Y^{n}_{s}|+|Z^{n}_{s}|+|Z_{s}|)^{\alpha}.
\]
Therefore by Lemma \ref{lm4.2} applied to $(Y^{n},Z^{n})$ (with
$X=Y$) there exists $p>1$ such that
\[
\sup_{n\ge1} E\left(\into
|f(s,Y^{n}_{s},Z^{n}_{s})-f(s,Y^{n}_{s},Z_{s})|\,ds\right)^{p}<\infty,
\]
which when combined with (\ref{eq4.17}) gives
\begin{equation}
\label{eq4.21} E\into
|f(s,Y^{n}_{s},Z^{n}_{s})-f(s,Y^{n}_{s},Z_{s})|\,ds\rightarrow0.
\end{equation}
By (\ref{eq4.12}), (\ref{eq4.21}) the right-hand side of
(\ref{eq4.2}) converges to zero. Therefore $\bar{Y}_{t}= Y_{t}$,
$t\in [0,T]$, and the proof of (ii) is complete.

Assertion (iii) follows from (ii) and Corollary \ref{col4.5},
whereas (iv) follows from (iii) and Lemma \ref{lm.min}.
\end{dow}
\medskip

We close this section with very useful theorem on monotone
convergence of semimartingales. The theorem generalizes
\cite[Theorem 3.1]{PengXu} (see also \cite[Theorem 2.1]{Peng}). In
the proof we will need the following two lemmas.

\begin{lm}
\label{LF1} Assume that $\{(Y^{n}, X^{n},K^{n},A^{n}\}$ is a
sequence of progressively measurable processes such that
\[
Y^{n}_{t}=-K_{t}^{n}+A^{n}_{t}+ X^{n}_{t},\quad t\in [0,T]
\]
for $n\in\mathbb{N}$ and
\begin{enumerate}
\item[\rm(a)]$A^{n}, K^{n}\in \mathcal{V}^{+,1},$
\item[\rm(b)]$\{dA^{n}\}$ is increasing, $A^{n}_{t}\rightarrow A_{t}$,
$t\in [0,T]$, $EA_{T}<\infty$,
\item[\rm(c)]$Y^{n}_{t}\nearrow Y_{t}$, $t\in [0,T]$, $Y^{1},Y$ are of
class \mbox{\rm(D)}, $Y_{0}=0$,
\item[\rm(d)]There exists a c\`adl\`ag process $X$ of
class \mbox{\rm(D)} such that $X_{0}=0$ and for some subsequence
$\{n'\}$, $X^{n'}_{\tau}\rightarrow X_{\tau}$ weakly in
$\mathbb{L}^{1}(\FF_{T})$ for every $\tau\in\mathcal{T}$.
\end{enumerate}
Then $Y\in\DM$, $A\in\mathcal{V}^{+}$, there exists
$K\in\mathcal{V}^{+}$ such that $K^{n'}_{\tau}\rightarrow
K_{\tau}$ weakly in $\mathbb{L}^{1}(\FF_{T})$ for every
$\tau\in\mathcal{T}$ and
\[Y_{t}=-K_{t}+A_{t}+X_{t},\quad t\in [0,T].\]
\end{lm}
\begin{dow}
Put $K_{t}=A_{t}+X_{t}-Y_{t}$, $t\in [0,T]$. By (c),
$Y^{n'}_{\tau}\rightarrow Y_{\tau}$ weakly in
$\mathbb{L}^{1}(\FF_{T})$ for every $\tau\in\mathcal{T}$, and
hence, by (b)--(d), $K^{n'}_{\tau}\rightarrow K_{\tau}$ weakly in
$\mathbb{L}^{1}(\FF_{T})$ for every $\tau\in\mathcal{T}$. Since
$K^{n}_{\sigma}\le K^{n}_{\tau}$,  for any
$\sigma,\tau\in\mathcal{T}$ such that $\sigma\le\tau$, it follows
that $K_{\sigma}\le K_{\tau}$, hence that $K$ is increasing.
Finally, by Lemmas  3.1 and  3.2 in \cite{PengXu}, $A,K,Y\in\DM$.
\end{dow}

\begin{lm}
\label{dombsde} Assume \mbox{\rm(H2)--(H5)}. Let
$L^{n},L\in\mathcal{V},\, g_{n}, g, \bar{f}\in\mathbb{L}^{1}(\FF)$
and let $(Y^{n},Z^{n})$, $(Y,Z)\in\DM\otimes M$ be processes such
that $t\mapsto f(t,Y^{n}_{t},Z^{n}_{t})$, $t\mapsto
f(t,Y_{t},Z_{t})\in\mathbb{L}^{1}(0,T)$ and
\[
Y^{n}_{t}=Y^{n}_{0}-\intot g_{n}(s)\,ds
-\intot f(s,Y^{n}_{s},Z^{n}_{s})\,ds
-\intot dL^{n}_{s}+\intot Z^{n}_{s}\, dB_{s}, \quad t\in [0,T],
\]
\[
Y_{t}=Y_{0}-\intot g(s)\,ds-\intot \bar{f}(s)\,ds
-\intot dL_{s}+\intot Z_{s}\, dB_{s}, \quad t\in [0,T].
\]
Then if
\begin{enumerate}
\item[\rm(a)]$E\sup_{n\ge 0}(L^{n})^{+}_{T}
+E\into |f(s,0,0)|\,ds<\infty$,
\item[\rm(b)]$\liminf_{n\rightarrow+\infty}
\int_{\sigma}^{\tau}(Y_{s}-Y^{n}_{s})\,dL^{n}_{s}\ge 0$ for
every $\sigma,\tau\in\mathcal{T},\, \sigma\le \tau,$
\item[\rm(c)]There exists $C\in\mathcal{V}^{1,+}$ such that
$|\Delta(Y_{t}-Y^{n}_{t})|\le|\Delta C_{t}|,\, t\in [0,T],$
\item [\rm(d)]There exist $\underline{y},\overline{y}
\in\mathcal{V}^{1,+}+\MM_{loc}$ of class \mbox{\rm(D)} such that
\[
\overline{y}_{t}\le Y_{t}\le \underline{y}_{t},\,t\in [0,T],
\quad E\into f^{+}(s,\overline{y}_{s},0)\,ds
+E\into f^{-}(s,\underline{y}_{s},0)\,ds<\infty,
\]
\item[\rm(e)]There exists $h\in\mathbb{L}^{1}(\FF)$ such that
$|g_{n}(s)|\le h(s)$ for a.e. $s\in [0,T]$,
\item [\rm(f)] $Y^{n}_{t}\rightarrow Y_{t}$, $t\in [0,T]$,
\end{enumerate}
then
\begin{equation}
\label{eq.pm} Z^{n}\rightarrow Z,\quad \lambda\otimes P
\mbox{-a.e.},\quad
\int_{0}^{\cdot}|f(s,Y^{n}_{s},Z^{n}_{s})-f(s,Y_{s},Z_{s})|\,ds
\rightarrow 0\quad\mbox{in ucp}
\end{equation}
and there exists a stationary sequence
$\{\tau_{k}\}\subset\mathcal{T}$ such that for every
$k\in\mathbb{N}$ and $p\in(0,2)$,
\begin{equation}
\label{eq.p0}
E\int_{0}^{\tau_{k}}|Z^{n}_{s}-Z_{s}|^{p}\,ds\rightarrow 0,
\end{equation}
and if $|\Delta C_{t}|=0$, $t\in [0,T]$, then the above
convergence holds for $p=2$, too. If additionally
$g_{n}\rightarrow g$ weakly in $\mathbb{L}^{1}(\FF)$ and
$L^{n}_{\tau}\rightarrow L_{\tau}$ weakly in
$\mathbb{L}^{1}(\FF_{T})$ for every $\tau\in\mathcal{T}$, then
$\bar{f}(s)=f(s,Y_{s},Z_{s})$ for a.e. $s\in [0,T]$.
\end{lm}
\begin{dow}
%We divide the proof into three steps.
%\smallskip\\
Step 1. (Reduction to the study of stopped processes) Assume that
$\{\tau_{k}\}\subset \mathcal{T}$ is a stationary sequence. Write
$Y^{n,k}_{t}=Y^{n}_{t\wedge\tau_{k}}$,
$L^{n,k}_{t}=L^{n}_{t\wedge\tau_{k}}$,
$Z^{n,k}=Z^{n}\mathbf{1}_{[0,\tau_{k}]}$,
$\eta^{n,k}=f(\cdot,Y^{n,k},Z^{n,k})\mathbf{1}_{[0,\tau_{k}]}$,
$g^{n,k}=g^{n}\mathbf{1}_{[0,\tau_{k}]}$,
$Y^{k}_{t}=Y_{t\wedge\tau_{k}}$, $L^{k}_{t}=L_{t\wedge\tau_{k}}$,
$Z^{k}=Z\mathbf{1}_{[0,\tau_{k}]}$ and
$\eta^{k}=\bar{f}\mathbf{1}_{[0,\tau_{k}]}$,
$g^{k}=g\mathbf{1}_{[0,\tau_{k}]}$.  Then
\[
Y^{n,k}_{t}=Y^{n,k}_{0}-\intot g^{n,k}(s)\,ds -\intot
f^{k}(s,Y^{n,k}_{s},Z^{n,k}_{s})\,ds-\intot dL^{n,k}_{s} +\intot
Z^{n,k}_{s}\,dB_{s}
\]
and
\[
Y^{k}_{t}=Y^{k}_{0}-\intot g^{k}(s)\,ds-\intot \bar{f}^{k}(s)\,ds
-\intot dL^{k}_{s}+\intot Z^{k}_{s}\,dB_{s}
\]
for $t\in[0,T]$. Since $\{\tau_{k}\}$ is stationary, it follows
from the above that it suffices to prove the lemma for stopped
processes $(Y^k,Z^k)$.

Step 2. (Localization) By Lemma \ref{lm4.55} there exist a
stationary sequence $\{\delta^{1}_{k}\}\subset\mathcal{T}$ and
constants $c_{k}$ such that $\underline{y}^{*}_{\delta^{1}_{k}}\le
c_{k}$ and $\bar{y}^{*}_{\delta^{1}_{k}}\le c_{k}$. Let $\bar{D}$
be a c\`adl\`ag version of the process $D_{t}=\sup_{n\ge
0}(L^{n})^{+}_{t}$, $t\in[0,T]$. Then $(L^{n})^{+}_{t}\le
\bar{D}_{t}$, $t\in [0,T]$ and by assumption (a),
$E\bar{D}_{T}<\infty$. Therefore by Lemma \ref{lm4.55} there exist
a stationary sequence $\{\delta^{2}_{k}\}\subset\mathcal{T}$ and
constants $b_{k}$ such that $\bar{D}_{\delta^2_k}\le b_{k}$. Put
$\tau_{k}=\min\{\delta^{1}_k,\delta^{2}_k,\delta^{3}_k\}$, where
\begin{align*}
\delta^{3}_{k}&=\inf\{t\in [0,T];
\intot |f(s,0,0)|\,ds+\intot|h(s)|\,ds \\
&\qquad\qquad+\intot f^{-}(s,\underline{y}_{s},0)\,ds +\intot
f^{+}(s,\bar{y}_{s},0)\,ds>k\}\wedge T.
\end{align*}
From Proposition \ref{prop.gen}, Proposition \ref{prop4.3} and the
definition of $\tau_{k}$ it follows that there exists $C$ not
depending on $n$ such that
\begin{align*}
E\left( \int_{0}^{\tau_{k}}|Z^{n}_{s}|^{2}\,ds
+(\int_{0}^{\tau_{k}}d|{L}^{n}|_{s})^{2}
+(\int_{0}^{\tau_{k}}|g_{n}(s)|\,ds)^{2}
+(\int_{0}^{\tau_{k}}|f(s,Y^{n}_{s},Z^{n}_{s})|\,ds)^{2}\right)\le
C.
\end{align*}
Similarly, by Lemma \ref{lm4.55}  there exist a stationary
sequence $\{\delta^{4}_{k}\}\subset\mathcal{T}$ and constants
$a_{k}$ such that $L^{-}_{\delta^{4}_{k}}\le a_{k}$, and
furthermore, putting
$\tau'_{k}=\min\{\delta^{1}_{k},\delta^{4}_{k},\delta^{5}_{k}\}$,
where $\delta^{5}_{k}$ is defined as $\delta^{3}_{k}$ but with
$f(s,0,0)$ replaced by $\bar{f}(s)$, we conclude from Propositions
\ref{prop.gen} and \ref{prop4.3} that
\begin{align*}
E\left(\int_{0}^{\tau'_{k}}|Z_{s}|^{2}\,ds
+(\int_{0}^{\tau'_{k}}d|{L}|_{s})^{2}+
(\int_{0}^{\tau'_{k}}|g(s)|\,ds)^{2}
+(\int_{0}^{\tau'_{k}}|\bar{f}(s)|\,ds)^{2}\right)<\infty.
\end{align*}

Step 3. By Step 1 and Step 2 we may assume that there exists $C$
not depending on $n$ such that
\begin{align}
\label{eq.p1}
E\left(\int_{0}^{T}|Z^{n}_{s}|^{2}\,ds
+(\int_{0}^{T}d|{L}^{n}|_{s})^{2}+(\int_{0}^{T} |g_{n}(s)|\,ds)^{2}
+(\int_{0}^{T}|f(s,Y^{n}_{s},Z^{n}_{s})|\,ds)^{2}\right)\le
C.
\end{align}
and $\underline{y}, \overline{y}$ are bounded, $h,\bar{f},g,
f(\cdot,0,0),  f^{+}(\cdot,\overline{y}_{\cdot},0),
f^{-}(\cdot,\underline{y}_{\cdot},0)\in \mathbb{L}^{1,2}(\FF),$
$Z\in M^{2}, L, C\in \mathcal{V}^{2}$. We also may assume that
$\mu\le 0$. By (d), (f), (H4), (H5) and the Lebesgue dominated
convergence theorem,
\begin{equation}
\label{m3.811} \int_{0}^{T}|f(s,Y^{n}_{s},0)-f(s,Y_{s},0)|\,ds
\rightarrow0
\end{equation}
By (d) and (H3),
\[
-f^{-}(s,\underline{y}_{s},0)\le f(s,Y^{n}_{s},0)
\le f^{+}(s,\bar{y}_{s},0),
\]
Therefore using (\ref{m3.811}) and the Lebesgue dominated
convergence theorem we get
\begin{equation}
\label{AABB}
\lim_{n\rightarrow+\infty}E\left(\int_{0}^{T}
|f(s,Y^{n}_{s},0) -f(s,Y_{s},0)|\,ds\right)^{2}=0.
\end{equation}
Observe that by (H2) and (\ref{eq.p1}), $f(s,Y^{n}_{s},Z_{s}^{n})
=f(s,Y^{n}_{s},0)+\eta^{n}(s)$ for some $\{\eta^{n}\}\subset
\mathbb{L}^{2}(\FF)$ such that
\begin{equation}
\label{eq.ppp2}
\sup_{n\ge 0} E\into |\eta^{n}(s)|^{2}\,ds<\infty.
\end{equation}
We have
\begin{align*}
&E\into|f(s,Y^{n}_{s},Z^{n}_{s})-\bar{f}(s)||Y^{n}_{s}-Y_{s}|\,ds\\
&\qquad\le E\into|f(s,Y^{n}_{s},0)-f(s,Y_{s},0)|
|Y^{n}_{s}-Y_{s}|\,ds+E\into|f(s,Y_{s},0)||Y^{n}_{s}-Y_{s}|\,ds\\
&\qquad\quad+ (E\into |\eta^{n}(s)|^{2}\,ds)^{1/2}(E\into
|Y^{n}_{s}-Y_{s}|^{2}\,ds)^{1/2}
+E\into|\bar{f}(s)||Y^{n}_{s}-Y_{s}|\,ds.
\end{align*}
Therefore from (d), (f), (\ref{AABB}), (\ref{eq.ppp2}) and the
fact that  $\underline{y},\overline{y}$ are bounded we get
\begin{equation}
\label{eq.p3} E\into
|f(s,Y^{n}_{s},Z^{n}_{s})-\bar{f}(s)||Y^{n}_{s}-Y_{s}|\,ds\rightarrow0.
\end{equation}
Since $\underline{y},\overline{y}$ are bounded and
$g\in\mathbb{L}^{1}(\FF)$, it follows from (d)--(f) that
\begin{equation}
\label{eq.p4}
E\into |g(s)-g_{n}(s)||Y^{n}_{s}-Y_{s}|\,ds
\le E\into (|h(s)|+|g(s)|)|Y^{n}_{s}-Y_{s}|\,ds\rightarrow 0.
\end{equation}
Let $\sigma,\tau\in\mathcal{T}$, $\sigma\le\tau$. Then by  It\^o's
formula and (c),
\begin{align}
\label{eq.cern} E\int_{\sigma}^{\tau}|Z_{s}-Z^{n}_{s}|^{2}\,ds
&\le E|Y_{\tau}-Y^{n}_{\sigma}|^{2}
+2\int_{\sigma}^{\tau}|g(s)-g^{n}(s)||Y_{s}-Y^{n}_{s}|\,ds\nonumber\\
&\quad+E\int_{\sigma}^{\tau}|f(s,Y_{s},Z_{s})-f(s,Y^{n}_{s},Z^{n}_{s})|
|Y_{s}-Y^{n}_{s}|\,ds\nonumber \\
&\quad +2\int_{\sigma}^{\tau}(Y_{s}-Y^{n}_{s})\,d(L_{s}-L^{n}_{s})
+2\sum_{\sigma<t\le\tau} |\Delta C_{t}|^{2}.
\end{align}
Therefore from (b), (d), (f), (\ref{eq.p3}), (\ref{eq.p4}) and
boundedness  of $\underline{y}, \overline{y}$ it may be concluded
that for any $N\in\mathbb{N}$ and any
$\sigma_{1},\tau_{1},...,\sigma_{N},\tau_{N}\in\mathcal{T}$ such
that $\sigma_{k}\le\tau_{k}$ for $k=1,\dots,N$ we have
\begin{equation}
\label{eq.p5}
\limsup_{n\rightarrow +\infty} \sum_{k=1}^{N}
E\int_{\sigma_{k}}^{\tau_{k}}|Z_{s}-Z^{n}_{s}|^{2}\,ds
\le \sum_{k=1}^{N}E\sum_{\sigma_{k}<t\le \tau_{k}}|\Delta C_{t}|^{2}.
\end{equation}
It suffices now to repeat arguments following Eq. 2.10 in
\cite[Theorem 2.1]{Peng} to show that $Z^{n}\rightarrow
Z,\,\lambda\otimes P-$a.e.. In view of (\ref{eq.p1}) this implies
that $E\into|Z^{n}_{s}-Z_{s}|\,ds\rightarrow 0$, which when
combined with (\ref{AABB}) and (H2) yields (\ref{eq.pm}). Next, by
(\ref{eq.p1}) and pointwise convergence of $Z^{n}$, $E\into
|Z^{n}_{s}-Z_{s}|^{p}\,ds\rightarrow 0$ for every $p\in (0,2)$,
and moreover, if $|\Delta C_{t}|=0$, $t\in[0,T]$, then by
(\ref{eq.p5}) with $N=1, \sigma_{1}=0,\tau_{1}=T$ it follows that
$E\into |Z^{n}_{s}-Z_{s}|^{2}\,ds\rightarrow 0$. Thus,
(\ref{eq.p0}) is satisfied, because we consider processes
$Y^n,Z^n$ stopped at $\tau_k$ (see Step 1). Finally, if
$\{g_{n}\}, \{L^{n}\}$ satisfy the additional assumptions then
$\intot f(s,Y^{n}_{s},Z^{n}_{s})\,ds\rightarrow
\intot\bar{f}(s)\,ds$ weakly in $\mathbb{L}^{1}(\FF_{T})$, which
when combined with (\ref{eq.pm}) implies that
$\bar{f}(s)=f(s,Y_{s},Z_{s})$ for a.e. $s\in [0,T]$.
\end{dow}
\medskip

We are now ready to prove the main result on monotone convergence
of semimartingales.
\begin{tw}
\label{monbsde} Let \mbox{\rm(H2)--(H5)} be satisfied,
$(Y^{n},Z^{n})\in\DM\otimes M$, $A^{n},K^{n}\in\mathcal{V}^{+}$,
$t\mapsto f(t,Y^{n}_{t},Z^{n}_{t})\in\mathbb{L}^{1}(0,T)$ and
\[
Y^{n}_{t}=Y^{n}_{0}-\intot g_{n}(s)\,ds-\intot
f(s,Y^{n}_{s},Z^{n}_{s})\,ds -\intot dK^{n}_{s}+\intot
dA^{n}_{s}+\intot Z^{n}_{s}\,dB_{s}
\]
for $t\in[0,T]$. Assume that
\begin{enumerate}
\item[\rm(a)]$dA^{n}\le dA^{n+1}$, $n\in\mathbb{N}$,
$\sup_{n\ge0}EA^{n}_{T}<\infty$,
\item[\rm(b)]$\liminf_{n\rightarrow+\infty}
\int_{\sigma}^{\tau}(Y_{s}-Y^{n}_{s})\,d(K^{n}_{s}-A^{n}_{s})\ge
0$ for every $\sigma,\tau\in\mathcal{T}$, $\sigma\le \tau$,
\item[\rm(c)]There exists $C\in\mathcal{V}^{1,+}$ such that
$|\Delta K^{n}_{t}|\le|\Delta C_{t}|$, $t\in [0,T]$,
\item[\rm(d)]There exist processes $\underline{y},
\overline{y}\in \mathcal{V}^{1}+\mathcal{M}_{loc}$ of class
\mbox{\rm(D)} such that
\[
E\into f^{+}(s,\overline{y}_{s},0)\,ds +E\into
f^{-}(s,\underline{y}_{s},0)\,ds<\infty,\quad \overline{y}_{t}\le
Y^{n}_{t} \le\underline{y}_{t},\quad t\in [0,T],
\]
\item[\rm(e)]$E\into |f(s,0,0)|\,ds<\infty$ and $|g_{n}(s)|\le h(s)$ for a.e.
$s\in[0,T]$ for some $h\in\mathbb{L}^{1}(\FF),$
\item[\rm(f)]$Y^{n}_{t}\nearrow Y_{t}$, $t\in [0,T]$.
\end{enumerate}
Then $Y\in\DM$, there exist
$K\in\mathcal{V}^{+},A\in\mathcal{V}^{1,+}, Z\in M,
g\in\mathbb{L}^{1}(\FF)$ such that
\[
Y_{t}=Y_{0}-\intot g(s)\,ds-\intot f(s,Y_{s},Z_{s})\,ds-\intot dK_{s}
+\intot dA_{s}+\intot Z_{s}\,dB_{s},\quad t\in [0,T]
\]
and
\[
Z^{n}\rightarrow Z,\quad \lambda\otimes P\mbox{-a.e.}, \quad \int_{0}^{\cdot}
|f(s,Y^{n}_{s},Z^{n}_{s})-f(s,Y_{s},Z_{s})|\,ds \rightarrow
0\quad\mbox{in ucp}.
\]
Moreover, there exists a stationary sequence
$\{\tau_{k}\}\subset\mathcal{T}$ such that for every $p\in(0,2)$,
\[
E\int_{0}^{\tau_{k}}|Z^{n}_{s}-Z_{s}|^{p}\,ds\rightarrow0,
\]
and if $|\Delta C_{t}|+|\Delta K_{t}|=0$, $t\in [0,T]$ then the
above convergence holds for $p=2$, too.
\end{tw}
\begin{dow}
First of all let us note that if there exists a stationary
sequence $\{\tau_{k}\}\subset\mathcal{T}$ such that assertions of
the theorem  hold on the interval $[0,\tau_{k}]$ for every $k\in
\mathbb{N}$ then they hold on $[0,T]$ as well (see Step 1 of the
proof of Lemma \ref{dombsde}). Let $A_{t}=\sup_{n\ge 1}
A^{n}_{t},\, t\in [0,T]$. By (a) and Lemma \ref{LF1},
$A\in\mathcal{V}^{1,+}$. By Lemma \ref{lm4.55} there exist a
stationary sequence $\{\delta^{3}_{k}\}\subset \mathcal{T}$ and
constants $c_{k}$ such that $A_{\delta^{3}_{k}}\le c_{k},\,
k\in\mathbb{N}$. Let $\delta^{1}_{k}, \delta^{2}_{k}$ be defined
as in the proof of Lemma \ref{dombsde}. Then the sequence
$\{\tau_{k}\}\subset\mathcal{T}$, where
$\tau_{k}=\min\{\delta^{1}_{k},\delta^{2}_{k},\delta^{3}_{k}\}$,
is stationary. By the remark at the beginning of the proof and the
definition of $\tau_{k}$ we may assume that
$\underline{y},\overline{y}$ are bounded, $A\in\mathcal{V}^{2,+}$,
$f(\cdot,0,0),h,f^{+}(\cdot,\overline{y}_{\cdot},0),
f^{-}(\cdot,\underline{y}_{\cdot},0)\in\mathbb{L}^{1,2}(\FF)$ and,
by Propositions \ref{prop.gen} and \ref{prop4.3}, that there
exists $C>0$ not depending on $n$ such that
\begin{equation}
\label{eq.mon1}
E\into |Z^{n}_{s}|^{2}\,ds+E(\into dK^{n}_{s})^{2}
+E(\into|f(s,Y^{n}_{s},Z^{n}_{s})|\,ds)^{2}\le C.
\end{equation}
We aim to apply  Lemma \ref{LF1} to the process $Y^{n}$. To do
this we have to show that up to a subsequence,
$\{f(\cdot,Y^{n}_{\cdot},Z^{n}_{\cdot})\}$ is weakly convergent in
$\mathbb{L}^{1}(\FF)$. Since by (\ref{eq.mon1}) and (H2),
$f(s,Y^{n}_{s},Z^{n}_{s})=f(s,Y^{n}_{s},0)+\eta^{n}(s)$ for some
sequence $\{\eta^{n}\}\subset\mathbb{L}^{2}(\FF)$ such that
$\sup_{n\ge 0}E\into|\eta^{n}(s)|^{2}\,ds<\infty$, it suffices to
prove the desired convergence for the sequence
$\{f(\cdot,Y^{n}_{\cdot},0)\}$. Let us observe that by (c), (e),
(H4), (H5) and the Lebesgue dominated convergence theorem,
\[
\into |f(s,Y^{n}_{s},0)-f(s,Y_{s},0)|\,ds\rightarrow0.
\]
By (c) and (H3),
\[
-f^{-}(s,\underline{y}_{s},0)\le f(s,Y^{n}_{s},0)\le
f^{+}(s,\overline{y}_{s},0).
\]
Since $f^{+}(\cdot,\overline{y}_{\cdot},0)$,
$f^{-}(\cdot,\underline{y}_{\cdot},0)\in\mathbb{L}^{1,2}(\FF)$, it
follows from the above and the Lebesgue dominated convergence
theorem that
\[
E\into |f(s,Y^{n}_{s},0)-f(s,Y_{s},0)|\,ds\rightarrow 0.
\]
Let us denote by $\eta$ the weak limit of $\{\eta^{n}\}$ in
$\mathbb{L}^{2}(\FF)$, by $g$ the weak limit of $\{g_{n}\}$ in
$\mathbb{L}^{1}(\FF)$ and by $Z$ the weak limit of $\{Z^{n}\}$ in
$\mathbb{L}^{2}(\FF)$. Then by Lemma \ref{LF1}, $Y\in\DM$ and
there exists $K\in\mathcal{V}^{+}$ such that
\[
Y_{t}=Y_{0}-\intot g(s)\,ds-\intot f(s,Y_{s},0)\,ds-\intot \eta(s)\,ds
-\intot dK_{s}+\intot dA_{s}+\intot Z_{s}\,dB_{s}.
\]
Applying now Lemma \ref{dombsde} we prove the theorem except for
the last assertion.  Finally, using monotonicity of $\{Y^{n}\},
\{dA^{n}\}$ one can show inequality (\ref{eq.cern}) with $|\Delta
K^{n}_{t}|+|\Delta K_{t}|$ in place of $|\Delta C_{t}|$, which
proves the last assertion.
\end{dow}

\nsubsection{BSDEs with one reflecting barrier} \label{sec5}

In this section we prove existence and uniqueness of solutions of
reflected BSDEs with one irregular barrier and data in
$\mathbb{L}^{p}$ with $p\in [1,2)$. We also prove that the
solutions can be approximated by penalization method and give a  comparison
result.

We will need the following  additional hypotheses.

\begin{enumerate}
\item[(H6)]$L, U$ are  progressively
measurable processes, $L_{t}\le U_{t}$ for a.e. $t\in[0,T]$.
\item[(H7)]There exists a semimartingale
$X\in\HH^{p}$ such that $L_{t}\le X_{t}$ for a.e. $t\in[0,T]$ and
$E(\int_{0}^{T}f^{-}(s,X_{s},0)\,ds)^{p}<\infty$.
\item[(H7*)] There exists a semimartingale $X$ of class (D) such that
$X\in\mathcal{V}^{1}+\MM^{q}_{c}$ for every $q\in(0,1)$, $L_{t}\le
X_{t}$, $t\in[0,T]$ and
$E\int_{0}^{T}f^{-}(s,X_{s},0)\,ds<\infty$.
\end{enumerate}

\begin{df}
We say that a triple $(Y,Z,K)$ of progressively measurable
processes  is a solution of RBSDE$(\xi,f+dV,L)$ if
\begin{enumerate}
\item [\rm(a)]$K\in \mathcal{V}^{+}$,
\item [\rm(b)]$Z\in M$ and the mapping $[0,T]\ni t\rightarrow
f(t,Y_{t},Z_{t})$ belongs to $\mathbb{L}^{1}(0,T),\, P$-a.s.,
\item [\rm(c)]$Y_{t}=\xi+\int_{t}^{T}f(s,Y_{s},Z_{s})\,ds+\intt dV_{s}
+\int_{t}^{T}dK_{s}-\int_{t}^{T}Z_{s}\,dB_{s},\quad t\in [0,T],$
\item [\rm(d)]$L_{t}\le Y_{t}$ for a.e.$\, t\in [0,T]$,
$\int_{0}^{T}(Y_{t-}-\hat{L}_{t-})\,dK_{t}=0$ for every
$\hat{L}\in \mathcal{D}$ such that $L_{t}\le \hat{L}_{t}\le
Y_{t},\, P$-a.s. for a.e. $t\in [0,T]$.
\end{enumerate}
\end{df}

Uniqueness of solutions of RBSDEs follows from the following
comparison results.

\begin{stw}\label{prop4.1}
Assume \mbox{\rm(H2)}. Let $(Y^{i},Z^{i},K^{i})$ be a solution of
RBSDE$(\xi^{i},f^{i}+dV^{1},L^{i})$, $i=1,2$. If
$(Y^{1}-Y^{2})^{+}\in\mathcal{D}^{q}$ for some $q>1$,
$\xi^{1}\le\xi^{2}$, $dV^{1}\le dV^{2}$, $L^{1}_{t}\le L^{2}_{t}$
for a.e. $t\in [0,T]$ and either \mbox{\rm(\ref{eq3.11})} or
\mbox{\rm(\ref{eq3.12})} is satisfied then $Y^{1}_{t}\le
Y^{2}_{t}$, $t\in [0,T]$.
\end{stw}
\begin{dow}
Without loss of generality we may assume that $\mu\le0$. Assume
that (\ref{eq3.11}) is satisfied. By the It\^o-Tanaka formula and
Proposition \ref{prop2.1}, for any $p\in(1,q)$ and
$\tau\in\mathcal{T}$ we have
\begin{align*}
&|(Y^{1}_{t\wedge \tau}-Y^{2}_{t\wedge \tau})^{+}|^{p}
+\frac{p(p-1)}{2}\int_{t\wedge\tau}^{\tau}\mathbf{1}_{\{Y^{1}_{s}
\neq Y^{2}_{s}\}}|(Y^{1}_{s}-Y^{2}_{s})^{+}|^{p-2}
|Z^{1}_{s}-Z^{2}_{s}|^{2}\,ds\nonumber\\
&\quad\le|(Y^{1}_{\tau}-Y^{2}_{\tau})^{+}|^{p}
+p\int_{t\wedge\tau}^{\tau}|(Y^{1}_{s}-Y^{2}_{s})^{+}|^{p-1}
(f^{1}(s,Y^{1}_{s},Z^{1}_{s})-f^{2}(s,Y^{2}_{s},Z^{2}_{s}))\,ds\nonumber\\
&\qquad+p\int_{t\wedge\tau}^{\tau}
|(Y^{1}_{s}-Y^{2}_{s})^{+}|^{p-1}\,d(V^{1}_{s}-V^{2}_{s})
+p\int_{t\wedge\tau}^{\tau}|(Y^{1}_{s-}-Y^{2}_{s-})^{+}|^{p-1 }\,
(dK^{1}_{s}-dK^{2}_{s})\nonumber\\
&\qquad-p\int_{t\wedge\tau}^{\tau}|(Y^{1}_{s}-Y^{2}_{s})^{+}|^{p-1}
(Z^{1}_{s}-Z^{2}_{s})\,dB_{s}.
\end{align*}
Since $L^1_t\le L^2_t$ for a.e. $t\in[0,T]$, $L^{1}_{t}\le
Y^{1}_{t}\wedge Y^{2}_{t}\le Y^{1}_{t}$ for a.e. $t\in [0,T]$. By
monotonicity of the function $x\mapsto p|x|^{p-1}\hat{x}$ and
property (d) of the definition of a solution of the reflected
BSDE,
\begin{align*}
&\int_{t\wedge\tau}^{\tau}|(Y^{1}_{s-}-Y^{2}_{s-})^{+}|^{p-1}\,
(dK^{1}_{s}-dK^{2}_{s})\le\int_{t\wedge\tau}^{\tau}
|(Y^{1}_{s-}-Y^{2}_{s-})^{+}|^{p-1}\,dK^{1}_{s}\\
&\qquad=\int_{t\wedge\tau}^{\tau}
|(Y^{1}_{s-}-Y^{2}_{s-})^{+}|^{p-1}\,\hat{\mbox{sgn}}
[(Y^{1}_{s-}-Y^{2}_{s-})^{+}]\,dK^{1}_{s}\\
&\qquad\le \int_{t\wedge\tau}^{\tau}|(Y^{1}_{s-}-Y^{1}_{s-}\wedge
Y^{2}_{s-})^{+}|^{p-1}\,\hat{\mbox{sgn}}[(Y^{1}_{s-}-Y^{1}_{s-}
\wedge Y^{2}_{s-})^{+}]\,dK^{1}_{s}=0,
\end{align*}
the last equality being a consequence of the fact that
\begin{align*}
&\int_{t\wedge\tau}^{\tau}\hat{\mbox{sgn}}
[(Y^{1}_{s-}-Y^{2}_{s-})^{+}]\,dK^{1}_{s}\\
& \qquad=\int_{t\wedge\tau}^{\tau} \mathbf{1}_{\{Y^{1}_{s-}\neq
Y^{1}_{s-}\wedge Y^{2}_{s-}\}} |Y^{1}_{s}-Y^{1}_{s-}\wedge
Y^{2}_{s-}|^{-1} d(\int_{0}^{s} (Y^{1}_{r-}-Y^{1}_{r-}\wedge
Y^{2}_{r-})\,dK^{1}_{r}).
\end{align*}
The rest of the proof runs as the proof of Proposition
\ref{prop3.4} (see the reasoning following (\ref{eq3.13})).
\end{dow}

\begin{wn}
\label{cor.how} Assume \mbox{\rm(H2), (Z)}. Let   $(Y^{i},Z^{i},
K^{i})$, be a solution of RBSDE$(\xi^{i},f^{i}+dV^{i}, L^{i})$
such that $(Y^{i},Z^{i})\in
\mathbb{L}^{q}(\FF)\otimes\mathbb{L}^{q}(\FF)$ for some
$q>\alpha$, $i=1,2$. If $(Y^{1}-Y^{2})^{+}$ is of class (D),
$\xi^{1}\le\xi^{2}$, $dV^{1}\le dV^{2}$, $L^{1}_{t}\le L^{2}_{t}$
for a.e. $t\in [0,T]$ and \mbox{\rm(\ref{eq3.11})} or
\mbox{\rm(\ref{eq3.12})} is satisfied then $Y^{1}_{t}\le
Y^{2}_{t}$, $t\in [0,T]$.
\end{wn}
\begin{dow}
As in the proof of Proposition \ref{prop4.1} one can reduce the
problem to the case of nonreflected BSDEs and then use Corollary
\ref{cor3.10} to conclude the result.
\end{dow}

\begin{uw}
\label{uwnii} Since the proof of Proposition \ref{prop4.1} and
Corollary \ref{cor.how} is based on the proof of Proposition
\ref{prop3.4}, arguments from Remark \ref{uwni} show that if $f$
does not depend on $z$ then  assumption (Z) and the assumptions
that $(Y^{i},Z^{i})\in
\mathbb{L}^{q}(\FF)\otimes\mathbb{L}^{q}(\FF)$ for some $q>\alpha$
 are superfluous in Corollary \ref{cor.how}.
\end{uw}

For an arbitrary sequence $\{x_{n}\}$ of elements of some linear
space we set
\[
\sigma_{k}(\{x_{n}\})=\frac{x_{1}+x_{2}+\ldots+x_{k}}{k}\,,\quad
k\ge 0.
\]

\begin{tw}
\label{th4.1} Let assumptions  \mbox{\rm(H2), (H3)} hold. Then
there exists at most one solution $(Y,Z,K)$ of RBSDE$(\xi,f+dV,L)$
such that $Y\in \mathcal{D}^{p}$ for some $p>1$.
\end{tw}
\begin{dow}
Follows immediately from Proposition \ref{prop4.1}.
\end{dow}

\begin{tw}\label{stw4}
Let assumptions  \mbox{\rm(H2), (H3), (Z)} hold. Then there exists
at most one solution $(Y,Z,K)$ of RBSDE$(\xi,f+dV,L)$ such that
$Y$ is of class \mbox{\rm(D)} and $Z\in\bigcup_{\beta>\alpha}
M^{\beta}$.
\end{tw}
\begin{dow}
Follows immediately from Corollary \ref{cor.how}.
\end{dow}
\medskip

We prove existence of solutions separately for data in $\BL^p$
with $p>1$ and for data in $\BL^1$.

\begin{tw}\label{th.main1}
Let $p>1$.
\begin{enumerate}
\item[\rm(i)] Assume  \mbox{\rm(H1)--(H6)}.
Then there exists a solution $(Y,Z,K)\in\mathcal{D}^{p}\otimes
M^{p}\otimes\mathcal{V}^{+,p}$ of RBSDE$(\xi,f+dV,L)$ iff
\mbox{\rm(H7)} is satisfied.
\item[\rm(ii)]Assume \mbox{\rm(H1)--(H7)}. For
$n\in\mathbb{N}$ let $(Y^{n},Z^{n})\in \DM^{p}\otimes M^{p}$ be a solution
of the BSDE
\begin{align}
\label{PPP}
Y^{n}_{t}&=\xi_{n}+\int_{t}^{T}f(s,Y^{n}_{s},Z^{n}_{s})\,ds+\intt
dV_{s}\nonumber\\
&\quad+n\int_{t}^{T}(Y^{n}_{s}-L_{s})^{-}\,ds
-\int_{t}^{T}Z^{n}_{s}\,dB_{s},\quad t\in[0,T],
\end{align}
and let $\xi_{1}\in\mathbb{L}^{p}(\FF_{T})$, $\xi_{n}\nearrow
\xi$. Then
\begin{equation}
\label{eq5.07} Y^{n}_{t}\nearrow Y_{t},\quad t\in [0,T],\quad
Z^{n}\rightarrow Z,\,\, \lambda\otimes P\mbox{-a.s.},
\end{equation}
for every $r\in [1,p)$, $q\in [1,2)$,
\begin{equation}
\label{eq5.03}
E(\into|Z^{n}_{s}-Z_{s}|^{q}\,ds)^{r/q}\rightarrow0,
\end{equation}
and for every stopping time $\tau\in\mathcal{T}$,
\begin{equation}
\label{eq5.04} K^{n}_{\tau}\rightarrow K_{\tau}\mbox{ weakly in }
\mathbb{L}^{p}(\FF_{T}),
\end{equation}
where
\begin{equation}
\label{eq5.08} K^{n}_{t}=n\int^t_0(Y^{n}_{s}-L_{s})^{-}\,ds.
\end{equation}
\end{enumerate}
\end{tw}
\begin{dow}
Necessity. Assume that there exists a solution $(Y,Z,K)\in
\DM^{p}\otimes M^{p}\otimes\mathcal{V}^{+,p}$ of
RBSDE$(\xi,f+dV,L)$. Then by Proposition \ref{prop.gen}, $E(\into
|f(s,Y_{s},Z_{s})|\,ds)^{p}<\infty$, which implies that $Y\in
\mathcal{H}^{p}$. Moreover, by the definition of a solution of
RBSDE$(\xi,f,L)$, $Y_{t}\ge L_{t}$ for a.e. $t\in [0,T]$, i.e.
(H7) is satisfied with $X=Y$.

Sufficiency. First of all let us note that by Theorem \ref{th3.1}
there exists a unique solution $(Y^{n},Z^{n})\in\DM^{p}\otimes
M^{p}$ of (\ref{PPP}). Using standard change of variable we may
reduce the proof to the case $\mu\le0$. Therefore in what follows
we assume that $\mu\le0$.  For simplicity of notation we will also
assume that $\xi_{n}=\xi$, $n\ge 0$. The proof in the general case
only requires some obvious changes.

Step 1. We will show that there exists a supersolution
$(\bar{X},\bar{Z})\in\DM^{p}\otimes M^{p}$ of BSDE$(\xi,f+dV)$
such that
\begin{equation}
\label{eq5.01} X_{t}\le\bar{X}_{t},\quad Y^{n}_{t}\le
\bar{X}_{t},\quad t\in [0,T]
\end{equation}
and
\begin{equation}
\label{eq5.02}
E(\into|f(s,\bar{X}_{s},\bar{Z}_{s})|\,ds)^{p}<\infty.
\end{equation}
Since $X\in\mathcal{H}^{p}$ and the Brownian filtration has the
representation property, there exist $C\in\mathcal{V}^{p},H\in
M^{p}$ such that
\[
X_{t}=X_{T}-\intt dC_{s}-\intt H_{s}\,dB_{s},\quad t\in [0,T].
\]
The above equation can be rewritten in the form
\begin{align*}
X_{t}&=X_{T}+\intt f(s,X_{s},H_{s})\,ds+\intt dV_{s}-\intt
(f^{+}(s,X_{s},H_{s})\,ds+dC^{+}_{s}+dV^{+}_{s})\\
&\quad+\intt (f^{-}(s,X_{s},H_{s})\,ds
+dC^{-}_{s}+dV_{s}^{-})-\intt H_{s}\,dB_{s},\quad t\in [0,T].
\end{align*}
By (H7) and (H2), $E(\into f^{-}(s,X_{s},H_{s})\,ds)^{p}<\infty$.
By Theorem \ref{th3.1} there exists a solution
$(\bar{X},\bar{Z})\in\DM^{p}\otimes M^{p}$ of the BSDE
\begin{align*}
\bar{X}_{t}&=X_{T}\vee\xi+\intt f(s,\bar{X}_{s},\bar{Z}_{s})\,ds
+\intt dV_{s}\\
&\quad+\intt (f^{-}(s,X_{s},H_{s})\,ds
+dC^{-}_{s}+dV^{-}_{s})-\intt \bar{Z}_{s}\,dB_{s},\quad t\in
[0,T].
\end{align*}
Moreover, by Proposition \ref{prop.gen},
\[
E\left(\into |f(s,\bar{X}_{s},\bar{Z}_{s})|\,ds
\right)^{p}<\infty,
\]
and by Proposition \ref{prop3.4},
\[
X_{t}\le\bar{X}_{t},\quad t\in [0,T].
\]
Of course, $(\bar{X},\bar{Z})$ is a supersolution of
BSDE$(\xi,f+dV)$. Since $L_{t}\le X_{t}$ for a.e. $t\in [0,T]$, it
follows from the last estimate that $L_{t}\le \bar{X}_{t}$ for
a.e. $t\in [0,T]$. Therefore
\begin{align*}
\bar{X}_{t}&=X_{T}\vee\xi+\intt
f(s,\bar{X}_{s},\bar{Z}_{s})\,ds+\intt dV_{s}+\intt
(f^{-}(s,X_{s},H_{s})\,ds+dC^{-}_{s}+dV^{-}_{s})\\
&\quad +n\intt(\bar{X}_{s}-L_{s})^{-}\,ds-\intt
\bar{Z}_{s}\,dB_{s},\quad t\in [0,T].
\end{align*}
Hence, by Proposition \ref{prop3.4}, $Y^{n}_{t}\le\bar{X}_{t}$,
$t\in[0,T]$, which completes the proof of Step 1.

Step 2.  By Proposition \ref{prop3.4}, $Y^{n}_{t}\le Y^{n+1}_{t}$,
$t\in [0,T]$ for every $n\in\mathbb{N}$. Therefore setting
$Y_{t}=\sup_{n\ge1}Y^{n}_{t}$, $t\in [0,T]$, we have
\begin{equation}
\label{m0}
Y^{n}_{t}\nearrow Y_{t},\quad t\in [0,T].
\end{equation}
Since the assumptions of Theorem \ref{monbsde} are satisfied with
$\underline{y}=\bar{X}$ and $\overline{y}=Y^{1}$,  there exists
$Z\in M$ such that
\begin{equation}
\label{eq5.9} Z^{n}\rightarrow Z,\quad\lambda\otimes P\mbox{-a.e.}
\end{equation}
and there exists $K\in\mathcal{V}^{+}$ such that
\[
Y_{t}=\xi+\int_{t}^{T}f(s,Y_{s},Z_{s})\,ds
+\intt dV_{s}+\int_{t}^{T}dK_{s}-\int_{t}^{T}Z_{s}\,dB_{s},\quad t\in [0,T].
\]
By (\ref{eq5.01}), (\ref{eq5.02}) and Propositions \ref{prop4.3}
and \ref{prop.gen}, there exists $C$ not depending on $n$ such
that
\begin{equation}
\label{m2}
E\sup_{t\le T}|Y^{n}_{t}|^{p}+E(\into|Z^{n}_{s}|^{2}\,ds)^{p/2}
+E(\into dK^{n}_{s})^{p}+E(\into |f(s,Y^{n}_{s},Z^{n}_{s})|\,ds)^{p}\le C.
\end{equation}
From (\ref{m0})--(\ref{m2}) it follows that $Y\in\DM^{p}$, $Z\in
M^{p}$, $f(\cdot,Y,Z)\in\mathbb{L}^{1,p}(\FF)$,
$K\in\mathcal{V}^{p,+}$ and (\ref{eq5.03}), (\ref{eq5.04}) hold
true.

Step 3. In the last part of the proof we show that $Y$ majorizes $L$ and
$K$ satisfies the minimality condition. Let $C$ be the constant
appearing on the right-hand side of (\ref{m2}). From (\ref{PPP})
one can easily deduce that there exists $c_{p}$ depending only on
$p$ such that
\[
E(\into (Y^{n}_{s}-L_{s})^{-}\,ds)^{p}\le c_{p}Cn^{-p},
\]
which when combined with (\ref{m0}) implies that
\[
L_{t}\le Y_{t}\quad\mbox{for a.e. }t\in [0,T].
\]
From (\ref{m2}) and the fact that the space
$\mathbb{L}^{2,p}(\FF)$ has the Banach-Saks property we conclude
that there exists a subsequence (still denoted by $\{n\}$) such
that
\[
\sigma_{n}(\{\int_{0}^{\cdot}Z^{n}_{s}\,dB_{s}\})
\rightarrow\int_{0}^{\cdot}Z_{s}\,dB_{s}\quad\mbox{in ucp}.
\]
By Theorem \ref{monbsde},
\begin{equation}
\label{m4}
\int_{0}^{\cdot}|f(s,Y^{n}_{s},Z^{n}_{s})-f(s,Y_{s},Z_{s})|\,ds
\rightarrow 0\quad\mbox{in ucp}.
\end{equation}
By (\ref{m0}) and (\ref{m4}),
\[
\sigma_{n}(\{\int_{0}^{\cdot}|f(s,Y^{n}_{s},Z^{n}_{s})
-f(s,Y_{s},Z_{s})|\,ds\})\rightarrow 0\mbox{ in ucp},\quad
\sigma_{n}(\{Y^{n}_{t}\})\rightarrow Y_{t},\,\,t\in [0,T],
\]
and hence
\begin{equation}
\label{m5}
\sigma_{n}(\{K^{n}_{t}\})\rightarrow K_{t},\quad t\in [0,T].
\end{equation}
Let $\hat{L}\in\DM$ be an arbitrary process such that $L_{t}\le
\hat{L}_{t}\le Y_{t}$ for a.e. $t\in [0,T]$ and let $\{\tau_{k}\}$
be an increasing sequence of succesive jumps of process $V,K$
(with the convention that $\tau_{0}\equiv 0$). Since $Y^{n}$ is
continuous on $(\tau_{k},\tau_{k+1})$,  it follows from Dini's
theorem that for every $k\in\mathbb{N}$, $Y^{n}(\omega)\rightarrow
Y(\omega)$ uniformly on compact sets in
$(\tau_{k}(\omega),\tau_{k+1}(\omega))$ for a.e.
$\omega\in\Omega$. Therefore by (\ref{m5}) and Helly's theorem,
\begin{equation}
\label{m6}
\int_{(\tau_{k},\tau_{k+1})}(Y^{[n/l]}_{s}-\hat{L}_{s})\,
d\sigma_{n}(\{K^{n}_{s}\})\rightarrow \int_{(\tau_{k},\tau_{k+1})}
(Y_{s-}-\hat{L}_{s-})\,dK_{s}
\end{equation}
for every $k,l\in\mathbb{N}$. On the other hand,
\begin{align*}
&(Y^{[n/l]}_{s}-\hat{L}_{s-})d\sigma_{n}(\{K^{n}_{s}\})
=\frac{1}{n}\sum_{k=1}^{n}k(Y^{[n/l]}_{s}
-\hat{L}_{s-})(Y^{k}_{s}-L_{s})^{-}\,ds\\
&\quad=\frac{1}{n}\sum_{k<[n/l]}k(Y^{[n/l]}_{s}-\hat{L}_{s-})
(Y^{k}_{s}-L_{s})^{-}\,ds+\frac{1}{n}\sum_{[n/l]\le k\le n}
k(Y^{[n/l]}_{s} -\hat{L}_{s-})(Y^{k}_{s}-L_{s})^{-}\,ds\\
&\quad \le \frac{2}{n}|Y^{*}_{T}|
\sum_{k<[n/l]}k(Y^{k}_{s}-L_{s})^{-}\,ds +\frac{1}{n}\sum_{[n/l]
\le
k\le n}k(Y^{k}_{s}-\hat{L}_{s})(Y^{k}_{s}-\hat{L}_{s})^{-}\,ds\\
&\quad\le\frac{2}{n}|Y^{*}_{T}|\sum_{k<[n/l]}k(Y^{k}_{s}-L_{s})^{-}\,ds.
\end{align*}
Hence
\[
(Y^{[n/l]}_{s}-\hat{L}_{s-})d\sigma_{n}(\{K^{n}_{s}\}) \le
2|Y^{*}_{T}|\frac{[n/l]}{n}d\sigma_{[n/l]}(\{K^{n}_{s}\}).
\]
By the above, (\ref{m6}) and  Helly's theorem,
\[
0\le \int_{(\tau_{k},\tau_{k+1})} (Y_{s-}-\hat{L}_{s-})\,dK_{s}
\le 2|Y^{*}_{T}|\frac{1}{l}\into dK_{s}
\]
for $k,l\in\mathbb{N}$, which implies that for every
$k\in\mathbb{N}$,
\begin{equation}
\label{m7}
\int_{(\tau_{k},\tau_{k+1})} (Y_{s-}-\hat{L}_{s-})\,dK_{s}=0.
\end{equation}
What is left is to show that
\begin{equation}
\label{m8}
\sum_{0<t\le T}(Y_{t-}-\hat{L}_{t-})\Delta K_{t}=0.
\end{equation}
But (\ref{m8}) is an immediate consequence of Lemma \ref{lm4.6}.
Indeed, if $\Delta K_{t}>0$ then $-\Delta Y_{t}-\Delta
V_{t}=\Delta K_{t}>0$, which implies that $Y_{t}+\Delta
V_{t}<Y_{t-}$\,. By the last inequality and (\ref{eq4.01}),
$Y_{t-}=\hat{L}_{t-}$. Thus, for every $t\in (0,T]$, if $\Delta
K_{t}>0$ then $Y_{t-}=\hat{L}_{t-}$\,, which forces (\ref{m8}). By
(\ref{m7}) and (\ref{m8}),
\begin{equation}
\label{eq5.06} \into(Y_{t-}-\hat{L}_{t-})\,dK_{t}=0.
\end{equation}
Since (\ref{eq5.06}) holds true for any process $\hat{L}\in\DM$
such that $L_{t}\le\hat{L}_{t}\le Y_{t}$ for a.e. $t\in[0,T]$, the
process $K$ satisfies the minimality condition.
\end{dow}

\begin{tw}\label{th.main2}
Let $p=1$.
\begin{enumerate}
\item[\rm(i)] Assume \mbox{\rm(H1)--(H6), (Z)}.
Then there exists a solution $(Y,Z,K)$ of RBSDE$(\xi,f+dV,L)$ such
that $(Y,Z,K)\in \mathcal{D}^{q}\otimes
M^{q}\otimes\mathcal{V}^{+,1}$ for every $q\in (0,1)$ and $Y$ is
of class \mbox{\rm(D)} iff \mbox{\rm(H7*)} is satisfied.
\item[\rm(ii)]Assume \mbox{\rm(H1)--(H6), (H7*), (Z)}. For
$n\in\mathbb{N}$ let $(Y^{n},Z^{n})$ be a solution of
\mbox{\rm(\ref{PPP})} with $\xi_n$ such that
$\xi_{1}\in\mathbb{L}^{p}(\FF_{T})$, $\xi_{n}\nearrow\xi$,
$(Y^{n},Z^{n})\in\bigcap_{q<1}\DM^{q}\otimes M^{p}$ and $Y^{n}$ is
of class \mbox{\rm(D)}. Then \mbox{\rm(\ref{eq5.07})} holds true
and there exists a stationary sequence $\{\tau_{k}\}$ of stopping
times such that for any $q\in[1,2)$, $p>1$,
\begin{equation}
\label{eq5.28}
E(\int_{0}^{\tau_{k}}|Z^{n}_{s}-Z_{s}|^{q}\,ds)^{p}\rightarrow0
\end{equation}
and for every stopping time $\tau\in\mathcal{T}$,
\begin{equation}
\label{eq5.29} K^{n}_{\tau\wedge\tau_{k}}\rightarrow
K_{\tau\wedge\tau_{k}}\quad\mbox{weakly in }
\mathbb{L}^{p}(\FF_{T}),
\end{equation}
where $K^n$ is defined by \mbox{\rm(\ref{eq5.08})}.
\end{enumerate}
\end{tw}
\begin{dow}
(i) Necessity. Follows from  Proposition \ref{prop.gen} by the
same method as in the proof of Theorem \ref{th.main1}.

Sufficiency. In much the same way as in the proof of Theorem
\ref{th.main1} one can show that there exists a supersolution
$(\bar{X},\bar{Z})$ of BSDE$(\xi,f+dV)$ such that
$(\bar{X},\bar{Z})\in\DM^{q}\otimes M^{q}$ for $q\in (0,1)$,
$\bar{X}$ is of class (D), (\ref{eq5.01}), (\ref{eq5.02}) with
$p=1$ are satisfied and $Y^{n}_{t}\le Y^{n+1}_{t}$, $t\in [0,T]$.
The only difference in the proof lies in the fact that we replace
the space $\DM^{p}\otimes M^{p}$ by the space of processes
$(Y,Z)\in\DM^{q}\otimes M^{q}$, $q\in (0,1)$, such that $Y$ is of
class (D), we replace $\mathcal{V}^{+,p}$ by $\mathcal{V}^{+,1}$
and we use Theorem \ref{th3.2} instead of Theorem \ref{th3.1}, and
Corollary \ref{cor3.10} instead of Proposition \ref{prop3.4}. By
Lemma \ref{lm4.55}, there exist a stationary sequence
$\{\delta^{1}_{k}\}\subset \mathcal{T}$ and constants $c_{k}$ such
that $\bar{X}^{*}_{\delta^{1}_{k}}+|V|_{\delta^{1}_{k}}
+Y^{1,*}_{\delta^{1}_{k}}\le c_{k}$. Let
$\tau_{k}=\delta^{1}_{k}\wedge\delta^{2}_{k}$, where
\[
\delta^{2}_{k}=\inf\{t\in [0,T]; \intot|f(s,0,0)|\,ds +\intot
|f(s,\bar{X}_{s},\bar{Z}_{s})|\,ds >k\}\wedge T.
\]
The sequence $\{\tau_{k}\}$ is stationary. Observe that the data
$(Y^{n}_{\tau_{k}}, f, V, L)$ satisfy the assumptions of Theorem
\ref{th.main1} on each interval $[0,\tau_{k}]$ for every $p>1$.
Using Theorem \ref{th.main1} and stationarity of the sequence
$\{\tau_{k}\}$ shows that there exists a triple
$(Y,Z,K)\in\DM\otimes M\otimes\mathcal{V}^{+}$ such that $t\mapsto
f(t,Y_{t},Z_{t})\in \mathbb{L}^{1}(0,T)$, (\ref{eq5.07}),
(\ref{eq5.28}), (\ref{eq5.29}) hold true,
\[
Y_{t}=\xi+\intt f(s,Y_{s},Z_{s})\,ds +\intt dV_{s}+\intt
dK_{s}-\intt Z_{s}\,dB_{s},\quad t\in [0,T]
\]
and $Y_{t}\ge L_{t}$ for a.e. $t\in [0,T]$,  $\into
(Y_{t-}-\hat{L}_{t-})\,dK_{t}=0$ for every $\hat{L}\in\DM$ such
that $L_{t}\le\hat{L}_{t}\le Y_{t}$ for a.e. $t\in [0,T]$.  The
proof is completed by showing integrability properties of $Y,Z,K$.
Integrability of $Y$ follows from (\ref{eq5.01}), monotonicity of
the sequence $\{Y^{n}\}$ and the fact that the processes
$Y^{1},\bar{X}$ belong to $\DM^{q}$ for $q\in (0,1)$ and are of
class (D). By integrability of $Y$ and Lemma \ref{lm4.2},  $Z\in
M^{q},\, q\in (0,1)$. Let us set
\begin{equation*}
\tau_k=\inf\{t\in [0,T];\int_0^t|Z_{s}|^{2}\,ds\ge k\}\wedge T.
\end{equation*}
Then
\[
K_{\tau_{k}}=Y_{0}-Y_{\tau_{k}}
-\int_{0}^{\tau_{k}}f(s,Y_{s},Z_{s})\,ds
+\int_{0}^{\tau_{k}}\,dB_{s}.
\]
Since $Y$  is of class (D), using Fatou's lemma, (H2), (Z) and the
fact that $Y_{t}\le\bar X_{t}$, $t\in [0,T]$, we conclude from the
above equality that
\[
EK_{T}\le EY^{+}_{0}+E\xi^{-}+E\int_{0}^{T}f^{-}(s,\bar
X_{s},0)\,ds+\gamma
E\int_{0}^{T}(g_{s}+|Y_{s}|+|Z_{s}|)^{\alpha}\,ds.
\]
Hence $K\in\mathcal{V}^{+,1}$ by (\ref{eq5.02}) with $p=1$ and
integrability of $(Y,Z)$.
\end{dow}

\nsubsection{BSDEs with two reflecting barriers} \label{sec6}

In this section we generalize results of Section \ref{sec5} to the
case of BSDEs with two irregular reflecting barriers and data in
$\mathbb{L}^{p}$ with $p\in [1,2)$.

The following natural hypotheses on the barriers generalize the
so-called Mokobodzki condition.

\begin{enumerate}
\item[(H8)]There exists $X\in\HH^{p}$ such that $L_{t}\le X_{t}\le U_{t}$
for a.e. $t\in [0,T]$ and
\[
E(\int_{0}^{T}|f(s,X_{s},0)|\,ds)^{p}<\infty.
\]
\item [(H8*)]There exists a semimartingale $X$ of class (D) such that
$X\in\mathcal{V}^{1}+\MM^{q}_{c}$ for every $q\in (0,1)$,
$L_{t}\le X_{t}\le U_{t}$ for a.e.  $t\in [0,T]$ and
\[
E\int_{0}^{T}|f(s,X_{s},0)|\,ds<\infty.
\]
\end{enumerate}

In the sequel the abbreviation $\underline{\mbox{R}}$BSDE stands
for reflected BSDE with lower obstacle
and $\overline{\mbox{R}}$BSDE stands for reflected BSDE with upper
obstacle.

\begin{df}
We say that a triple $(Y,Z,R)$ of progressively measurable
processes  is a solution of RBSDE$(\xi,f+dV,L,U)$ if
\begin{enumerate}
\item[\rm(a)]$R\in \mathcal{V}$,
\item[\rm(b)]$Z\in M$, the mapping $[0,T]\ni t\mapsto
f(t,Y_{t},Z_{t})$ belongs to $\mathbb{L}^{1}(0,T)$, $P$-a.s.,
\item[\rm(c)]$Y_{t}=\xi+\int_{t}^{T}f(s,Y_{s},Z_{s})\,ds+\intt dV_{s}
+\int_{t}^{T}dR_{s}-\int_{t}^{T}Z_{s}\,dB_{s}$, $t\in[0,T]$,
\item[\rm(d)]$L_{t}\le Y_{t}\le U_{t}$ for a.e. $t\in[0,T]$ and
\begin{equation}
\label{eq6.01} \int_{0}^{T}(Y_{t-}-\hat{L}_{t-})\,dR^{-}_{t}
=\int_{0}^{T}(\check{U}_{t-}-Y_{t-})\,dR^{-}_{t}=0
\end{equation}
for every $\hat{L},\check{U}\in\mathcal{D}$ such that $L_{t}\le
\hat{L}_{t}\le Y_{t}\le\check{U}_{t}\le U_{t}$ for a.e. $t\in
[0,T]$.
\end{enumerate}
\end{df}

\begin{stw}\label{propd4.1}
Assume \mbox{\rm(H2)}. Let $(Y^{i},Z^{i},R^{i})$ be a solution of
RBSDE$(\xi^{i},f^{i}+dV^{i},L^{i},U^{i})$, $i=1,2$. If
$(Y^{1}-Y^{2})^{+}\in\mathcal{D}^{q}$ for some $q>1$,
$\xi^{1}\le\xi^{2}$, $dV^{1}\le dV^{2}$, $L^{1}_{t}\le
L^{2}_{t},\, U^{1}_{t}\le U^{2}_{t}$ for a.e. $t\in [0,T]$ and
either \mbox{\rm(\ref{eq3.11})} or \mbox{\rm(\ref{eq3.12})} is
satisfied then $Y^{1}_{t}\le Y^{2}_{t}$, $t\in [0,T]$.
\end{stw}
\begin{dow}
Without loss of generality we may assume that $\mu\le0$. Let us
fix $p\in(1,q)$ and assume that (\ref{eq3.11}) is satisfied. By
the It\^o-Tanaka formula and  Proposition \ref{prop2.1}, for every
$\tau\in\mathcal{T}$,
\begin{align*}
&|(Y^{1}_{t\wedge \tau}-Y^{2}_{t\wedge \tau})^{+}|^{p}
+\frac{p(p-1)}{2}\int_{t\wedge\tau}^{\tau}\mathbf{1}_{\{Y^{1}_{s}
\neq Y^{2}_{s}\}}|(Y^{1}_{s}-Y^{2}_{s})^{+}|^{p-2}
|Z^{1}_{s}-Z^{2}_{s}|^{2}\,ds\nonumber\\
&\quad\le|(Y^{1}_{\tau}-Y^{2}_{\tau})^{+}|^{p}
+p\int_{t\wedge\tau}^{\tau}|(Y^{1}_{s}-Y^{2}_{s})^{+}|^{p-1}
(f^{1}(s,Y^{1}_{s},Z^{1}_{s})-f^{2}(s,Y^{2}_{s},Z^{2}_{s}))\,ds\nonumber\\
&\qquad+p\int_{t\wedge\tau}^{\tau}|(Y^{1}_{s-}-Y^{2}_{s-})^{+}|^{p-1}
\,(dV^{1}_{s}-dV^{2}_{s})+p\int_{t\wedge\tau}^{\tau}
|(Y^{1}_{s-}-Y^{2}_{s-})^{+}|^{p-1 }\,(dR^{1}_{s}-dR^{2}_{s})\nonumber\\
&\qquad-p\int_{t\wedge\tau}^{\tau}|(Y^{1}_{s}-Y^{2}_{s})^{+}|^{p-1}
(Z^{1}_{s}-Z^{2}_{s})\,dB_{s}.
\end{align*}
Since $L^1_t\le L^2_t$, $U^{1}_{t}\le U^{2}_{t}$  for a.e.
$t\in[0,T]$, $L^{1}_{t}\le Y^{1}_{t}\wedge Y^{2}_{t}\le
Y^{1}_{t}$, $Y^{2}_{t}\le Y^{1}_{t}\vee Y^{2}_{t}\le U^{2}_{t}$
for a.e $t\in [0,T]$. By monotonicity of the function $x\mapsto
p|x|^{p-1}\hat{x}$ and property (d) of the definition of a
solution of the reflected BSDE,
\begin{align*}
&\int_{t\wedge\tau}^{\tau}|(Y^{1}_{s-}-Y^{2}_{s-})^{+}|^{p-1}\,
(dR^{1}_{s}-dR^{2}_{s})\le \int_{t\wedge\tau}^{\tau}
|(Y^{1}_{s-}-Y^{2}_{s-})^{+}|^{p-1}\,d(R^{1,+}_{s}+R^{2,-}_{s})\\
&\qquad=\int_{t\wedge\tau}^{\tau}|(Y^{1}_{s-}-Y^{2}_{s-})^{+}|^{p-1}\,
\hat{\mbox{sgn}}[(Y^{1}_{s-}-Y^{2}_{s-})^{+}]\,dR^{1,+}_{s}\\
&\qquad\quad+\int_{t\wedge\tau}^{\tau}|(Y^{1}_{s-}-Y^{2}_{s-})^{+}|^{p-1}
\,\hat{\mbox{sgn}}[(Y^{1}_{s-}-Y^{2}_{s-})^{+}]\,dR^{2,-}_{s}\\
&\qquad\le\int_{t\wedge\tau}^{\tau}|(Y^{1}_{s-}-Y^{1}_{s-}\wedge
Y^{2}_{s-})^{+}|^{p-1}\,\hat{\mbox{sgn}}[(Y^{1}_{s-}-Y^{1}_{s-}
\wedge Y^{2}_{s-})^{+}]\,dR^{1,+}_{s}\\
&\qquad\quad+ \int_{t\wedge\tau}^{\tau}|(Y^{1}_{s-}\vee
Y^{2}_{s-}-
Y^{2}_{s-})^{+}|^{p-1}\,\hat{\mbox{sgn}}[(Y^{1}_{s-}\vee
Y^{2}_{s-}-Y^{2}_{s-})^{+}]\,dR^{2,-}_{s}=0,
\end{align*}
the last equality being a consequence of the fact that
\begin{align*}
&\int_{t\wedge\tau}^{\tau}\hat{\mbox{sgn}}[(Y^{1}_{s-}-Y^{1}_{s}\wedge
Y^{2}_{s-})^{+}]\,dR^{1,+}_{s}\\
&\qquad =\int_{t\wedge\tau}^{\tau} \mathbf{1}_{\{Y^{1}_{s-}\neq
Y^{1}_{s-}\wedge Y^{2}_{s-}\}} |Y^{1}_{s}-Y^{1}_{s-}\wedge
Y^{2}_{s-}|^{-1} d(\int_{0}^{s} (Y^{1}_{r-}-Y^{1}_{r-}\wedge
Y^{2}_{r-})\,dR^{1,+}_{r})
\end{align*}
and
\begin{align*}
&\int_{t\wedge\tau}^{\tau} \hat{\mbox{sgn}}[(Y^{1}_{s-}\vee
Y^{2}_{s-}-Y^{2}_{s-})^{+}]\,dR^{2,-}_{s}\\
&\qquad =\int_{t\wedge\tau}^{\tau} \mathbf{1}_{\{Y^{1}_{s-}\vee
Y^{2}_{s-}\neq Y^{2}_{s-}\}} |Y^{1}_{s-}\vee Y^{2}_{s-}-
Y^{2}_{s-}|^{-1} d(\int_{0}^{s} (Y^{1}_{r-}\vee
Y^{2}_{r-}-Y^{2}_{r-})\,dR^{2,-}_{r}).
\end{align*}
The rest of the proof runs as the proof of Proposition
\ref{prop3.4} (see the reasoning following (\ref{eq3.13})).
\end{dow}

\begin{wn}
\label{cor.ad21} Assume that \mbox{\rm(H2), (Z)} are satisfied.
For $i=1,2$ let $(Y^{i},Z^{i}, R^{i})$ be a solution of
RBSDE$(\xi^{i},f^{i}+dV^{i},L^{i}, U^{i})$ such that
$(Y^{i},Z^{i})\in \mathbb{L}^{q}(\FF)\otimes\mathbb{L}^{q}(\FF)$
for some $q>\alpha$. If $(Y^{1}-Y^{2})^{+}$ is of class
\mbox{\rm(D)}, $\xi^{1}\le\xi^{2}$, $dV^{1}\le dV^{2}$,
$L^{1}_{t}\le L^{2}_{t},\, U^{1}_{t}\le U^{2}_{t}$ for a.e. $t\in
[0,T]$ and \mbox{\rm(\ref{eq3.11})} or \mbox{\rm(\ref{eq3.12})} is
satisfied then $Y^{1}_{t}\le Y^{2}_{t}$, $t\in [0,T]$.
\end{wn}
\begin{dow}
The proof is analogous to that of Proposition \ref{propd4.1}, the
only difference being in the fact that in the latter part of the
proof we now refer to the proof of Corollary \ref{cor3.10} instead
of the proof of Proposition \ref{prop3.4}.
\end{dow}

\begin{tw}
\label{thd4.1} Assume \mbox{\rm(H2), (H3)}. Then there exists at
most one solution $(Y,Z,R)$ of RBSDE$(\xi,f+dV,L,U)$ such that
$Y\in \mathcal{D}^{p}$ for some $p>1$.
\end{tw}
\begin{dow}
Follows immediately from Proposition \ref{propd4.1}.
\end{dow}

\begin{tw}\label{stwdd4}
Assume \mbox{\rm(H2), (H3), (Z)}. Ten there exists at most one
solution $(Y,Z,R)$ of RBSDE$(\xi,f+dV,L,U)$ such that $Y$ is of
class \mbox{\rm(D)} and $Z\in\bigcup_{\beta>\alpha} M^{\beta}$.
\end{tw}
\begin{dow}
Follows immediately from Corollary \ref{cor.ad21}.
\end{dow}

\begin{tw}
\label{th.main3} Let $p>1$. Assume that \mbox{\rm(H1)--(H6)} are
satisfied.
\begin{enumerate}
\item [\rm(i)]There exists a solution
$(Y,Z,R)\in\DM^{p}\otimes M^{p}\otimes\mathcal{V}^{p}$ of
RBSDE$(\xi,f+dV,L,U)$ iff \mbox{\rm(H8)} is satisfied.
\item [\rm(ii)]Let $(Y^{n,n},Z^{n,n})\in \DM^{p}\otimes M^{p}$ be
a solution of the BSDE
\begin{align*}
Y^{n,n}_{t}&=\xi_{n}+\intt f(s,Y^{n,n}_{s},Z^{n,n}_{s})\,ds
+\intt dV_{s}\\
&\quad+n\intt(Y^{n,n}_{s}-L_{s})^{-}\,ds -n\intt
(Y^{n,n}_{s}-U_{s})^{+}\,ds-\intt Z^{n,n}_{s}\, dB_{s},\quad t\in
[0,T]
\end{align*}
with $\xi_{n}$ such that there exist
$\xi^{1}_{n},\xi^{2}_{n}\in\mathbb{L}^{p}(\FF_{T})$ with the
property that $\xi^{1}_{n}\le\xi_{n}\le \xi_{n}^{2}$,
$\xi^{1}_{n}\nearrow\xi,\,\, \xi^{2}_{n}\searrow\xi$. Then
\[
Y^{n,n}_{t}\rightarrow Y_{t},\quad t\in [0,T],\quad
Z^{n,n}\rightarrow Z, \quad\lambda\otimes P-\mbox{a.e.}
\]
and for every $q\in [1,2)$, $r\in [1,p)$,
\begin{equation}
\label{eq6.1}
E(\into|Z^{n,n}_{s}-Z_{s}|^{q}\,ds)^{r/q}\rightarrow0.
\end{equation}
\item[\rm(iii)]Let $(\bar{Y}^{n},\bar{Z}^{n},\bar{A}^{n})
\in\DM^{p}\otimes M^{p}\otimes \mathcal{V}^{+,p}$ be a solution of
$\overline{R}$BSDE$(\bar{\xi}_{n},\bar{f}_{n}+dV,U)$ with
\[
\bar{f}_{n}(t,y,z)=f(t,y,z)+n(y-L_{t})^{-}
\]
and $\bar{\xi}_{n}\in\mathbb{L}^{p}(\FF_{T})$ such that
$\bar{\xi}_{n}\nearrow\xi$. Then
\[
\bar{Y}_{t}^{n}\nearrow Y_{t},\quad t\in [0,T], \quad
\bar{Z}^{n}\rightarrow Z,\quad\lambda\otimes P\mbox{-a.e.}
\]
and
\[
d\bar{A}^{n}\le d\bar{A}^{n+1},\, n\in\mathbb{N}, \quad
\bar{A}^{n}_{t}\nearrow R^{-}_{t},\quad t\in [0,T],
\]
for every $q\in [1,2)$, $r\in [1,p)$,
\[
E(\into|\bar{Z}^{n}_{s}-Z_{s}|^{q}\,ds)^{r/q}\rightarrow0,
\]
and for every $\tau\in\mathcal{T}$,
\[
\bar{K}^{n}_{\tau}\rightarrow R^{+}_{\tau}\quad \mbox{weakly in }
\mathbb{L}^{p}(\FF_{T}),
\]
where $\bar{K}^{n}_{t}=n\intot(\bar{Y}^{n}_{s}-L_{s})^{-}\,ds$.
\item[\rm(iv)]Let $(\underline{Y}^{m},\underline{Z}^{m},\underline{K}^{m})
\in\DM^{p}\otimes M^{p}\otimes \mathcal{V}^{+,p}$ be a solution of
$\underline{R}$BSDE$(\underline{\xi}_{m}, \underline{f}_{m}+dV,L)$
with
\[
\underline{f}_{m}(t,y,z)=f(t,y,z)-m(y-U_{t})^{+}
\]
and $\underline{\xi}_{m}\in\mathbb{L}^{p}(\FF_{T})$ such that
$\underline{\xi}_{m}\searrow\xi$. Then
\[
\underline{Y}^{m}_{t}\searrow Y_{t},\quad t\in [0,T],
\quad\underline{Z}^{m}\rightarrow Z,\quad \lambda\otimes
P\mbox{-a.e.}
\]
and
\[
d\underline{K}^{m}\le d\underline{K}^{m+1},\, m\in\mathbb{N},\quad
\underline{K}^{m}_{t}\nearrow R^{+}_{t},\quad t\in [0,T],
\]
for every $q\in [1,2)$, $r\in [1,p)$,
\[
E(\into|\underline{Z}^{m}_{s}-Z_{s}|^{q}\,ds)^{r/q}\rightarrow0,
\]
and for every $\tau\in\mathcal{T}$,
\[
\underline{A}^{m}_{\tau}\rightarrow R^{-}_{\tau} \quad
\mbox{weakly in }\mathbb{L}^{p}(\FF_{T}),
\]
where $\underline{A}^{m}_{t}=\intot
m(\underline{Y}^{m}_{s}-U_{s})^{+}\,ds$.
\item[\rm(v)]If $L,U, V$ are continuous and $L_{T}\le\xi\le U_{T}$
then as $n,m\rightarrow+\infty$,
\[
E\sup_{0\le t\le T}|\underline{A}^{m}_{t}-R^{-}_{t}|^{p}
+E\sup_{0\le t\le T}|\bar{K}^{n}_{t}-R^{+}_{t}|^{p}\rightarrow 0,
\]
\[
E\sup_{0\le t\le T}|\bar{Y}^{n}_{t}-Y_{t}|^{p} +E\sup_{0\le t\le
T}|\underline{Y}^{m}_{t}-Y_{t}|^{p} +E\sup_{0\le t\le
T}|Y^{n,n}_{t}-Y_{t}|^{p}\rightarrow0
\]
and
\begin{align}
\label{eq6.03} &E(\into|Z^{n,n}_{s}-Z_{s}|^{2}\,ds)^{p/2}
+E(\into|\underline{Z}^{m}_{s}-Z_{s}|^{2}\,ds)^{p/2}\nonumber\\
&\qquad+E(\into|\bar{Z}^{n}_{s}-Z_{s}|^{2}\,ds)^{p/2}\rightarrow0.
\end{align}
\item[\rm(vi)]For every $\check{U},\hat{L}\in\DM$ such that
$L_{t}\le\hat{L}_{t}\le Y_{t}\le \check{U}_{t}\le U_{t}$ for  a.e. $t\in [0,T]$,
\[
\Delta R^{+}_{t}=(\hat{L}_{t-}-Y_{t}-\Delta V_{t})^{+}, \quad
\Delta R^{-}_{t}=(Y_{t}-\check{U}_{t-}+\Delta V_{t})^{+},\quad
t\in(0,T].
\]
\end{enumerate}
\end{tw}
\begin{dow}
Necessity. Let $(Y,Z,R)\in\DM^{p}\otimes M^{p}\otimes
\mathcal{V}^{p}$ be a solution of RBSDE$(\xi,f+dV,L,U)$. Then
$L_{t}\le Y_{t}\le U_{t}$ for a.e. $t\in [0,T]$ and by Proposition
\ref{prop.gen},
\[
E(\into|f(s,Y_{s},Z_{s})|\,ds)^{p}<\infty,
\]
which implies that $Y\in\mathcal{H}^{p}$. Thus, (H8) is satisfied
with $X=Y$.

Sufficiency. To shorten notation we give the proof under the
assumption that
$\bar{\xi}_{n}=\underline{\xi}_{m}=\xi^{1}_{n}=\xi^{2}_{n}=\xi$.
The proof in the general case is analogous. It only requires some
obvious changes. By Theorem \ref{th3.1}, for every
$n,m\in\mathbb{N}$ there exists a unique solution
$(Y^{n,m},Z^{n,m})\in\DM^{p}\otimes M^{p}$ of the BSDE
\begin{align}
\label{m3.0}
\nonumber Y^{n,m}_{t}&=\xi+\intt f(s,Y^{n,m}_{s},Z^{n,m}_{s})\,ds
+\intt dV_{s}+n\intt(Y^{n,m}_{s}-L_{s})^{-}\,ds\\
&\quad-m\intt(Y^{n,m}_{s}-U_{s})^{+}\,ds-\intt Z^{n,m}_{s}\,
dB_{s},\quad t\in [0,T].
\end{align}
Set
\begin{equation}
\label{eq6.33}
A^{n,m}_{t}=m\intot(Y^{n,m}_{s}-U_{s})^{+}\,ds,\quad K^{n,m}_{t}
=n\intot(Y^{n,m}-L_{s})^{-}\,ds,\quad t\in [0,T].
\end{equation}
Step 1. We first show that
\begin{equation}
\label{m3.1} \sup_{n,m\ge1}\left(E(\into dK^{n,m}_{s})^{p}
+E(\into dA^{n,m}_{s})^{p}\right)<\infty.
\end{equation}
Since $X\in\HH^{p}$ and the Brownian filtration has the
representation property, there exist $C\in\mathcal{V}^{p}, H\in
M^{p}$ such that
\[
X_{t}=X_{0}-\intot dC_{s}-\intot H_{s}\,dB_{s}.
\]
The above formula may be rewritten in the form
\[
X_{t}=X_{T}+\intt f(s,X_{s},H_{s})\,ds+\intt dV_{s}+\intt
dK^{'}_{s}-\intt dA^{'}_{s}-\intt H_{s}\,dB_{s},
\]
where
\[
K^{'}_{t}=\intot (f^{-}(s,X_{s},H_{s})\,ds+dC_{s}^{-}+dV^{-}_{s}),
\quad A^{'}_{t}=\intot (f^{+}(s,X_{s},H_{s})\,ds+dC_{s}^{+}+dV^{+}_{s}).
\]
By Theorem \ref{th3.1}, for every $m\in\mathbb{N}$ there exists a
solution $(\bar{X}^{m},\bar{H}^{m})\in\DM^{p}\otimes M^{p}$ of the
BSDE
\begin{align}
\label{gw17} \bar{X}^{m}_{t}&=X_{T}\vee \xi+\intt
f(s,\bar{X}^{m}_{s},\bar{H}^{m}_{s})\,ds +\intt dV_{s}+\intt
dK^{'}_{s}\nonumber\\
&\quad-\intt m(\bar{X}^{m}_{s}-U_{s})^{+}\,ds -\intt
\bar{H}^{m}_{s}\,dB_{s},\quad t\in [0,T].
\end{align}
Since $L_{t}\le X_{t}\le U_{t}$ for a.e. $t\in [0,T]$, we have
\begin{align*}
X_{t}&=X_{T}+\intt f(s,X_{s},Z_{s})\,ds +\intt dV_{s}
-m\intt(X_{s}-U_{s})^{+}\,ds\\
&\quad +\intt dK^{'}_{s}-\intt dA^{'}_{s}-\intt
H_{s}\,dB_{s},\quad t\in [0,T].
\end{align*}
Hence, by Proposition \ref{prop3.4}, $\bar{X}^{m}_{t}\ge X_{t}$
for a.e. $t\in [0,T]$, which implies that $\bar{X}^{m}_{t}\ge
L_{t}$ for a.e. $t\in [0,T]$. Consequently,
\begin{align*}
\bar{X}^{m}_{t}&=X_{T}\vee \xi +\intt
f(s,\bar{X}^{m}_{s},\bar{H}^{m}_{s})\,ds+\intt dV_{s} +\intt
dK^{'}_{s}\\&\quad+n\intt(\bar{X}^{m}_{s}-L_{s})^{-}\,ds
-m\intt(\bar{X}^{m}_{s}-U_{s})^{+}\,ds -\intt
\bar{H}^{m}_{s}\,dB_{s},\quad t\in [0,T].
\end{align*}
Applying once again Proposition \ref{prop3.4} we see that
$\bar{X}^{m}_{t}\ge Y^{n,m}_{t}$, $t\in [0,T]$, for every
$n,m\in\mathbb{N}$. Thus,
\begin{equation}
\label{m3.2}
dA^{n,m}=m(Y^{n,m}_{s}-U_{s})^{+}\,ds
\le m(\bar{X}^{m}_{s}-U_{s})^{+}\,ds,\quad n,m\in\mathbb{N}.
\end{equation}
Observe now that $(-\bar{X}^{m},-\bar{H}^{m})$ is a supersolution
of BSDE$(-\xi\vee X_{T},\tilde{f}-dV-dK^{'})$ with
\[
\tilde{f}(t,y,z)=-f(t,-y,-z).
\]
Since $K^{'}\in\mathcal{V}^{p}$, $-\bar{X}^{m}_{t}\le-X_{t}$,
$t\in [0,T]$, $X\in\DM^{p}$ and $E(\into
\tilde{f}^{-}(s,-X_{s},0)\,ds)^{p}<\infty$, it follows from
Proposition \ref{prop4.3} that there exists $C>0$ not depending on
$n$ such that $E(\into dA^{n,m}_{s})^{p}\le C$. The same
conclusion can be drawn for $K^{n,m}$. To see this it suffices to
consider a solution $(\underline{X}^{n},\underline{H}^{n})\in
\DM^{p}\otimes M^{p}$ of the BSDE
\begin{align}
\label{eq6.7} \underline{X}^{n}_{t}&=X_{T}\wedge \xi
+\intt f(s,\underline{X}^{n}_{s},\underline{H}^{n}_{s})\,ds
+\intt dV_{s} \nonumber\\
&\quad +\intt n(\underline{X}^{n}_{s}-L_{s})^{-}\,ds -\intt
dA^{'}_{s}-\intt \underline{H}^{n}_{s}\,dB_{s},\quad t\in[0,T]
\end{align}
and then repeat (with some obvious changes) arguments following
(\ref{gw17}).
\medskip\\
Step 2. We will show that there exists a triple $(Y,Z,R)\in
\DM^{p}\otimes M^{p}\otimes \mathcal{V}^{p}$ which is in some
sense a limit of triple $(\bar{Y}^{n},\bar{Z}^{n},
\bar{K}^{n}-\bar{A}^{n})$ and
\begin{equation}
\label{0}
Y_{t}=\xi+\intt f(s,Y_{s},Z_{s})\,ds+\intt dV_{s}+\intt dR_{s}
-\intt Z_{s}\,dB_{s},\quad t\in [0,T].
\end{equation}
Let us first observe that
\begin{equation}
\label{m3.3}
d\bar{A}^{n}\le d\bar{A}^{n+1},\quad n\in\mathbb{N}.
\end{equation}
Indeed, by Theorem \ref{th.main1}, for every $t\in [0,T]$,
$\bar{A}^{n}_{t},\,\bar{A}^{n+1}_{t}$ are weak limits in
$\mathbb{L}^{p}(\FF_{T})$ of  $\{A^{n,m}_{t}\}$ and
$\{A^{n+1,m}_{t}\}$, respectively. This implies (\ref{m3.3})
because by Proposition \ref{prop3.4}, $Y^{n,m}_{t}\le
Y^{n+1,m}_{t}$, $t\in [0,T]$, for every $n,m\in\mathbb{N}$, and
consequently
\[
dA^{n,m}=m(Y^{n,m}_{t}-U_{t})^{+}\,dt \le
m(Y^{n+1,m}_{t}-U_{t})^{+}\,dt=dA^{n+1,m}.
\]
Set
\begin{equation}
\label{temp1} A_{t}=\sup_{n\ge1}\bar{A}^{n}_{t},\quad t\in[0,T].
\end{equation}
By (\ref{m3.1}), (\ref{m3.3}) and \cite[Lemma 2.2]{Peng},
$A\in\mathcal{V}^{+,p}$. Now observe that by Proposition
\ref{prop3.4} and Theorem \ref{th.main1},
\begin{equation}
\label{star234}
\bar{y}_{t}\le Y^{n,m}_{t}\le \underline{y}_{t},\quad
t\in[0,T],\quad n,m\in\mathbb{N},
\end{equation}
where $(\bar{y},\bar{z},\bar{k})\in\DM^{p}\otimes
M^{p}\otimes\mathcal{V}^{+,p}$ (resp.
$(\underline{y},\underline{z},\underline{k})\in\DM^{p}\otimes
M^{p}\otimes\mathcal{V}^{+,p}$) is a solution of
$\overline{\mbox{R}}$BSDE$(\xi,f+dV,U)$ (resp.
$\underline{\mbox{R}}$BSDE$(\xi,f+dV,L)$). By Theorem
\ref{th.main1},
\begin{equation}
\label{m3.5} Y^{n,m}_{t}\searrow\bar{Y}^{n}_{t},\quad t\in
[0,T],\quad n\in\mathbb{N}.
\end{equation}
Hence
\begin{equation}
\label{m3.6}
\bar{y}_{t}\le\bar{Y}^{n}_{t}\le \underline{y}_{t},
\quad t\in [0,T],\quad n\in\mathbb{N}.
\end{equation}
By Proposition \ref{prop4.1}, for every $n\in\mathbb{N}$,
$\bar{Y}^{n}_{t}\le\bar{Y}^{n+1}_{t}$ for $t\in [0,T]$. Therefore
setting $Y_{t}=\sup_{n\ge1}\bar{Y}^{n}_{t}$, $t\in[0,T]$, we see
that
\begin{equation}
\label{temp2}
Y_{t}=\lim_{n\rightarrow +\infty} \bar{Y}^{n}_{t},\quad t\in [0,T].
\end{equation}
From Step 1 and what has already been proved in Step 2 we conclude
that assumptions (a), (c)-(f) of Theorem \ref{monbsde} are
satisfied. Let $\sigma,\tau\in\mathcal{T}$ be such that $\sigma\le
\tau$. Then by (\ref{m3.3}), (\ref{temp1}), (\ref{m3.6}) and
(\ref{temp2}),
\begin{align*}
\int_{\sigma}^{\tau}(Y_{s}-\bar{Y}^{n}_{s})\,
d(\bar{K}^{n}_{s}-\bar{A}^{n}_{s}+V_{s})&
\ge -\int_{\sigma}^{\tau}(Y_{s}-\bar{Y}^{n}_{s})\, d(\bar{A}^{n}_{s}-V_{s})\\
&\ge-\int_{\sigma}^{\tau}(Y_{s}-\bar{Y}^{n}_{s})\,d(A_{s}-V_{s})\rightarrow
0.
\end{align*}
This shows that assumption (b) of Theorem \ref{monbsde} is
satisfied as well. By Theorem \ref{monbsde} and (\ref{m3.6}),
(\ref{m3.1}) there exists a quadruple $(Y,Z,K,A)\in \DM^{p}\otimes
M\otimes \mathcal{V}^{p,+}\otimes \mathcal{V}^{p,+}$ such that
\begin{equation}
\label{m3.15} \bar{Z}^{n}\rightarrow Z,\quad \lambda\otimes
P\mbox{-a.e.}
\end{equation}
and (\ref{0}) is satisfied with $R=K-A$.
Using (\ref{m3.1}), Proposition \ref{prop3.1} and Proposition
\ref{prop.gen} we conclude from (\ref{m3.0}) that there exists $C$
not depending on $n,m$ such that
\begin{equation}
\label{eq.ret}
E(\into|f(s,Y^{n,m}_{s},Z^{n,m}_{s})|\,ds)^{p}+E(\into |Z^{n,m}_{s}|^{2})^{p/2}\le C.
\end{equation}
From this we deduce that
\begin{equation}
\label{00} \bar{K}^{n}_{\tau}\rightarrow K_{\tau}
\quad\mbox{weakly in }\mathbb{L}^{p}(\FF_{T})
\end{equation}
for every $\tau\in\mathcal{T}$,
$f(\cdot,Y_{\cdot},Z_{\cdot})\in\mathbb{L}^{1,p}(\FF),\, Z\in M^{p}$ and  (\ref{eq6.1})
holds true for every $q\in [1,2), r\in [1,p)$.

Step 3. We will show that $K,A$ satisfy the minimality condition.
From (\ref{m3.0}), (\ref{m3.1}) and (\ref{eq.ret}) it may be
concluded that there exists $C$ not depending on $n,m$ such that
\[
E(\into (Y^{n,m}_{s}-L_{s})^{-}\,ds)^{p} +E(\into
(Y^{n,m}_{s}-U_{s})^{+}\,ds)^{p}\le C(n^{-p}+m^{-p}).
\]
By the above inequality, (\ref{m3.5}) and (\ref{temp2}),
\begin{equation}
\label{m3.17} L_{t}\le Y_{t}\le U_{t}\quad\mbox{for a.e. }t\in
[0,T].
\end{equation}
Let us fix $\check{U}\in\DM$ such that $Y_{t}\le \check{U}_{t}\le
U_{t}$ for  a.e. $t\in [0,T]$. By (\ref{m3.3}) and  (\ref{temp1}),
$d\bar{A}^{n}\rightarrow dA$ in the variation norm. Hence
\[
\into (\check{U}_{s-}-Y_{s-})\,d\bar{A}^{n}_{s}\rightarrow
\into (\check{U}_{s-}-Y_{s-})\,dA_{s}.
\]
On the other hand,
\[
0\le \into (\check{U}_{s-}-Y_{s-})\,d\bar{A}^{n}_{s} \le \into
(\check{U}_{s-}-\bar{Y}^{n}_{s-})\,d\bar{A}^{n}_{s}=0,
\]
because the triple $(\bar{Y}^{n},\bar{Z}^{n},\bar{A}^{n})$ is a
solution of $\overline{\mbox{R}}$BSDE$(\xi,\bar{f}_{n}+dV,U)$.
Consequently,
\begin{equation}
\label{eq6.2}\into (\check{U}_{s-}-Y_{s-})\,dA_{s}=0,
\end{equation}
i.e. $A$ satisfies the minimality condition. Suppose now that
$\Delta A_{t}>0$. Since $d\bar{A}^{n}\rightarrow dA$ in the
variation norm, there exists $n_{0}\in\mathbb{N}$ such that
$\Delta\bar{A}^{n}_{t}>0$ for every $n\ge n_{0}$. Since
$\bar{A}^{n}$ satisfies the minimality condition,
$\Delta\bar{A}^{n}_{t}=\bar{Y}^{n}_{t}-\bar{Y}^{n}_{t-}+\Delta
V_{t}=\bar{Y}^{n}_{t}-\check{U}_{t-}+\Delta V_{t}$, so using once
again the fact that $d\bar{A}^{n}\rightarrow dA$ we get
\begin{equation}
\label{kto}
\Delta A_{t}=Y_{t}-\check{U}_{t-}+\Delta V_{t}.
\end{equation}
Since $A$ satisfies the minimality condition,
$Y_{t-}=\check{U}_{t-}$\,. Hence
\[
Y_{t}-\check{U}_{t-}=\Delta Y_{t}=-\Delta K_{t}-\Delta
V_{t}+\Delta A_{t}=-\Delta K_{t}+Y_{t}-\check{U}_{t-}\,,
\]
which forces $\Delta K_{t}=0$. Thus, we have shown that for every
$t\in(0,T]$, if $\Delta K_{t}>0$ then $\Delta A_{t}=0$. Let us fix
$\hat{L}\in\DM$ such that $L_{t}\le \hat{L}_{t}\le Y_{t}$ for a.e.
$t\in [0,T]$. By Lemma \ref{lm4.6},
\begin{equation}
\label{m3.19}
Y_{t-}=\hat{L}_{t-}\vee(Y_{t}+\Delta V_{t}-\Delta A_{t}),\quad t\in (0,T].
\end{equation}
Suppose that $\Delta K_{t}>0$. Then
$Y_{t-}=\hat{L}_{t-}\vee(Y_{t}+\Delta V_{t})$ since  $\Delta
A_t=0$. On the other hand, $\Delta K_{t}=-\Delta Y_{t}-\Delta
V_{t}$, which implies that $Y_{t-}>Y_{t}+\Delta V_{t}$. The last
inequality when combined with (\ref{m3.19}) shows that
$Y_{t-}=\hat{L}_{t-}$\,. Thus, for every $t\in(0,T]$, if $\Delta
K_{t}>0$ then $Y_{t-}=\hat{L}_{t-}$\,. Consequently,
\begin{equation}
\label{m3.20} \sum_{0<t\le T}(Y_{t-}-\hat{L}_{t-})\Delta K_{t}=0.
\end{equation}
Let $\{\tau_{k}\}$ be an increasing sequence of successive jumps
of $K,A,V$ with the convention that $\tau_{0}\equiv 0$. Then by
Dini's theorem, $\bar{Y}^{n}\rightarrow Y$  uniformly on compact subsets of
$(\tau_{k},\tau_{k+1})$ for each $k\in\mathbb{N}$. Since
$\{\bar{A}^{n}\}$ converges pointwise, in much the same way as in
Step 3 of the proof of Theorem \ref{th.main1} one can show that
for every $k\in\mathbb{N}$,
\[
\int_{(\tau_{k},\tau_{k+1})}(Y_{t-}-\hat{L}_{t-})\,dK_{t}=0.
\]
Combining this with (\ref{m3.20}) we get
\begin{equation}
\label{eq6.3} \into(Y_{t-}-\hat{L}_{t-})\,dK_{t}=0,
\end{equation}
which proves that $K$ satisfies the minimality condition.

Step 4.  We will show (i). Write $R_{t}=K_{t}-A_{t}$. Since
$K,A\in\mathcal{V}^{+,p}$, it follows from the minimality property
of the Jordan decomposition of signed measures that $dR^{+}\le
dK$, $dR^{-}\le dA$. From this and (\ref{eq6.2}), (\ref{eq6.3}) we
get (\ref{eq6.01}), which when combined with (\ref{0}),
(\ref{m3.17}) and integrability properties of $Y,Z,K,A$ proved in
Step 2 shows that  the triple $(Y,Z,R)\in\DM^{p}\otimes
M^{p}\otimes \mathcal{V}^{p}$ is a solution of
RBSDE$(\xi,f+dV,L,U)$.

Step 5. We will show (iii) and (iv). By what has already been
proved in Steps 1-3, assertion (iii) will be proved once we prove
that $A_{t}=R^{-}_{t}$, $K_{t}=R^{+}_{t}$, $t\in [0,T]$. To prove
these equalities let us first note that by the same method as in
the proof of (\ref{m3.2}), but with the process $X$ replaced by
$Y$, the processes $K', A'$ replaced by $R^{+}, R^{-}$ and
$(\bar{X}^{m},\bar{H}^{m})$ replaced by the solution
$(\tilde{Y}^{m},\tilde{Z}^{m})$ of the  BSDE
\begin{align*}
\tilde{Y}^{m}_{t}&=\xi+\intt f(s,\tilde{Y}^{m}_{s},\tilde{Z}^{m}_{s})\,ds
+\intt dV_{s}+\intt dR^{+}_{s}\\
&\quad-m\intt(\tilde{Y}^{m}_{s}-U_{s})^{+}\,ds
-\intt\tilde{Z}^{m}_{s}\,dB_{s},\quad t\in [0,T]
\end{align*}
one can show that
\begin{equation}
\label{m3.22}
dA^{n,m}\le d\tilde{A}^{m}\quad n,m\in\mathbb{N},
\end{equation}
where $\tilde{A}^m_{t}=m\intot(\tilde{Y}^{m}_{s}-U_{s})^{+}\,ds,\,
t\in [0,T]$.  By Theorems \ref{th4.1} and \ref{th.main1}, the
sequence $\{(\tilde{Y}^{m},\tilde{Z}^{m},\tilde{A}^{m})\}$
converges in the sense of Theorem \ref{th.main1} to $(Y,Z,R^{-})$.
Since the arithmetic mean preserves inequalities,
%(this relation will be the only crucial thing in our reasoning?)
without loss of generality we may assume that
\[
A^{n,m}_{t}\rightarrow \bar{A}^{n}_{t},\quad
\tilde{A}^{m}_{t}\rightarrow R^{-}_{t},\quad t\in [0,T]
\]
(see the reasoning in Step 3 of the proof of Theorem
\ref{th.main1}). Therefore from (\ref{m3.3}), (\ref{temp1}), and
(\ref{m3.22}) it follows that $dA\le dR^{-}$. By the minimality
property of the Jordan decomposition, $R^{-}_{t}=A_{t}$, $t\in
[0,T]$, and consequently $R^{+}_{t}=K_{t}$, $t\in [0,T]$, which
completes the proof of (iii). Applying similar arguments to
$\{(\underline{Y}^{m},\underline{Z}^{m},
\underline{K}^{m},\underline{A}^{m})\}$ proves (iv).

Step 6. We will show (ii). By Proposition \ref{prop3.4},
\begin{equation}
\label{mest} \bar{Y}^{n}_{t}\le Y^{n,n}_{t}\le
\underline{Y}^{n}_{t},\quad t\in [0,T].
\end{equation}
By (\ref{mest}) and (iii), (iv),
\begin{equation}
\label{m3.23}
Y^{n,n}_{t}\rightarrow Y_{t},\quad t\in [0,T].
\end{equation}
Therefore the assumptions (a), (c)-(f) of Lemma \ref{dombsde} are
satisfied. By (\ref{star234}), (\ref{m3.17}), (\ref{m3.23}) for
every $\sigma,\tau\in\mathcal{T}$ such that $\sigma\le \tau$,
\begin{align*}
&\int_{\sigma}^{\tau}(Y_{s}-Y^{n,n}_{s})\,d(K^{n,n}-A^{n,n})+V)_{s}
=n\int_{\sigma}^{\tau}(Y_{s}-Y^{n,n}_{s})(Y^{n,n}_{s}-L_{s})^{-}\,ds\\
&\qquad
-n\int_{\sigma}^{\tau}(Y_{s}-Y^{n,n}_{s})(Y^{n,n}_{s}-U_{s})^{+}\,ds
+\int_{\sigma}^{\tau}(Y_{s}-Y^{n,n}_{s})\,dV_{s}\\
&\quad\ge
n\int_{\sigma}^{\tau}(L_{s}-Y^{n,n}_{s})(Y^{n,n}_{s}-L_{s})^{-}\,ds
-n\int_{\sigma}^{\tau}(U_{s}-Y^{n,n}_{s})(Y^{n,n}_{s}-U_{s})^{+}\,ds\\
&\qquad+\int_{\sigma}^{\tau}(Y_{s}-Y^{n,n}_{s})\,dV_{s}\ge
\int_{\sigma}^{\tau}(Y_{s}-Y^{n,n}_{s})\,dV_{s}\rightarrow 0.
\end{align*}
This shows that assumption (b) of Lemma \ref{dombsde} is satisfied
as well. Consequently, $Z^{n,n}\rightarrow Z$, $\lambda\otimes
P$-a.e., which together with (\ref{eq.ret}) implies (\ref{eq6.1}).

Step 7. We will show (v) and (vi). Observe that (vi) follows
immediately from (\ref{eq6.2}), (\ref{kto}), (\ref{eq6.3}) and the
equalities  $K=R^{+}, A=R^{-}$ proved in Step 5. From (vi) and the
fact that $V$ is continuous it follows that the processes $K,A,Y$
are continuous. By Theorem \ref{th.main1} the  processes
$\bar{Y}^{n}$, $\bar{K}^{n}$, $\underline{Y}^{m}$,
$\underline{A}^{m}$ are continuous as well. Using Dini's theorem,
integrability of $\bar{Y}^{1}$, $\bar{K}^{1}$,
$\underline{Y}^{1}$, $\underline{A}^{1}$, $A$, $K$, $Y$,
(\ref{mest}) and the Lebesgue dominated convergence theorem shows
the desired convergence of the sequences $\{\bar{Y}^{n}\}$,
$\{\bar{K}^{n}\}$, $\{\underline{Y}^{m}\}$,
$\{\underline{A}^{m}\}$. To prove (\ref{eq6.03}) let us first
observe that by It\^o's formula, (H2) and (H3),
\begin{align*}
\int_{0}^{T}|Z^{n,n}_{t}-Z_{t}|^{2}\,dt &\le
2\lambda\int_{0}^{T}|Y^{n,n}_{t}-Y_{t}||Z^{n,n}_{t}-Z_{t}|\,dt
+2\int_{0}^{T}|Y^{n,n}_{t}-Y_{t}|\,dR^{+}_{t}\\
&\quad+2\int_{0}^{T}|Y^{n,n}_{t}-Y_{t}|\,dR^{-}_{t}+\sup_{0\le
t\le T}|\intt (Z^{n}_{s}-Z^{m}_{s})
(Y^{n}_{s}-Y^{m}_{s})\,dB_{s}|.
\end{align*}
Hence
\begin{align*}
&E(\int_{0}^{T}|Z^{n,n}_{t}-Z_{t}|^{2}\,dt)^{p/2}\le
C(p,\lambda)\bigg(E\sup_{0\le t\le T}|Y^{n,n}_{t}-Y_{t}|^{p}\\
&\quad+(E\sup_{0\le t\le T}|Y^{n,n}_{t}-Y_{t}|^{p})^{1/2}
(E|R^{+}_{T}|^{p})^{1/2} +(E\sup_{0\le t\le
T}|Y^{n,n}_{t}-Y_{t}|^{p})^{1/2} (E|R^{-}_{T}|^{p})^{1/2}\bigg),
\end{align*}
which together with  uniform convergence of $\{Y^{n,n}\}$ implies
(\ref{eq6.03}). Since the proof of the other convergences in (v)
is similar, we omit it.
\end{dow}

\begin{tw}
\label{th.main4} Assume that \mbox{\rm(H1)--(H6), (Z)} are
satisfied with $p=1$.
\begin{enumerate}
\item [\rm(i)]There exists a solution
$(Y,Z,R)\in\DM^{q}\otimes M^{q}\otimes \mathcal{V}^{1}$, $q\in
(0,1)$, such that $Y$ is of class \mbox{\rm(D)} iff
\mbox{\rm(H8*)} is satisfied.
\item[\rm(ii)]Let $(Y^{n,n},Z^{n,n})\in\DM^{q}\otimes M^{q},\,q\in
(0,1)$, such that $Y^{n,n}$ is of class \mbox{\rm(D)} be a
solution of the BSDE
\begin{align*}
Y^{n,n}_{t}&=\xi+\intt f(s,Y^{n,n}_{s},Z^{n,n}_{s})\,ds
+\intt dV_{s}+n\intt(Y^{n,n}_{s}-L_{s})^{-}\,ds\\
&\quad-n\intt(Y^{n,n}_{s}-U_{s})^{+}\,ds-\intt Z^{n,n}_{s}\,
dB_{s}, \quad t\in [0,T].
\end{align*}
Then
\[
Y^{n,n}_{t}\rightarrow Y_{t},\quad t\in [0,T],\quad
Z^{n,n}\rightarrow Z, \quad \lambda\otimes P\mbox{-a.e.}
\]
and there exists a stationary sequence $\{\tau_{k}\}\subset
\mathcal{T}$ such that for every $q\in[1,2)$, $r>1$,
\[
E(\int_{0}^{\tau_{k}}|Z^{n,n}_{s}-Z_{s}|^{q}\,ds)^{r}\rightarrow
0.
\]
\item [\rm(iii)] Let
$(\bar{Y}^{n},\bar{Z}^{n},\bar{A}^{n})\in \DM^{q}\otimes M^{q}
\otimes \mathcal{V}^{+,1}$, $q\in (0,1)$, such that $\bar{Y}^{n}$
is of class \mbox{\rm(D)} be a solution of
$\overline{R}$BSDE$(\bar{\xi},\bar{f}_{n}+dV,U)$ with
\[
\bar{f}_{n}(t,y,z)=f(t,y,z)+n(y-L_{t})^{-}.
\]
Then
\[
\bar{Y}_{t}^{n}\nearrow Y_{t},\quad t\in [0,T],\quad \bar{Z}^{n}
\rightarrow Z,\quad \lambda\otimes P\mbox{-a.e.},
\]
\[
d\bar{A}^{n}\le d\bar{A}^{n+1},\, n\in\mathbb{N}, \quad
\bar{A}^{n}_{t}\nearrow R^{-}_{t},\quad t\in [0,T]
\]
and there exists a stationary sequence
$\{\tau_{k}\}\subset\mathcal{T}$ such that for every
$\tau\in\mathcal{T}$ and $q\in[1,2)$, $r>1$,
\[
E(\int_{0}^{\tau_{k}} |\bar{Z}^{n}_{s}-Z_{s}|^{q}\,ds)^{r}
\rightarrow 0,\quad \bar{K}^{n}_{\tau_{k}\wedge\tau}\rightarrow
R^{+}_{\tau_{k}\wedge \tau}\quad \mbox{weakly in }
\mathbb{L}^{r}(\FF_{T}),
\]
where $\bar{K}^{n}_{t}=n\intot(\bar{Y}^{n}_{s}-L_{s})^{-}\,ds.$
\item [\rm(iv)] Let
$(\underline{Y}^{m},\underline{Z}^{m},\underline{K}^{m})\in
\DM^{q}\otimes M^{q}\otimes \mathcal{V}^{+,1}$, $q\in(0,1)$, such
that $\underline{Y}^{m}$ is of class \mbox{\rm(D)} be a solution
of $\underline{R}$BSDE$(\underline{\xi},\underline{f}_{m}+dV,L)$
with
\[
\underline{f}_{m}(t,y,z)=f(t,y,z)-m(y-U_{t})^{+}.
\]
Then
\[
\underline{Y}^{m}_{t}\searrow Y_{t},\quad t\in [0,T], \quad
\underline{Z}^{m}\rightarrow Z,\quad \lambda\otimes P\mbox{-a.e.},
\]
\[
d\underline{K}^{m}\le d\underline{K}^{m+1},\, m\in\mathbb{N},
\quad \underline{K}^{m}_{t}\nearrow R^{+}_{t},\quad t\in [0,T]
\]
and there exists a stationary sequence $\{\tau_{k}\}\subset
\mathcal{T}$ such that for every $\tau\in\mathcal{T}$ and
$q\in[1,2)$, $r>1$,
\[
E(\int_{0}^{\tau_{k}} |\underline{Z}^{m}_{s}-Z_{s}|^{q}\,ds)^{r}
\rightarrow 0,\quad \underline{A}^{m}_{\tau_{k}\wedge\tau}
\rightarrow R^{-}_{\tau_{k}\wedge\tau}\quad \mbox{weakly in }
\mathbb{L}^{r}(\FF_{T}),
\]
where
$\underline{A}^{m}_{t}=m\intot(\underline{Y}^{m}_{s}-U_{s})^{+}\,ds$.
\item[\rm(v)]If $L,U,V$ are continuous and $L_{T}\le\xi\le U_{T}$
then there exists a stationary sequence $\{\tau_{k}\}\subset
\mathcal{T}$ such that for every $p>1$,
\[
E\sup_{0\le t\le \tau_{k}}|\underline{A}^{m}_{t}-R^{-}_{t}|^{p}
+E\sup_{0\le t\le
\tau_{k}}|\bar{K}^{n}_{t}-R^{+}_{t}|^{p}\rightarrow0,
\]
\[
E\sup_{0\le t\le \tau_{k}}|\bar{Y}^{n}_{t}-Y_{t}|^{p}
+E\sup_{0\le t\le \tau_{k}}|\underline{Y}^{m}_{t}-Y_{t}|^{p}
+E\sup_{0\le t\le \tau_{k}}|Y^{n,n}_{t}-Y_{t}|^{p}\rightarrow 0
\]
and
\[
E(\int_{0}^{\tau_{k}}|Z^{n,n}_{s}-Z_{s}|^{2}\,ds)^{p/2}
+E(\int_{0}^{\tau_{k}}|\underline{Z}^{m}_{s}-Z_{s}|^{2}\,ds)^{p/2}
+E(\int_{0}^{\tau_{k}}|\bar{Z}^{n}_{s}-Z_{s}|^{2}\,ds)^{p/2}
\rightarrow0.
\]
Moreover,
\begin{equation}
\label{eq6.29} E\sup_{0\le t\le
T}|\underline{A}^{m}_{t}-R^{-}_{t}| +E\sup_{0\le t\le
T}|\bar{K}^{n}_{t}-R^{+}_{t}|\rightarrow 0,
\end{equation}
and for every $q\in(0,1)$,
\begin{equation}
\label{eq6.30} E\sup_{0\le t\le T}|\bar{Y}^{n}_{t}-Y_{t}|^{q}
+E\sup_{0\le t\le T}|\underline{Y}^{m}_{t}-Y_{t}|^{q} +E\sup_{0\le
t\le T}|Y^{n,n}_{t}-Y_{t}|^{q}\rightarrow0.
\end{equation}
\item [\rm(vi)]For every $\check{U}, \hat{L}\in\DM$ such that
$\hat{L}_{t}\le Y_{t}\le \check{U}_{t}$ for $t\in [0,T]$,
\[
\Delta R^{+}_{t}=(\hat{L}_{t-}-Y_{t}-\Delta V_{t})^{+}, \quad
\Delta R^{-}_{t}=(Y_{t}-\check{U}_{t-}+\Delta V_{t})^{+}, \quad
t\in(0,T].
\]
\end{enumerate}
\end{tw}
\begin{dow}
First let us note that by Theorem \ref{th.main2} there exist
solutions $(Y^{n,n},Z^{n,n})$,
$(\bar{Y}^{n},\bar{Z}^{n},\bar{A}^{n})$,
$(\underline{Y}^{m},\underline{Z}^{m},\underline{K}^{m})$ of
equations of assertions (ii)--(iv) of the theorem having the
required integrability properties. By Theorem \ref{th3.2} there
exists a solution $(\bar{X}^{m},\bar{H}^{m})$ of (\ref{gw17}) and
a solution $(\underline{X}^{n},\underline{H}^{n})$ of
(\ref{eq6.7}) such that $(\bar{X}^{m},\bar{H}^{m})$,
$(\underline{X}^{n},\underline{H}^{n})\in\DM^{q}\otimes M^{q}$ for
$q\in (0,1)$ and  $\bar{X}^{m},\,\underline{X}^{n}$ are of class
(D). Using this in much the same way as in Step 1 of the proof of
Theorem \ref{th.main3} we show that
\begin{equation}
\label{AHA} \sup_{n,m\ge1} \left(E\into dA^{n,m}_{s} +E\into
dK^{n,m}_{s}\right)<\infty,
\end{equation}
where $A^{n,m}, K^{n,m}$ are defined by (\ref{eq6.33}). By
Corollary \ref{cor3.10},
\begin{equation}
\label{AHA1} \bar{y}_{t}\le Y^{n,m}_{t}\le\underline{y}_{t},\quad
t\in[0,T],
\end{equation}
where $Y^{n,m}$ is defined by (\ref{m3.0}) and
$(\bar{y},\bar{z},\bar{k})$,
$(\underline{y},\underline{z},\underline{k})\in\DM^{q}\otimes
M^{q}\otimes\mathcal{V}^{+,1}$ are solutions of
$\overline{\mbox{R}}$BSDE$(\xi,f+dV,U)$ and
$\underline{\mbox{R}}$BSDE$(\xi,f+dV,L)$), respectively, such that
$\bar{y},\underline{y}$ are of class (D). Since
$X\in\mathcal{V}^{1}+\MM^{q}$, $q\in (0,1)$, there exist
$C\in\mathcal{V}^{1}$ and $H\in M^{q}$ such that
\[
X_{t}=X_{0}-\intot dC_{s}-\intot H_{s}\,dB_{s},\quad t\in [0,T].
\]
Let
\[
\delta^{1}_{k}=\inf\{t\in [0,T], \intot |f(s,0,0)|\,ds
+\intot |f(s,X_{s},0)|\,ds>k\}\wedge T.
\]
By Lemma \ref{lm4.55} there exist a stationary sequence
$\{\delta^{2}_{k}\}\subset\mathcal{T}$ and constants $c_{k}$ such
that $|V|_{\delta^2_{k}}+\bar{y}^{*}_{\delta^2_{k}}
+\underline{y}^{*}_{\delta^2_{k}}+|C|_{\delta^2_{k}}\le c_{k}$ for
$k\in\BN$. Let us put
$\tau_{k}=\delta^{1}_{k}\wedge\delta^{2}_{k}$ and observe that by
the definition of $\tau_{k}$, the data
$(\bar{Y}^{n}_{\tau_{k}},f,V,L,U)$ satisfy the assumptions of
Theorem \ref{th.main3} on the interval $[0,\tau_{k}]$ for every
$p>1$. Therefore on each interval $[0,\tau_k]$ the sequence
$\{(\bar{Y}^{n},\bar{Z}^{n},\bar{K}^{n},\bar{A}^{n})\}$ converges
in the sense of Theorem \ref{th.main3} to the unique solution
$(Y^{k},Z^{k}, R^{k,+},R^{k,-})$ of RBSDE$(\xi^{k},f+dV,L,U)$ with
$\xi^{k}=\lim_{n\rightarrow+\infty}\bar{Y}^{n}_{\tau_{k}}$. By
stationarity of the sequence $\{\tau_{k}\}$, there exists a triple
$(Y,Z,R)\in\DM\otimes M\otimes\mathcal{V}$ such that $t\mapsto
f(t,Y_{t},Z_{t})\in\mathbb{L}^{1}(0,T)$,
\[
Y_{t}=\xi+\intt f(s,Y_{s},Z_{s})\,ds
+\intt dV_{s}+\intt dR_{s}-\intt Z_{s}\,dB_{s},\quad t\in [0,T],
\]
\[
L_{t}\le Y_{t}\le U_{t}\,\quad\mbox{for a.e. }t\in [0,T]
\]
and
\[
\into (Y_{t-}-\hat{L}_{t-})\,dR^{+}_{t} =\into
(\check{U}_{t-}-Y_{t-})\,dR^{-}_{t}=0
\]
for every $\hat{L}, \check{U}\in\DM$ such that
$L_{t}\le\hat{L}_{t}\le Y_{t}\le\check{U}_{t}\le U_{t}$  for a.e.
$t\in [0,T]$. Actually, by (\ref{AHA}) and (\ref{AHA1}),
$R,Y\in\DM^{q}$ for $q\in(0,1)$, which implies that $Z\in M^{q}$
by Lemma \ref{lm3.1}. All the desired in (iv) and (v) convergences
of the sequence
$\{(\bar{Y}^{n},\bar{Z}^{n},\bar{K}^{n},\bar{A}^{n})\}$, except
for the case $q=1$ in (iv) and $q\in(0,1)$ in (v), and all the
desired properties of the limits follow now from Theorem
\ref{th.main3} and stationarity of $\{\tau_{k}\}$. If $L,U,V$ are
continuous and $L_{T}\le \xi\le U_{T}$ then by Theorem
\ref{th.main3} and stationarity of $\{\tau_{k}\}$, the processes
$Y,R$ are continuous. Therefore using the fact that
$Y\in\mathcal{S}^{q}$, $q\in (0,1)$,
$R^{+}\in\mathcal{V}^{+,1}_{c}$, the monotone character of the
convergence of $\{\bar{Y}^{n}\}$ and $\{\bar{K}^{n}\}$, Dini's
theorem and the Lebesgue dominated convergence theorem one can
show the required convergence in assertion (iv) for $q=1$ and in
assertion (v) for  $q\in(0,1)$. The same reasoning may be applied
to the quadruple
$(\underline{Y}^{m},\underline{Z}^{m},\underline{K}^{m},
\underline{A}^{m})$. Finally, since
\[
\bar{Y}^{n}_{t}\le Y^{n,n}_{t}\le\underline{Y}^{n}_{t},\quad t\in
[0,T],
\]
the data $(Y^{n,n}_{\tau_{k}},f,V,L,U)$  satisfy the assumptions
of Theorem \ref{th.main3} on $[0,\tau_{k}]$. Therefore the
required in (v) convergences of $\{Y^{n,n}\},\{Z^{n,n}\}$ follow
from the above inequality and Theorem \ref{th.main3}.
\end{dow}

\nsubsection{Nonintegrable solutions of reflected BSDEs}
\label{sec7}

In this section we consider reflected BSDEs with monotone
generator and barriers satisfying only the standard Mokobodzki
condition. In the case of BSDEs with one reflecting barrier this
means that we assume (H1)--(H6) ((H1)--(H6) and (Z) in case $p=1$)
and that $L\in\mathbb{L}^{\infty,p}(\FF)$ ($L$ is of class (D) in
case $p=1$). In the case of two barriers this means that we assume
(H1)--(H6) ((H1)--(H6) and (Z) in case $p=1$) and the standard
Mokobodzki condition (M) ((M*) in case $p=1$) formulated later on.
Theorems \ref{th.main1}, \ref{th.main2}, \ref{th.main3} and
\ref{th.main4} say that in general we can not expect existence of
$\mathbb{L}^{p}$ solutions. Nevertheless we show that under the
standard Mokobodzki condition there exist solutions having weaker
integrability properties. In fact, it may happen that some
components of the solution are not in $\BL^p$ for any $p>0$ (see
\cite[Example 7.3]{Kl4}).

We begin with BSDEs with one reflecting barrier.
\begin{tw}
\label{th.main5} Let $p\ge 1$. Assume that \mbox{\rm(H1)--(H6)}
are satisfied and  $L\in\mathbb{L}^{p}(\FF)$ in case $p>1$ and
that \mbox{\rm(H1)--(H6), (Z)} are satisfied and $L$ is of class
\mbox{\rm(D)} if $p=1$. Then there exists a solution $(Y,Z,K)$ of
RBSDE$(\xi,f+dV,L)$ such that $(Y,Z,K)\in\mathcal{D}^{p}\otimes
M\otimes \mathcal{V}^{+}$ if $p>1$ and
$(Y,Z,K)\in\mathcal{D}^{q}\otimes M\otimes \mathcal{V}^{+}$, $q\in
(0,1)$, $Y$ is of class \mbox{\rm(D)} if $p=1$, and all the
statements of assertion (ii) of Theorem \ref{th.main2} hold true.
\end{tw}
\begin{dow}
In the proof the basic text relates to the case where $p>1$; the
statements in parentheses relate to the case $p=1$. Let
$(\tilde{Y}^{n},\tilde{Z}^{n})\in \DM^{p}\otimes M^{p}$ (resp.
$(\tilde{Y}^{n},\tilde{Z}^{n})\in \DM^{q}\otimes M^{q},\, q\in
(0,1)$, such that $\tilde{Y}^{n}$ of class (D)) be a solution of
the  BSDE
\[
\tilde{Y}^{n}_{t}=\xi
+\intt f^{+}(s,\tilde{Y}^{n}_{s},\tilde{Z}^{n}_{s})\,ds
+\intt n(\tilde{Y}^{n}-L_{s})^{-}\,ds
+\intt dV_{s}-\intt \tilde{Z}^{n}_{s}\,dB_{s},\quad t\in [0,T].
\]
By Proposition \ref{prop3.4} (resp. Corollary \ref{cor3.10}),
$\tilde{Y}^{n}_{t}\le \tilde{Y}^{n+1}_{t}$ and $Y^{n}_{t}\le
\tilde{Y}^{n}_{t}$, $t\in [0,T]$, for $n\in\mathbb{N}$.
Consequently,
\[
Y^{1}_{t}\le Y^{n}_{t}\le \tilde{Y}^{n}_{t},\quad t\in [0,T],
\quad n\in\mathbb{N}.
\]
Let us observe that the data $(\xi,f^{+},L)$ satisfy assumptions
(H1)--(H7) (resp. (H1)-(H6), (H7*)) with $X=R(L)$, where $R(L)$ is
a c\`adl\`ag version of Snell's envelope of the process $L$. By
Theorem \ref{th.main1} (resp. Theorem \ref{th.main2}),
$\tilde{Y}^{n}_{t}\nearrow \tilde{Y}_{t},\, t\in [0,T]$, where
$\tilde{Y}\in\DM^{p}$ (resp. $\tilde{Y}\in\DM^{q},\, q\in (0,1)$,
$\tilde{Y}$ is of class (D)) is the first component of the
solution of RBSDE$(\xi,f^{+}+dV,L)$. Hence
\begin{equation}
\label{TT2}
Y^{1}_{t}\le Y^{n}_{t}\le\tilde{Y}_{t},
\quad t\in [0,T],\quad n\in\mathbb{N}.
\end{equation}
By Lemma \ref{lm4.55} there exist a stationary sequence
$\{\delta^{1}_{k}\}\subset \mathcal{T}$ and constants $c_{k}$ such
that
\[
|V|_{\delta^{1}_{k}}+Y^{1,*}_{\delta^{1}_{k}}
+\tilde{Y}_{\delta^{1}_{k}}\le c_{k},
\quad k\in\mathbb{N}.
\]
Put $\tau_{k}=\delta^{1}_{k}\wedge \delta^{2}_{k}$, where
\[
\delta^{2}_{k}=\inf\{t\in [0,T]; X^{*}_{t}
+\int_{0}^{t}f^{-}(s,X_{s},0)\,ds>k\}\wedge T.
\]
Since for every $p>1$ the data $(Y^{n}_{\tau_{k}},f,L)$ satisfy
the assumptions of Theorem \ref{th.main1} on each interval
$[0,\tau_{k}]$, the theorem follows.
\end{dow}

\begin{uw}
If $p>1$ then by Theorem \ref{th4.1} the solution of Theorem
\ref{th.main5} is unique in the class $\mathcal{D}^{p}\otimes
M\otimes \mathcal{V}^{+}$. We do not know whether in general the
solution is unique in case $p=1$. However, if $p=1$, then by
Remark \ref{uwnii}, the solution is unique in its class if $f$
does not depend on $z$.
\end{uw}

Each of the following conditions is called the Mokobodzki
condition.

\begin{enumerate}
\item[(M)] There exists $X\in\HH^{p}$ such that
$L_{t}\le X_{t}\le U_{t}$ for a.e. $t\in [0,T]$.
\item[(M*)]There exists
$X\in\mathcal{V}^{1}+\mathcal{M}^{q}_{c}$, $q\in (0,1)$, such that
$X$ is of class (D) and $L_{t}\le X_{t}\le U_{t}$ for a.e. $t\in
[0,T]$.
\end{enumerate}

\begin{tw}
\label{th.main6} Let $p>1$ and let assumptions
\mbox{\rm(H1)--(H6)} and \mbox{\rm(M)} hold.
\begin{enumerate}
\item[\rm(i)]
There exists a solution $(Y,Z,R)$ of RBSDE$(\xi,f+dV,L,U)$ such
that $(Y,Z,R)\in\DM^{p}\otimes M\otimes \mathcal{V}$.
\item[\rm(ii)]Let $(Y^{n,m},Z^{n,m})\in\DM^{p}\otimes M^{p}$ be
a solution of \mbox{\rm(\ref{m3.0}}) and let the triples
$(\bar{Y}^{n},\bar{Z}^{n},\bar{A}^{n})$,
$(\underline{Y}^{m},\underline{Z}^{m},
\underline{K}^{m})\in\DM^{p}\otimes M\otimes\mathcal{V}^{+}$ be
solutions of $\overline{R}$BSDE $(\xi,\bar{f}_{n}+dV,U)$
and $\underline{R}$BSDE$(\xi,
\underline{f}_{m}+dV,L)$ of Theorem \ref{th.main4}. Then
assertions (ii)--(vi) of Theorem \ref{th.main4} apart from
\mbox{\rm(\ref{eq6.29}), (\ref{eq6.30})} hold true.
\end{enumerate}
\end{tw}
\begin{dow}
Existence of solutions $(Y^{n,m},Z^{n,m})$ and
$(\bar{Y}^{n},\bar{Z}^{n},\bar{A}^{n})$,
$(\underline{Y}^{m},\underline{Z}^{m},\underline{K}^{m})$ follow
from Theorem \ref{th.main5} and Theorem \ref{th3.1}. Let
$\bar{y},\underline{y}$ be the first components of the solutions
of $\overline{\mbox{R}}$BSDE$(\xi,f+dV,U)$ and
$\underline{\mbox{R}}$BSDE$(\xi,f+dV,L)$, respectively. By Theorem
\ref{th.main5}, these solutions exist, are unique and
$\bar{y},\underline{y}\in\DM^{p}$. By Proposition \ref{prop3.4},
\[
\bar{y}_{t}\le Y^{n,m}_{t}\le \underline{y}_{t},
\quad t\in [0,T],\quad n,m\in\mathbb{N}.
\]
By Lemma \ref{lm4.55} there exist a stationary sequence
$\{\delta^{1}_{k}\}\subset \mathcal{T}$ and constants $c_{k}$ such
that
\[
|V|_{\delta^{1}_{k}}+|C|_{\delta^{1}_{k}}+\bar{y}^{*}_{\delta^{1}_{k}}
+\underline{y}^{*}_{\delta^{1}_{k}}\le c_{k},\quad k\in\BN,
\]
where $C$ is the finite variation part of the Doob-Meyer
decomposition of the process $X$. Put
$\tau_{k}=\delta^{1}_{k}\wedge\delta^{2}_{k}$, where
\[
\delta^{2}_{k}=\inf\{t\in [0,T];
\intot|f(s,X_{s},0)|\,ds>k\}\wedge T.
\]
Since for every $p>1$ the data
$(Y^{n,n}_{\tau_{k}},\bar{Y}^{n}_{\tau_{k}},
\underline{Y}^{m}_{\tau_{k}},f,L,U)$ satisfy the assumptions of
Theorem \ref{th.main4} on each interval $[0,\tau_{k}]$, applying
Theorem \ref{th.main3} gives the desired results.
\end{dow}
\medskip

Investigation of BSDEs with two reflecting barriers in case $p=1$
is more complicated than in case $p>1$, because if $p=1$ then in
general we cannot use Corollary \ref{cor.how} to compare
solutions, and in consequence we not know whether the sequences
$\{\bar{Y}^{n}\}, \{\underline{Y}^{m}\}$ of Theorem \ref{th.main6}
are monotone. To apply Corollary \ref{cor.how} one have to know
that $\bar{Z}^{n},\underline{Z}^{m}\in\mathbb{L}^{q}$ for some
$q>\alpha$, where $\alpha$ comes from condition (Z).
Unfortunately, in general $\bar{Z}^{n}, \underline{Z}^{m}$ are not
in $\mathbb{L}^{q}$ unless (H7*) is satisfied. If $f$ does not
depend on $z$ then by Remark \ref{uwnii}, to compare elements of
the sequences $\{\bar{Y}^{n}\}, \{\underline{Y}^{m}\}$ it suffices
to know that $\bar{Y}^{n},\underline{Y}^{m}\in\mathbb{L}^{q}(\FF)$
for some $q>\alpha$ and they are of class (D), which is the case.
If $f$ depends on $z$, to overcome the monotonicity difficulties
we replace solutions $\bar{Y}^{n},\underline{Y}^{m}$ of reflected
BSDEs by limits of their penalizations.

To be more specific, let $(Y^{n,m},Z^{n,m})$ denote a solution of
(\ref{m3.0})  and let $A^{n,m}$
be defined by (\ref{eq6.33}). By Theorem \ref{th.main5},
$\{(Y^{n,m},Z^{n,m},A^{n,m})\}_m$ converges, in the sense of
Theorem \ref{th.main2},
%converges in the sense specified by (in?) Theorem ..
to some process
$(\bar{Y}^{n},\bar{Z}^{n},\bar{A}^{n})\in\DM^{q}\otimes
M\otimes\mathcal{V}$, $q\in (0,1)$, such that $\bar{Y}^{n}$ is of
class \mbox{\rm(D)}. Similarly, if we denote by
$(Y^{n,m},Z^{n,m})$ a solution of (\ref{m3.0}) and we define $K^{n,m}$ by (\ref{eq6.33})
then $\{(Y^{n,m},Z^{n,m},A^{n,m})\}_n$ converges, in the sense of
Theorem \ref{th.main2}, to some process
$(\underline{Y}^{m},\underline{Z}^{m},
\underline{K}^{m})\in\DM^{q}\otimes M\otimes\mathcal{V}$, $q\in
(0,1)$, such that $\underline{Y}^{m}$ is of class \mbox{\rm(D)}.
As we shall see in the proof of the following theorem the
sequences $\{\bar{Y}^{n}\}_n$, $\{\underline{Y}^{m}\}_m$ defined
this way are monotone.

\begin{tw}
Let $p=1$ and let assumptions \mbox{\rm(H1)--(H6), (Z)} and
\mbox{\rm(M)} hold.
\begin{enumerate}
\item[\rm(i)]
There exists a solution $(Y,Z,R)$ of RBSDE$(\xi,f+dV,L,U)$ such
that $(Y,Z,R)\in\DM^{q}\otimes M\otimes\mathcal{V}$ for $q\in
(0,1)$ and $Y$ is of class \mbox{\rm(D)}.
\item[\rm(ii)]Let $(\bar{Y}^{n},\bar{Z}^{n},\bar{A}^{n}),
(\underline{Y}^{m},\underline{Z}^{m},\underline{K}^{m})$ denote
processes defined in the paragraph preceding the theorem. Then
assertions (ii)--(vi) of Theorem \ref{th.main4} apart from
\mbox{\rm(\ref{eq6.29}), (\ref{eq6.30})} hold true.
\end{enumerate}
\end{tw}
\begin{dow}
The proof runs as the proof Theorem \ref{th.main6} apart from the
fact that we use Theorem \ref{th3.2} instead of Theorem
\ref{th3.1}  and Corollary \ref{cor3.10} instead of Proposition
\ref{prop3.4}, and now we consider  $\bar{y},\underline{y}\in
\DM^{q}$, $q\in(0,1)$, such that  $\bar{y},\underline{y}$ are of
class (D) and are limits of penalizations for
$\overline{\mbox{R}}$BSDE$(\xi,f+dV,U)$ and
$\underline{\mbox{R}}$BSDE$(\xi,f+dV,L)$, respectively. The only
additional fact we have to prove is the monotonicity of the
sequences $\{\bar{Y}^{n}\}, \{\underline{Y}^{m}\}$. But the
monotonicity follows immediately from the fact that
%by the construction?
\[
\bar{Y}^{n}_{t}=\lim_{m\rightarrow +\infty} Y^{n,m}_{t},\quad
\underline{Y}^{m}_{t}=\lim_{n\rightarrow+\infty} Y^{n,m}_{t},\quad
t\in [0,T],
\]
and by Corollary \ref{cor3.10}, $Y^{n,m}_{t}\le Y^{n+1,m}_{t}$,
$Y^{n,m}_{t}\ge Y^{n,m+1}_{t}$, $t\in[0,T]$, $n\in\mathbb{N}$.
\end{dow}

\end{document}